\input amstex
\documentstyle{amsppt}
\loadmsbm

\nologo

\TagsOnRight

\NoBlackBoxes

\define\acc{\operatorname{acc}}

\define\supp{\operatorname{supp}}
\define\Lip{\operatorname{Lip}}
\define\diam{\operatorname{diam}}
\define\dist{\operatorname{dist}}

\def\floor{\mathbin{\hbox{\vrule height1.2ex width0.8pt depth0pt
        \kern-0.8pt \vrule height0.8pt width1.2ex depth0pt}}}

\font\letter=cmss10 

\font\normalsmall=cmss10 scaled 500

\font\normal=cmss10 scaled 700

\font\normalbig=cmss10

\define\ini{\operatorname{i}}

\define\termi{\operatorname{t}}

\define\smallsmallG{\text{\normalsmall G}}
\define\smallsmallvertexi{\text{\normalsmall i}}
\define\smallsmallvertexj{\text{\normalsmall j}}

\define\smalledge{\text{\normal e}}
\define\smallE{\text{\normal E}}
\define\smallV{\text{\normal V}}
\define\smallG{\text{\normal G}}
\define\smallvertexi{\text{\normal i}}
\define\smallvertexj{\text{\normal j}}

\define\edge{\text{\normalbig e}}

\define\E{\text{\normalbig E}}
\define\V{\text{\normalbig V}}
\define\G{\text{\normalbig G}}
\define\vertexi{\text{\normalbig i}}
\define\vertexj{\text{\normalbig j}}

\define\Haus{\text{\normal H}}

\define\scon{\text{\normal con}}

\define\sdyncon{\text{\normal dyn-con}}

\define\dyn{\text{\normal dyn}}

\define\vector{\text{\normal vec}}

\define\Hol{\text{\normal H\"ol}}

\define\erg{\text{\normal erg}}

\define\rel{\text{\normal rel}}

\define\specradsmall{\text{\normal spec-rad}}

\define\radius{\text{\normal rad}}

\define\distance{\text{\letter d}}

\define\LDistance{\text{\letter L}}

\font\tenscr=callig15 scaled 800
\font\sevenscr=callig15
\font\fivescr=callig15
\skewchar\tenscr='177 \skewchar\sevenscr='177 \skewchar\fivescr='177
\newfam\scrfam \textfont\scrfam=\tenscr \scriptfont\scrfam=\sevenscr
\scriptscriptfont\scrfam=\fivescr
\def\scr#1{{\fam\scrfam#1}}

\font\tenscri=suet14
\font\sevenscri=suet14
\font\fivescri=suet14
\skewchar\tenscri='177 \skewchar\sevenscri='177 \skewchar\fivescri='177
\newfam\scrifam \textfont\scrifam=\tenscri \scriptfont\scrifam=\sevenscri
\scriptscriptfont\scrifam=\fivescri
\def\scri#1{{\fam\scrifam#1}}

%
%
%
%

\hsize = 6.07 true in
\vsize = 9.23 true in

\topmatter
\title
Dynamical multifractal zeta-functions
and
fine
multifractal spectra
of graph-directed self-conformal constructions
\endtitle
\endtopmatter

\centerline{\smc V\. Mijovi\'c}
\centerline{Department of Mathematics}
\centerline{University of St\. Andrews}
\centerline{St\. Andrews, Fife KY16 9SS, Scotland}
\centerline{e-mail: {\tt vm27\@st-and.ac.uk}}

\medskip

\centerline{\smc L\. Olsen}
\centerline{Department of Mathematics}
\centerline{University of St\. Andrews}
\centerline{St\. Andrews, Fife KY16 9SS, Scotland}
\centerline{e-mail: {\tt lo\@st-and.ac.uk}}

\topmatter
\abstract{We introduce 
multifractal pressure
and
dynamical 
multifractal zeta-functions
 providing precise information 
of  a very general class of multifractal spectra, including, for example, the 
fine
multifractal spectra of 
graph-directed self-conformal 
 measures 
 and the fine multifractal spectra of ergodic Birkhoff averages of continuous functions
 on graph-directed self-conformal sets.
}
\endabstract
\endtopmatter

\bigskip

\centerline{\smc Contents}
\medskip
\roster
 \item"1." Introduction.
 \item"2." The setting, Part 1: 
Graph-directed self-conformal sets 
and 
graph-directed self-conformal measures.
 \item"3." The setting, Part 2:
Pressure and dynamical zeta-functions.
 \item"4." Statements of the  main results.
 \item"5." Applications:
 Multifractal spectra of measures
 and
  multifractal spectra of ergodic Birkhoff averages.
 \item"6." 
 Proofs. The map $M_{n}$.
   \item"7."
   Proofs. The measures $\Pi$ and $\Pi_{n}$.
   \item"8."
   Proofs. The modified multifractal pressure.
     \item"9."
   Proof of Theorem 4.4. 
      \item"10."
   Proof of Theorem 4.6.  
 \item"" References
\endroster

\medskip

\footnote""
{
\!\!\!\!\!\!\!\!
2000 {\it Mathematics Subject Classification.} 
Primary: 28A78.
Secondary: 37D30, 37A45.\newline
{\it Key words and phrases:} 
multifractals,
zeta functions.
pressure,
Bowen's formula,
large deviations,
Hausdorff dimension,
graph-directed self-conformal sets
}

\leftheadtext{L\. Olsen}

\rightheadtext{Dynamical multifractal zeta-functions}

\heading{1. Introduction.}\endheading

For a Borel measure $\mu$ on $\Bbb R^{d}$ 
and
a positive number $\alpha$,
let us consider
the set  of
those points
$x$ in $\Bbb R^{d}$ for which the measure
$\mu(B(x,r))$ of the ball
$B(x,r)$ with center $x$ and radius $r$ behaves like
$r^{\alpha}$ for small $r$,
i\.e\. the set
 $$
   \Bigg\{
   x\in K
   \,\Bigg|\,
   \lim_{r\searrow 0}
   \frac{\log\mu(B(x,r))}{\log r}
   =
   \alpha
   \Bigg\}\,.
   \tag1.1
 $$
If the intensity of the measure $\mu$ varies very widely, it may
happen that the sets in (1.1)
display a
fractal-like character for a range of values of $\alpha$. In this case
it is natural to study
the Hausdorff dimensions of the sets in (1.1)
 as $\alpha$
varies.
We therefore define the  fine multifractal spectrum of $\mu$
by
 $$
 f_{\mu}(\alpha)
 =
 \dim_{\Haus}
   \Bigg\{
   x\in K
   \,\Bigg|\,
   \lim_{r\searrow 0}
   \frac{\log\mu(B(x,r))}{\log r}
   =
   \alpha
   \Bigg\}\,.
 \tag1.2
$$
where $\dim_{\Haus}$ denotes the Hausdorff dimension;
here and below
we use the following convention,
namely, we define the Hausdorff
of the empty set to be $-\infty$, i\.e\. we put
 $$
 \dim_{\Haus}\varnothing
 =
 -
 \infty\,.
 $$
The fine multifractal spectrum is one of the 
the two main ingredients
in multifractal analysis. 
The second main ingredient is the  
Renyi dimensions.
Renyi
dimensions quantify the varying intensity of a measure by analyzing its moments at different scales. 
Formally, 
for $q\in\Bbb R$,
the $q$'th
Renyi dimensions 
$\tau_{\mu}(q)$
of $\mu$ is defined by
 $$
 \tau_{\mu}(q)
 =
 \lim_{r\searrow 0}
 \frac{\dsize \log\int\limits_{\,\,\,K}\mu(B(x,r))^{q-1}\,d\mu(x)}{-\log r}\,,
 $$
provided the limit exists. 
One of the main problems
in multifractal analysis is to 
 understand the
 multifractal spectrum and the Renyi dimensions,
 and
their relationship with each other.
During the past 20 years
there has been an enormous interest 
in
computing the multifractal spectra of measures
in 
the mathematical literature
and within the last
15 years the multifractal spectra of various classes of measures
in Euclidean space $\Bbb R^{d}$ 
exhibiting some degree of self-similarity have been computed 
rigorously, 
see
the textbooks [Fa,Pe]
and the references therein.

Dynamical zeta-functions were introduced by Artin \& Mazur
in the mid 1960's [ArMa]
based on an analogy with the number theoretical zeta-functions
associated with a function field
over a finite ring.
Subsequently Ruelle [Rue1,Rue2]
associated
zeta-functions
to 
certain statistical mechanical models in one dimensions.
During the past 35 years 
many parallels have been 
drawn between
number theory 
zeta-functions, 
dynamical zeta-functions,
and statistical mechanics zeta-functions.
However, much more recently and
motivated by the 
powerful techniques
provided by 
the use of the Artin-Mazur zeta-functions in number theory
and the use of the Ruelle
zeta-functions in dynamical systems,
Lapidus and collaborators
(see the intriguing books by Lapidus \& van Frankenhuysen [Lap-vF1,Lap-VF2]
and the references therein)
have recently
introduced and pioneered
to use of
zeta-functions in fractal geometry.
Inspired by this development,
within the past 4--5 years
several authors have 
paralleled this development
by
introducing
zeta-functions
into  multifractal geometry.
Indeed, in 2009,
Lapidus and collaborators introduced 
various intriguing 
{\it geometric} multifractal zeta-functions [LapRo,LapLe-VeRo]
designed to provide information about the multifractal spectrum
$f_{\mu}(\alpha)$
of 
self-similar measures $\mu$,
and
many
connections with multifractal spectra were
suggested
and in some cases proved;
for example, in simplified cases 
the multifractal spectrum of a self-similar measure could be recovered from a 
zeta-function.
The 
key idea in [LapRo,LapLe-VeRo] is both simple and attractive:
while traditional zeta-functions are defined  by 
\lq\lq summing over all data",
the
multifractal zeta-functions
in
[LapRo,LapLe-VeRo] 
are defined by only
\lq\lq summing over data that are multifractally relevant".
This idea is also the
leitmotif
in this work (as well as in earlier  work 
[Bak,MiOl,Ol4]), see,
in particular,
the first remark following the definition of the zeta-function
$\zeta_{C}^{\dyn,U}(\varphi;\cdot)$  in
Section 4.
Ideas similar to those in
[LapRo,LapLe-VeRo]
have very recently
been revisited and
investigated 
in
[Bak,MiOl] where
the authors
introduce related
{\it geometric}
multifractal zeta-functions 
tailored to study
the
multifractal spectra
of self-conformal measures
and a number of connections with
very general types of multifractal spectra were established.

We also point out that
the work by Lapidus et al\.  [LapRo,LapLe-VeRo] was followed 
shortly
afterwards
by
the introduction 
of 
a different type of multifractal zeta-function
by Levy-Vehel \& Mendivil
[Le-VeMe]
tailored to
provide information about the 
Renyi dimensions 
$\tau_{\mu}(q)$
of self-similar measures $\mu$.
Ideas related to those in [Le-VeMe]
have also been investigated
in [Ol2,Ol3]
where
the author
introduce
multifractal zeta-functions 
designed
 to study
the 
Renyi dimensions
and the (closely related) multifractal  Minkowski volume
of self-conformal measures.


In addition to the distinctively
{\it geometric}
approaches
in 
[Bak,Le-VeMe,MiOl,Ol2,Ol3],
it has been a major challenge 
to 
introduce 
and develop 
a
natural and meaningful
theory of
{\it dynamical}
multifractal zeta-functions
paralleling
the existing
 powerful 
theory
 of 
{\it dynamical} zeta-functions introduced and developed by Ruelle [Rue1,Rue2]
and others, see, for example,  
the surveys and books 
[Bal1,Bal2,ParPo1,ParPo2]
and the references therein.
In particular,
in the setting of self-conformal constructions,
[Ol4]
 introduced
a family
of 
{\it dynamical}
multifractal zeta-functions
designed
 to
 provide
precise information of very general classes of multifractal spectra, 
including, for example, the multifractal 
spectra of self-conformal measures and the multifractal spectra of ergodic 
Birkhoff averages 
of continuous functions.

However,
recently it has been recognised that
while self-conformal constructions
provide a useful and  important framework
 for studying fractal and multifractal geometry,
the
more general notion of 
{\it graph-directed} self-conformal constructions
provide 
a substantially more flexible 
and useful framework, see, for example, [MaUr]
for an elaboration of this.
In recognition of this viewpoint,
the purpose 
of this paper is to 
developed
a {\it dynamical}
theory of
multifractal zeta-functions in the setting of
{\it graph-directed} self-conformal constructions.


In Section 2--3 we briefly recall the 
definitions
of self-conformal constructions
and the accompanying
pressure and 
dynamical zeta-functions.
In Section 4
we provide our main definitions
of the
multifractal
pressure and
the
multifractal 
dynamical zeta-functions
 and we 
 state our main 
  results.
 Section 5 contains a number of examples, including,
 multifractal spectra
 of graph-directed
 self-conformal measures
 and different types of 
 multifractal spectra of ergodic averages of continuous functions on
 graph-directed
 self-conformal sets.
Finally, the proofs are presented in Sections 6--10.

\bigskip

%
%
%

\centerline{\smc 2. The setting, Part 1:}

\centerline{\smc
Graph-directed
self-conformal sets and 
graph-directed
self-conformal measures.}

\medskip

{\bf 2.1. Notation from symbolic dynamics.}
We first recall the 
 notation
and terminology  from symbolic dynamics 
that will be used in this paper.
Fix
a finite directed multigraph
$\G=(\V,\E)$
where
$\V$ denotes the set of vertices of $\G$ and
$\E$ denotes the set of edges of $\G$.
We will always assume that the graph $\G$ is strongly connected.
For an edge $\edge\in\E$,
we 
write
$\ini(\edge)$
for the initial vertex of $\edge$ 
and we 
write
$\termi(\edge)$ for the terminal vertex of $\edge$.
For $\vertexi,\vertexj\in\V$, write
 $$
 \aligned
 \E_{\smallvertexi}
&=
\Big\{\edge\in\E\,\Big|\,\ini(\edge)=\vertexi\Big\}\,,\\
 \E_{\smallvertexi,\smallvertexj}
&=
\Big\{\edge\in\E\,\Big|\,\ini(\edge)=\vertexi\,\,\text{and}\,\,\termi(\edge)=\vertexj\Big\}\,;
\endaligned
\tag2.1
 $$  
i\.e\.
 $\E_{\smallvertexi}$ is the family of all edges starting at $\vertexi$;
 and
 $ \E_{\smallvertexi,\smallvertexj}$
 is the 
 is the family of all edges starting at $\vertexi$
 and ending at $\vertexj$.
Also, for a positive integer $n$,
we
 write
 $$
 \aligned
 \Sigma_{\smallG}^{n}
&=
\Big\{
\edge_{1}\ldots\edge_{n}
\,\Big|\,
\edge_{i}\in\E\,\,\text{for $1\le i\le n$,}\\
&\qquad\qquad
   \quad\,\,\,\,\,\,
\termi(\edge_{1})=\ini(e_{2})\,,\\
&\qquad\qquad
   \quad\,\,\,\,\,\,
\termi(e_{i-1})=\ini(\edge_{i})
\,\,\text{and}\,\,
\termi(\edge_{i})=\ini(e_{i+1})
\,\,\text{for $1< i < n$,}\\
&\qquad\qquad
   \quad\,\,\,\,\,\,
\termi(e_{n-1})=\ini(\edge_{n})
\Big\}\\
 \Sigma_{\smallG}^{*}
&=
\bigcup_{m}\,\Sigma_{\smallG}^{m}\,,\\
 \Sigma_{\smallG}^{\Bbb N}
&=
\Big\{
\edge_{1}\edge_{2}\ldots\,
\,\Big|\,
\edge_{i}\in\E\,\,\text{for $1\le i$,}\\
&\qquad\qquad
   \quad\,\,\,\,\,\,
\termi(\edge_{1})=\ini(e_{2})\,,\\
&\qquad\qquad
   \quad\,\,\,\,\,\,
\termi(e_{i-1})=\ini(\edge_{i})
\,\,\text{and}\,\,
\termi(\edge_{i})=\ini(e_{i+1})
\,\,\text{for $1< i $}
\Big\}\,;
 \endaligned
 \tag2.2
 $$
i\.e\. $\Sigma_{\smallG}^{n}$ is the family of all
finite strings
$\bold i=\edge_{1}\ldots \edge_{n}$
consisting of finite paths in $\G$
of length $n$;
$\Sigma_{\smallG}^{*}$ is the family of all finite strings
$\bold i=\edge_{1}\ldots \edge_{m}$
with $m\in\Bbb N$ consisting of finite paths in $\G$;
and 
$\Sigma_{\smallG}^{\Bbb N}$ is the family of all
infinite
strings
$\bold i=\edge_{1}\edge_{2}\ldots $
consisting of infinite paths in $\G$.
For 
a finite string 
$\bold i=\edge_{1}\ldots \edge_{n}\in\Sigma_{\smallG}^{n} $, we write
 $$
 \ini(\bold i)
 =
 \ini(\edge_{1})\,,\,\,\,\,
 \termi(\bold i)
 =
 \termi(\edge_{n})\,,
 \tag2.3
 $$
and for an infinite string
$\bold i=\edge_{1} \edge_{2}\ldots\in\Sigma_{\smallG}^{\Bbb N} $, we write
 $$
 \ini(\bold i)
 =
 \ini(\edge_{1})\,.
 \tag2.4
 $$
Next,
for an infinite string 
$\bold i=\edge_{1}\edge_{2}\ldots\in\Sigma_{\smallG}^{\Bbb N} $
and a positive integer $n$, we will write
$\bold i|n
=
\edge_{1}\ldots \edge_{n}$.
In addition, for
a positive integer $n$
and
a finite string 
$\bold i=\edge_{1}\ldots \edge_{n}\in\Sigma_{\smallG}^{n} $
with length equal to  $n$,
 we will write
$|\bold i|
=
n$, and we let $[\bold i]$ denote the cylinder 
generated by $\bold i$, i\.e\.
  $$
  [\bold i]
 =
 \Big\{
 \bold j\in\Sigma_{\smallG}^{\Bbb N}
 \,\Big|\,
 \bold j|n=\bold i
 \Big\}\,.
 \tag2.5
 $$
Finally, let $S:\Sigma_{\smallG}^{\Bbb N}\to\Sigma_{\smallG}^{\Bbb N}$ denote the shift map, i\.e\.
 $$
 S(\edge_{1}\edge_{2}\ldots)
 =
 \edge_{2}\edge_{3}\ldots\,.
 $$

\bigskip

{\bf 2.2. Graph-directed
self-conformal sets and 
graph-directed 
self-conformal measures.}
Next, we recall the definition of
graph-directed
 self-conformal 
sets and measures. 
A graph-directed conformal iterated function system with probabilities 
is a list
 $$
 \Big(\,
 \V,\,
 \E,\,
 (V_{\smallvertexi})_{\smallvertexi\in\smallV},\,
 (X_{\smallvertexi})_{\smallvertexi\in\smallV},\,
 (S_{\smalledge})_{\smalledge\in\smallE},\,
 (p_{\smalledge})_{\smalledge\in\smallE}\,
 \Big)
 $$
where

\bigskip

\roster
\item"$\bullet$"
For each $\vertexi\in\V$
we have:
$V_{\smallvertexi}$ is an open, connected subset of $\Bbb R^{d}$.
\item"$\bullet$" 
For each $\vertexi\in\V$
we have:
$X_{\smallvertexi}\subseteq V_{\smallvertexi}$ is a compact set with
$X_{\smallvertexi}^{\circ\,-}=X_{\smallvertexi}$.
\item"$\bullet$" 
For each $\vertexi,\vertexj\in\V$
and
$\edge\in E_{\smallvertexi,\smallvertexj}$
we have:
$S_{\smalledge}:V_{\smallvertexj}\to V_{\smallvertexi}$
is a contractive
$C^{1+\gamma}$ diffeomorphism with
$0<\gamma<1$
such that
$S_{\smalledge}(X_{\smallvertexj})\subseteq X_{\smallvertexi}$.
\item"$\bullet$"
The Conformality Condition.
For $\vertexi,\vertexj\in\V$,
$\edge\in E_{\smallvertexi,\smallvertexj}$ and $x\in V_{\smallvertexj}$, let
$$
(DS_{\smalledge})(x):\Bbb R^{d}\to\Bbb R^{d}
$$
denote the derivative of $S_{\smalledge}$ at $x$. 
For each $\vertexi,\vertexj\in\V$
and
$\edge\in E_{\smallvertexi,\smallvertexj}$
we have:
$(DS_{\smalledge})(x)$ is a contractive similarity map, i\.e\.
there exists
$s_{\smalledge}(x)\in(0,1)$ such that
$|(DS_{\smalledge})(x)u-(DS_{\smalledge})(x)v|
 =
 s_{\smalledge}(x)|u-v|$
for all $u,v\in\Bbb R^{d}$.
\item"$\bullet$" 
For each $\vertexi\in\V$
we have:
$(p_{\smalledge})_{\smalledge\in\smallE_{\smallsmallvertexi}}$ is a probability vector.
\endroster

\bigskip

\noindent
It follows from [Hu] that there exists a unique
list $(K_{\smallvertexi})_{\smallvertexi\in\smallV}$ of
non-empty compact sets $K_{\smallvertexi}\subseteq X_{\smallvertexi}$ such that
 $$
 K_{\smallvertexi}
 \,=\,
 \bigcup_{\smalledge\in \smallE_{\smallsmallvertexi}}\,
 S_{\smalledge}K_{\termi(\smalledge)}\,,
 \tag2.6
 $$
and a unique
list $(\mu_{\smallvertexi})_{\smallvertexi\in\smallV}$ 
of probability measures with
$\supp\mu_{\smallvertexi}=K_{\smallvertexi}$ such that
 $$
 \mu_{\smallvertexi}
 \,=\,
 \sum_{\smalledge\in E_{\smallsmallvertexi}}\,
 p_{\smalledge}\,\mu_{\termi(\smalledge)}\circ S_{\smalledge}^{-1}\,.
 \tag2.7
 $$
The sets $(K_{\smallvertexi})_{\smallvertexi\in\smallV}$
and measures $(\mu_{\smallvertexi})_{\smallvertexi\in\smallV}$
are called the self-conformal sets and self-conformal
measures associated with the list 
$ (\,
 \V,\,
 \E,\,
 (V_{\smallvertexi})_{\smallvertexi\in\smallV},\,
 (X_{\smallvertexi})_{\smallvertexi\in\smallV},\,
 (S_{\smalledge})_{\smalledge\in\smallE},\,
 (p_{\smalledge})_{\smalledge\in\smallE}\,
)$, respectively.
We will frequently assume that the so-called Open Set condition (OSC) 
 is satisfied. The OSC is defined as follows:

\bigskip

\roster
\item"$\bullet$"
The Open Set Condition:
There exists a list
$(U_{\smallvertexi})_{\smallvertexi\in\smallV}$ of open non-empty and bounded sets
$U_{\smallvertexi}\subseteq X_{\smallvertexi}$ with 
$S_{\smalledge}\big(\,U_{\smallvertexj}\,\big)
 \subseteq
 U_{\smallvertexi}$
for all $\vertexi,\vertexj\in\V$ and all $\edge\in\E_{\smallvertexi,\smallvertexj}$ such that
 $S_{\smalledge_{1}}\big(\,U_{\termi(\smalledge_{1})}\,\big)
 \,\cap\,
 S_{\smalledge_{2}}\big(\,U_{\termi(\smalledge_{2})}\,\big)
 =
 \varnothing$
for all $\vertexi\in\V$ and all $\edge_{1},\edge_{2}\in\E_{\smallvertexi}$ with 
$\edge_{1}\not=\edge_{2}$.
\endroster 
 
\bigskip

\noindent
For $\bold i=\edge_{1}\ldots\edge_{n}\in\Sigma_{\smallG}^{n}$, we write
 $$
 \aligned
 S_{\bold i}
&=
 S_{\smalledge_{1}}\cdots S_{\smalledge_{n}}\,,\\
 K_{\bold i}
&=
 S_{\smalledge_{1}}\cdots S_{\smalledge_{n}}
 \big(K_{\termi(\smalledge_{n})}\big)\,,\\
p_{\bold i}
&=
 p_{\smalledge_{1}}\cdots p_{\smalledge_{n}}\,,\\ 
 \endaligned
 \tag2.8
 $$ 
and
we define the projection
$\pi:\Sigma_{\smallG}^{\Bbb N}\to \Bbb R^{d}$ by
 $$
 \big\{\,\pi(\bold i)\,\big\}
 =
 \bigcap_{n}
 \,
 K_{\bold i|n}
 \tag2.9
 $$
for $\bold i=\edge_{1}\edge_{2}\ldots\in\Sigma_{\smallG}^{\Bbb N}$.
Finally, we define
$\Lambda:\Sigma_{\smallG}^{\Bbb N}\to\Bbb R$ by
 $$
 \align
 \Lambda(\bold i)
 \,
&=\,
 \log
 \big|
 \,
 \big(DS_{\smalledge_{1}}\big)\big(\pi_{\termi(\smalledge_{1})}(S\bold i)\big)
 \,
 \big|
 \tag2.10
 \endalign
 $$
for $\bold i=\edge_{1}\edge_{2}\ldots\in\Sigma_{\smallG}^{\Bbb N}$;
loosely speaking the 
map $\Lambda$ represents the local change of scale as one goes from
$\pi_{\termi(\smalledge_{1})}(S\bold i)$ to
$\pi_{\ini(\smalledge_{1})}(\bold i)$.

\bigskip


\centerline{\smc 3. The setting, Part 2:}

\centerline{\smc
Pressure and dynamical zeta-functions.}

\medskip

Throughout 
 this section, 
and in the remaining parts of the paper, we will used the
following notation.
 Namely,
if $(a_{n})_{n}$ is a sequence of 
complex numbers
and if $f$ is the power
series
defined by $f(z)=\sum_{n}a_{n}z^{n}$ for $z\in\Bbb C$,
then we will denote the 
radius
of convergence of $f$ by $\sigma_{\radius}(f)$, i\.e\. we 
write
 $$
 \sigma_{\radius}(f)
 =
 \text{
 \lq\lq the radius of convergence of $f$"
 }\,.
 $$

Our
definitions
and
 results are
 motivated by the 
 notion of pressure
 from the thermodynamic formalism
  and
 the
 dynamical zeta-functions introduced by
 Ruelle [Rue1,Rue2]; see, also 
 [Bal1,\allowlinebreak
 Bal2,\allowlinebreak
 ParPo1,\allowlinebreak
 ParPo2].
 In addition,
Bowen's formula 
expressing the 
Hausdorff dimension of 
a self-conformal set
in terms of the pressure 
(or the dynamical zeta-function)
of the scaling map $\Lambda$ in (2.5)
also plays a leitmotif 
in our work.
Because of this we
now
recall the definition 
of pressure and dynamical zeta-function,
and the statement of Bowen's formula.
Let
$\varphi:\Sigma_{\smallG}^{\Bbb N}\to\Bbb R$
be a
continuous function;
here, and below, we equip
$\Sigma_{\smallG}^{\Bbb N}$
with the product topology
with discrete factors
and
all statements 
about continuity involving
$\Sigma_{\smallG}^{\Bbb N}$
will always refer to this topology.
 The
pressure
of 
 $\varphi$
 is defined 
by
 $$
 \align
 P(\varphi)
 &=
\lim_{n}
  \,\,
  \frac{1}{n}
 \,\,
 \log
   \sum
  \Sb
\bold i\in\Sigma_{\smallsmallG}^{n}
  \endSb
  \,\,
  \sup_{\bold u\in[\bold i]}
  \,\,
   \exp
 \,\,
 \sum_{k=0}^{n-1}\varphi S^{k}\bold u\,,
 \tag3.1
  \endalign
 $$
 see
 [Bo2]
or
[ParPo2];
we note that
it is well-known that the limit in (3.1) exists.
Also,
the dynamical zeta-function of $\varphi$ is defined by
$$
  \zeta^{\dyn}(\varphi;z)
  =
   \sum_{n}
   \,\,
   \frac{z^{n}}{n}
 \left(
   \sum
  \Sb
\bold i\in\Sigma_{\smallsmallG}^{n}
  \endSb
  \,\,
  \sup_{\bold u\in[\bold i]}
  \,\,
 \exp
 \,\,
 \sum_{k=0}^{n-1}\varphi S^{k}\bold u
\right)
\tag3.2
 $$
for those complex numbers $z$ for which the series converge, see
[ParPo2].
We now  list two easily established and well-known
properties of 
the
pressure
$P(\varphi)$
and of 
the radius of convergence
$\sigma_{\radius}\big(\,\zeta^{\dyn}(\varphi;\cdot)\,\big)$
of
the power-series  
$\zeta^{\dyn}(\varphi;\cdot)$.
While both results are well-known and easily proved
(see, for example, [Bar,Fa]), 
we have decided to list them 
since 
they play an important part 
in 
the discussion of our results.

\bigskip

\proclaim{Theorem A (see, for example, [Bar,Fa]). Radius of convergence}
Fix a continuous function
$\varphi:\Sigma_{\smallG}^{\Bbb N}\to\Bbb R$.
Then we have
  $$
  -
  \log\sigma_{\radius}\big(\,\zeta^{\dyn}(\varphi;\cdot)\,\big)
  =
  P(\varphi)\,.
  $$
\endproclaim

\bigskip

\proclaim{Theorem B (see, for example, [Bar,Fa]). 
Continuity and monotony properties of the pressure}
Fix a  
continuous function
$\Phi:\Sigma_{\smallG}^{\Bbb N}\to\Bbb R$
with $\Phi<0$.
Then the function
 $
 t\to P(t\Phi),
 $
where $t\in\Bbb R$, 
is continuous, strictly decreasing and convex 
with
$\lim_{t\to-\infty}P(t\Phi)=\infty$
and
 $\lim_{t\to\infty}P(t\Phi)=-\infty$.
 In particular, there is a unique real number 
 $s$ such that
  $$
  P(s\Phi)=0\,;
  $$
 alternatively,
 $s$ 
  is the unique real number such that
  $$
  \sigma_{\radius}\big(\,\zeta^{\dyn}(s\Phi;\cdot)\,\big)
  =
  1\,.
  $$
 \endproclaim

\bigskip

\noindent
The main importance of the 
pressure
(for the purpose of this exposition)  is that 
it provides a beautiful formula for the 
Hausdorff dimension of
a
graph-directed
 self-conformal set satisfying the OSC.
This
result
was first noted by 
[Bo1]
(in the setting of quasi-circles)
and is the content of the next result.

\bigskip

\proclaim{Theorem C (see, for example, [Bar,Fa]). Bowen's formula}
Let
 $(K_{\smallvertexi})_{\smallvertexi\in\smallV}$ be the 
 graph-directed
 self-conformal sets defined by (2.6)
 and
 let $\Lambda:\Sigma_{\smallG}^{\Bbb N}\to\Bbb R$
 be the scaling function defined by (2.10).
 Let
$s$ be the unique real number
such that
  $$
  P(s\Lambda)=0\,;
  $$
 alternatively,
 $s$ 
  is the unique real number such that
  $$
  \sigma_{\radius}\big(\,\zeta^{\dyn}(s\Lambda;\cdot)\,\big)
  =
  1\,.
  $$  
If the OSC is satisfied, then we have
  $$
  \dim_{\Haus}K_{\smallvertexi}=s
  $$ 
 for all $\vertexi\in\V$. 
 \endproclaim

\bigskip

\noindent
It is reasonable to
expect that
any meaningful
theory of dynamical multifractal zeta-functions
should
produce
multifractal analogues of Bowen's equation.
In the next section
we 
will
develop such a theory
for graph-directed self-conformal constructions
(extending the 
theory for self-conformal constructions introduced
 in [Ol4]).

  \bigskip


\heading{4. Statements of the main results.}\endheading

We denote the family of Borel probability measures on 
$\Sigma_{\smallG}^{\Bbb N}$
and the 
 family of shift invariant Borel probability measures on 
$\Sigma_{\smallG}^{\Bbb N}$
by $\Cal P(\Sigma_{\smallG}^{\Bbb N})$
and
$\Cal P_{S}(\Sigma_{\smallG}^{\Bbb N})$, respectively, i\.e\.
we write
 $$
 \align
 \Cal P(\Sigma_{\smallG}^{\Bbb N})
&=
 \Big\{
 \mu\,\Big|\,
 \text{
 $\mu$ is a Borel probability measures on 
$\Sigma_{\smallG}^{\Bbb N}$
 }
 \Big\}\,,\\
 \Cal P_{S}(\Sigma_{\smallG}^{\Bbb N})
&=
 \Big\{
 \mu\,\Big|\,
 \text{
 $\mu$ is a shift invariant Borel probability measures on 
$\Sigma_{\smallG}^{\Bbb N}$
 }
 \Big\}\,;\\
 \endalign
 $$
we will always
equip
$\Cal P(\Sigma_{\smallG}^{\Bbb N})$ 
and
$\Cal P_{S}(\Sigma_{\smallG}^{\Bbb N})$ with the weak topologies.
We now fix a metric space $X$ and a
continuous map
$U:\Cal P(\Sigma_{\smallG}^{\Bbb N})\to X$.
The multifractal zeta-function
framework
developed 
 in this paper
depend on the space $X$ and the map $U$;
judicious
choices  of $X$ and $U$ 
will provide
 important 
examples, including,
multifractal spectra of 
graph-directed self-conformal measures
(see Section 5.1)
and
a variety of
multifractal spectra 
of
ergodic averages
of continuous
functions on graph-directed self-conformal sets
(see Section 5.2).
Next,
 for a positive integer $n$,
 let
$L_{n}:\Sigma_{\smallG}^{\Bbb N}\to\Cal P(\Sigma_{\smallG}^{\Bbb N})$ be defined 
by
 $$
 L_{n}\bold i
 =
 \frac{1}{n}\sum_{k=0}^{n-1}\delta_{S^{k}\bold i}\,;
 \tag4.1
 $$
 recall, that $S:\Sigma_{\smallG}^{\Bbb N}\to\Sigma_{\smallG}^{\Bbb N}$ denotes the shift map.
 We can now define the 
multifractal pressure and zeta-function
associated with the space $X$ and the map $U$.

 \bigskip
 
 \proclaim{Definition.
 The multifractal pressure 
 $ \underline P_{C}^{U}(\varphi)$
 and
 $ \overline P_{C}^{U}(\varphi)$
associated with the space $X$ and the map $U$}
Let $X$ be a metric space and let $U:\Cal P(\Sigma_{\smallG}^{\Bbb N})\to X$ be 
continuous with respect to the weak topology.
Fix a continuous function
$\varphi:\Sigma_{\smallG}^{\Bbb N}\to\Bbb R$.
    For $C\subseteq X$,
we define the lower and upper
mutifractal pressure
of 
 $\varphi$
associated with the space $X$ and the map $U$ by
 $$
 \aligned
  \underline P_{C}^{U}(\varphi)
 &=
 \,
  \liminf_{n}
 \,
  \,\,
  \frac{1}{n}
 \,\,
 \log
   \sum
  \Sb
  \bold i\in\Sigma_{\smallsmallG}^{n}\\
  {}\\
  UL_{n}[\bold i]\subseteq C
  \endSb
   \,\,
  \sup_{\bold u\in[\bold i]}
  \,\,
 \exp
 \,\,
 \sum_{k=0}^{n-1}\varphi S^{k}\bold u\,,\\
   \overline P_{C}^{U}(\varphi)
 &=
  \limsup_{n}
  \,\,
  \frac{1}{n}
 \,\,
 \log
   \sum
  \Sb
 \bold i\in\Sigma_{\smallsmallG}^{n}\\
  {}\\
  UL_{n}[\bold i]\subseteq C
  \endSb
 \,\,
  \sup_{\bold u\in[\bold i]}
  \,\,
 \exp
 \,\,
 \sum_{k=0}^{n-1}\varphi S^{k}\bold u\,.
 \endaligned
 \tag4.2
 $$
If $ \underline P_{C}^{U}(\varphi)$ and $ \overline P_{C}^{U}(\varphi)$ coincide, then we write
$P_{C}^{U}(\varphi)$ for their common value, i\.e\. we write
$ P_{C}^{U}(\varphi)
=
 \underline P_{C}^{U}(\varphi)
 =
  \overline P_{C}^{U}(\varphi)$.

\endproclaim

 \bigskip
 
 \proclaim{Definition.
 The dynamical multifractal zeta-function 
 $\zeta_{C}^{\dyn,U}(\varphi;\cdot)$
 associated with the space $X$ and the map $U$}
 Let $X$ be a metric space and let $U:\Cal P(\Sigma_{\smallG}^{\Bbb N})\to X$ be 
continuous with respect to the weak topology.
Fix a continuous function
$\varphi:\Sigma_{\smallG}^{\Bbb N}\to\Bbb R$.
    For $C\subseteq X$,
we define the 
dynamical
multifractal
 zeta-function $\zeta_{C}^{\dyn,U}(\varphi;\cdot)$
associated with the space $X$ and the map $U$ by
 $$
  \zeta_{C}^{\dyn,U}(\varphi;z)
  =
   \sum_{n}
   \,\,
   \frac{z^{n}}{n}
 \left(
   \sum
  \Sb
 \bold i\in\Sigma_{\smallsmallG}^{n}\\
  {}\\
  UL_{n}[\bold i]\subseteq C
  \endSb
  \,\,
  \sup_{\bold u\in[\bold i]}
  \,\,
 \exp
 \,\,
 \sum_{k=0}^{n-1}\varphi S^{k}\bold u
\right)
\tag4.3
 $$
for those complex numbers $z$ for which the series converges.

\endproclaim

\bigskip

\noindent
{\bf Remark.}
Comparing the
definition of the pressure (3.1)
(the dynamical zeta-function (3.2))
and
the
definition of the multifractal pressure (4.2)
(the dynamical multifractal zeta-function (4.3)),
it is clear that 
the 
multifractal pressure
(the dynamical multifractal zeta-function)
is obtained by only summing over those 
strings 
$\bold i\in\Sigma_{\smallG}^{n}$
that are multifractally relevant, i\.e\.
those 
strings 
$\bold i\in\Sigma_{\smallG}^{n}$
for which
$UL_{n}[\bold i]\subseteq C$.

\bigskip

\noindent
{\bf Remark.}
It is clear that if $C=X$, then the 
multifractal 
\lq\lq constraint"
$UL_{n}[\bold i]\subseteq C$
is vacuously
satisfied,
and that,
in this case,
the 
multifractal pressure and 
dynamical multifractal zeta-function
reduce to the usual pressure and the usual
dynamical zeta-function, i\.e\.
 $$
 \underline P_{X}^{U}(\varphi)
 =
 \overline P_{X}^{U}(\varphi)
 =
  P(\varphi)
 $$
and
 $$
 \zeta_{X}^{\dyn,U}(\varphi;\cdot)
 =
 \zeta^{\dyn}(\varphi;\cdot)\,.
 $$

\bigskip

Before developing
the theory of the multifractal pressure and the 
multifractal 
zeta-functions further
we make to following two simple observations.
Firstly, we note in Proposition 4.2 below that
the expected
relationship between 
the multifractal pressure and the radius of 
convergence of the multifractal zeta-function
holds.
Secondly,
we would expect any dynamically meaningful
 theory of
dynamical multifractal zeta-functions to
lead to
multifractal Bowen formulas.
For this to hold, we must, at the very least,
ensure that
there  are unique solutions to 
the 
relevant
multifractal Bowen equations. i\.e\.
 we must ensure
 that
 if $\Phi:\Sigma_{\smallG}^{\Bbb N}$ is a continuous function, then
there is are unique real numbers
 $\,\,\scri f(C)$ and
 $\,\,\scri F\,(C)$
solving the following
 multifractal Bowen 
 equations, namely,
 $$
 \aligned
 \lim_{r\searrow 0}\overline P_{B(C,r)}^{U}(\,\,\,\scri f(C)\,\Phi)
 &=
 0\,,\\
 \overline P_{C}^{U}(\,\scri F\,(C)\,\Phi)
 &=
 0\,.
 \endaligned
 \tag4.4
 $$ 
That there are 
unique numbers 
 $\,\,\scri f(C)$ and
 $\,\,\scri F\,(C)$
 satisfying (4.4)
 is our second preliminary observation, see Proposition 4.3.
However, before
stating and proving Proposition 4.2 and Proposition 4.3, we first 
show that the limit 
$\lim_{r\searrow 0}\overline P_{B(C,r)}^{U}(\,\,\,\scri f(C)\,\Phi)$ in (4.4)
exists. This is the content of the next proposition.

\bigskip

\proclaim{Proposition 4.1}
Let $X$ be a metric space and let $U:\Cal P(\Sigma_{\smallG}^{\Bbb N})\to X$ be 
continuous with respect to the weak topology.
Let $C\subseteq X$ be a  subset of $X$.
Fix a continuous function $\varphi:\Sigma^{\Bbb N}\to\Bbb R$.
Then the following limits
 $$
 \align
&\lim_{r\searrow 0}\underline P_{B(C,r)}^{U}(\varphi)\,,\\
&\lim_{r\searrow 0}\overline P_{B(C,r)}^{U}(\varphi)\,,
\endalign
 $$
 exist.
\endproclaim 
\noindent{\it Proof}\newline
\noindent
This follows immediately from the fact that the function
$r\to \overline P_{B(C,r)}^{U}(\varphi)$ is monotone.
\hfill$\square$

\bigskip

\noindent
We can now state Proposition 4.2 and Proposition 4.3.

\bigskip

\proclaim{Proposition 4.2. Radius of convergence}
Let $X$ be a metric space and let $U:\Cal P(\Sigma_{\smallG}^{\Bbb N})\to X$ be 
continuous with respect to the weak topology.
Let $C\subseteq X$ be a  subset of $X$.
Fix a continuous function $\varphi:\Sigma^{\Bbb N}\to\Bbb R$.
We have
 $$
 -
 \log
 \sigma_{\radius}
 \big(
 \,
 \zeta_{C}^{\dyn,U}(\varphi;\cdot)
 \,
 \big)
 =
 \overline P_{C}^{U}(\varphi)\,.
 $$
\endproclaim 
\noindent{\it Proof}\newline
\noindent
This follows immediately from the fact that
if
$(a_{n})_{n}$ is a sequence of complex
numbers
and if $f$ denotes the power series defined by
$f(z)=\sum_{n}a_{n}z^{n}$, then
$\sigma_{\radius}(f)
 =
 \frac{1}{\limsup_{n}|a_{n}|^{\frac{1}{n}}}$.
\hfill$\square$

\bigskip

\proclaim{Proposition 4.3. 
Continuity and monotonicity
of the multifractal pressure}
Let $X$ be a metric space and let $U:\Cal P(\Sigma_{\smallG}^{\Bbb N})\to X$ be 
continuous with respect to the weak topology.
Let $C\subseteq X$ be a  subset of $X$.
Fix a continuous
function
$\Phi:\Sigma_{\smallG}^{\Bbb N}\to\Bbb R$
with
$\Phi<0$.
Then the functions
$t\to
 \lim_{r\searrow 0}\overline P_{B(C,r)}^{U}(t\Phi)$
and 
$t
 \to
 \overline P_{C}^{U}(t\Phi)$,
 where
 $t\in\Bbb R$, 
are
continuous,
strictly decreasing and convex with
 $\lim_{t\to-\infty}
\lim_{r\searrow 0}\overline P_{B(C,r)}^{U}(t\Phi)
 =
\infty$
and
$\lim_{t\to\infty}
\lim_{r\searrow 0}\overline P_{B(C,r)}^{U}(t\Phi)
 =
-\infty$,
and
$\lim_{t\to-\infty}
\overline P_{C}^{U}(t\Phi) 
 =
\infty$
and
$\lim_{t\to\infty}
\overline P_{C}^{U}(t\Phi) 
 =
-\infty$.
In particular, there are unique real numbers
$\,\,\scri f(C)$
and
$\scri F\,(C)$ 
such that
 $$
 \align
 \lim_{r\searrow 0}\overline P_{B(C,r)}^{U}(\,\,\,\scri f(C)\,\Phi)
 &=
 0\,,
 \tag4.5\\
 \overline P_{C}^{U}(\,\scri F\,(C)\,\Phi)
&=
0\,;
\tag4.6
 \endalign
 $$
alternatively,
$\,\,\scri f(C)$
and
$\scri F\,(C)$ 
 are the  unique real numbers 
 such that
 $$
 \align
 \lim_{r\searrow 0}
 \sigma_{\radius}
 \big(
 \,
 \zeta_{B(C,r)}^{U}(\,\,\,\scri f(C)\,\Phi;\cdot)
 \,
 \big)
 &=
 1\,,\\
 \sigma_{\radius}
 \big(
 \,
 \zeta_{C}^{U}(\,\scri F\,(C)\,\Phi;\cdot)
 \,
 \big)
&=
1\,.
 \endalign
 $$ 
\endproclaim
\noindent{\it Proof}\newline 
\noindent
This statement is not difficult to prove and,
for the sake of brevity,
we have therefore decided 
to omit the proof.
\hfill$\square$

\bigskip

We will now state our main results.
The results are divided into two 
parts:
the first part consists of
Theorem 4.4 and Corollary 4.5,
and the second part consists
of
Theorem 4.6 and Corollary 4.7.
The motivation for this is the following.
Let $\Lambda$ denote the scaling map in (2.10).
For judicious 
choices of $X$ and $U$,
we are clearly attempting to relate the solutions
 $\,\,\scri f(C)$ and
 $\,\,\scri F\,(C)$
of the multifractal Bowen
 equations (4.5) and (4.6)
to various multifractal spectra.
The following simple example serves
to illustrate this.
Namely,
let $(\mu_{\smallvertexi})_{\smallvertexi\in\smallV}$ be the graph-directed self-conformal measures
in (2.7).
Next, define $X$ and $U:\Cal P(\Sigma_{\smallG}^{\Bbb N})\to X$
as follows.
Let $X=\Bbb R$.
Finally,
define $\Phi:\Cal P(\Sigma_{\smallG}^{\Bbb N})\to\Bbb R$
by
$$
\Phi(\bold i)
=
\log p_{\ini(\bold i)}
 $$
and let
 $$
 U\mu
  =
 \frac{\int\Phi\,d\mu}{\int\Lambda\,d\mu}\,.
 \tag4.7
 $$
Note that if $\bold i\in\Sigma_{\smallG}^{n}$, then
$UL_{n}[\bold i]
 =
 \{
 \frac{\log p_{\bold i}}{\log |DS_{\bold i}(\pi\bold u)|}
 \,|\,
 \bold u\in\Sigma_{\smallG}^{\Bbb N}
 \,\,
 \text{with}
 \,\,
 \termi(\bold i)=\ini(\bold u)
 \}$.
 It therefore follows that
 if $\varphi:\Sigma_{\smallG}^{\Bbb N}\to\Bbb R$ is a continuous function, then
  $$
  \align
  \overline P_{C}^{U}(\varphi)
 &=
  \limsup_{n}
  \frac{1}{n}
  \log
  \left(
  \,\,\,\,\,\,
    \sum
  \Sb
   \bold i\in\Sigma_{\smallsmallG}^{n}\\
   {}\\
   \forall
 \bold u\in\Sigma_{\smallsmallG}^{\Bbb N}
 \,\,
 \text{with}
 \,\,
 \termi(\bold i)=\ini(\bold u)
 \,\,:\,\,
 \frac{\log p_{\bold i}}{\log |DS_{\bold i}(\pi\bold u)|}
\in
C
  \endSb
  \,\,\,\,\,\,
  \sup_{\bold u\in[\bold i]}
  \,\,
  \exp
  \sum_{k=0}^{n-1}
  \varphi S^{k}\bold u
  \right)\,,
  \tag4.8\\
  &{}\\
   \overline P_{B(C,r)}^{U}(\varphi)
 &=
  \limsup_{n}
  \frac{1}{n}
  \log
  \left(
    \sum
  \Sb
   \bold i\in\Sigma_{\smallsmallG}^{n}\\
   {}\\
   \forall
 \bold u\in\Sigma_{\smallsmallG}^{\Bbb N}
 \,\,
 \text{with}
 \,\,
 \termi(\bold i)=\ini(\bold u)
 \,\,:\,\,
 \frac{\log p_{\bold i}}{\log |DS_{\bold i}(\pi\bold u)|}
\in
B(C,r)
  \endSb
  \sup_{\bold u\in[\bold i]}
  \,\,
  \exp
  \sum_{k=0}^{n-1}
  \varphi S^{k}\bold u
  \right)\,,
  \tag4.9
  \endalign
  $$
and 
  $$
  \zeta_{C}^{\dyn,U}(\varphi;z)
  =
  \sum_{n}
  \,\,
  \frac{z^{n}}{n}
  \,\,
  \left(
  \sum
  \Sb
   \bold i\in\Sigma_{\smallsmallG}^{n}\\
   {}\\
   \forall
 \bold u\in\Sigma_{\smallsmallG}^{\Bbb N}
 \,\,
 \text{with}
 \,\,
 \termi(\bold i)=\ini(\bold u)
 \,\,:\,\,
 \frac{\log p_{\bold i}}{\log |DS_{\bold i}(\pi\bold u)|}
\in
C
  \endSb
  \sup_{\bold u\in[\bold i]}
  \,\,
  \exp
  \sum_{k=0}^{n-1}
  \varphi S^{k}\bold u
  \right)
  \,.
  \tag4.10
  $$ 
For $\alpha\in\Bbb R$ and $C=\{\alpha\}$,
we are now attempting to relate 
the solutions
 $\,\scr f\,\,(\alpha)$ and
 $\scr F\,\,(\alpha)$
of the following multifractal Bowen
 equations
 $$
 \align
 \lim_{r\searrow 0}\overline P_{B(C,r)}^{U}(\,\,\scr f\,\,(\alpha)\,\Lambda)
 &=
 0\,,
 \tag4.11\\
 \overline P_{C}^{U}(\,\scr F\,\,(\alpha)\,\Lambda)
 &=
 0\,,
 \tag4.12
 \endalign
 $$ 
to the multifractal spectrum $f_{\mu_{\smallsmallvertexi}}(\alpha)$ of
$\mu_{\smallvertexi}$;
observe that the existence and uniqueness of the solutions
 $\,\scr f\,\,(\alpha)$ and
 $\scr F\,\,(\alpha)$
 to (4.11) and (4.12) follow from
Proposition 4.3.
However, it is clear that if
$C=\{\alpha\}$, then
the
sum in (4.8) may be empty,
and the 
solution
 $\scr F\,\,(\alpha)$
 to the equation
 $ \overline P_{C}^{U}(\,\scr F\,\,(\alpha)\,\Lambda)
 =
 0$
is therefore
 $-\infty$, i\.e\.
$\scr F\,\,(\alpha)=-\infty$.
Hence, it may happen that
 $$
 \scr F\,\,(\alpha)=-\infty
 <
 f_{\mu_{\smallsmallvertexi}}(\alpha)\,.
 \tag4.13
 $$
It follows from this discussion that
if $C=\{\alpha\}$, 
then the 
pressure (4.8)
and  the zeta-function (4.10)
do not, in general, encode
sufficient
information
allowing us to recover the multifractal spectrum
$f_{\mu_{\smallvertexi}}(\alpha)$.
The reason for the strict inequality in (4.13) 
is, of course, clear:
even though
there are no
strings
$\bold i\in\Sigma_{\smallG}^{*}$
and
$\bold u\in\Sigma_{\smallG}^{\Bbb N}$
with
$\termi(\bold i)=\ini(\bold u)$
for which
the ratio
$ \frac{\log p_{\bold i}}{\log |DS_{\bold i}(\pi\bold u)|}$ equals $\alpha$,
there may
nevertheless
be
many sequences $(\bold i_{n})_{n}$ of strings
$\bold i_{n}\in\Sigma^{*}$
for which the sequence
of
sets
$\{
\frac{\log p_{\bold i_{n}}}{\log |DS_{\bold i_{n}}(\pi\bold u)|}
\,|\,
\bold u\in\Sigma_{\smallG}^{\Bbb N}
\,\,
\text{with}
\,\,
\termi(\bold i_{n})=\ini(\bold u)
\}$
\lq\lq shrinks"
to the singleton $\{\alpha\}$.
In order to capture this,
we can
proceed in two equally natural ways.
Either, we can 
consider
a family
of enlarged
\lq\lq target" sets 
shrinking to the original 
main \lq\lq target" $\{\alpha\}$;
this approach will be referred to as the
shrinking target approach.
Or, alternatively,
we can consider
a
fixed enlarge 
\lq\lq target" set
and regard this as our
original main \lq\lq target";
this approach will be referred to as the
fixed target approach.
Indeed,
the statement of our main  results below
is divided into two 
parts
paralleling the above discussion, namely:
the first part 
(consisting of
Theorem 4.4 and Corollary 4.5)
presents our results in the shrinking target setting,
and the 
second part 
consisting of 
Theorem 4.6 and Corollary 4.7)
presents our results in the fixed target setting.

\bigskip

\newpage

{\bf Statement of main results in the shrinking target setting.}
In the
shrinking target
setting,
Theorem 4.4  provide
a
variational principle for the 
multifractal pressure
and 
Corollary 4.5
provide a
variational principle
for
the solution
$\,\,\scri f(C)$ to the multifractal Bowen equation (4.5).
Below we 
denote the entropy
of a measure $\mu\in\Cal P(\Sigma_{\smallG}^{\Bbb N})$
by $h(\mu)$.

\bigskip

\proclaim{Theorem 4.4. 
The shrinking target variational principle
for the multifractal pressure}
Let $X$ be a metric space and let $U:\Cal P(\Sigma_{\smallG}^{\Bbb N})\to X$ be 
continuous with respect to the weak topology.
Let $C\subseteq X$ be a  subset of $X$.
Fix a continuous function $\varphi:\Sigma_{\smallG}^{\Bbb N}\to\Bbb R$.

\roster
\item"(1)"
We have
 $$
 \align
 \lim_{r\searrow0}
\underline P_{B(C,r)}^{U}(\varphi)
&=
 \lim_{r\searrow0}
\overline P_{B(C,r)}^{U}(\varphi)
=
\sup
\Sb
\mu\in\Cal P_{S}(\Sigma_{\smallsmallG}^{\Bbb N})\\
{}\\
U\mu\in \overline C
\endSb
\Bigg(
h(\mu)+\int\varphi\,d\mu
\Bigg)\,.\\
\endalign
$$ 
\item"(2)"
We have
 $$
 \align
 \lim_{r\searrow0}
-\log
\sigma_{\radius}\big(\,\zeta_{B(C,r)}^{\dyn,U}(\varphi;\cdot)\,\big)
&=
\sup
\Sb
\mu\in\Cal P_{S}(\Sigma_{\smallsmallG}^{\Bbb N})\\
{}\\
U\mu\in \overline C
\endSb
\Bigg(
h(\mu)+\int\varphi\,d\mu
\Bigg)\,.
\endalign
$$ 
\endroster
\endproclaim

\bigskip

\noindent
Theorem 4.4 is proved in Section 9
using techniques from 
large deviation theory
developed
 in Sections 6--8.
 Observe that if we let
 $C=X$ in Theorem 4.4, then
 the multifractal pressure 
 equals the
  usual pressure, i\.e\.
 $\overline P_{B(C,r)}^{U}(\varphi)
 =
 \underline P_{B(C,r)}^{U}(\varphi)
 =
 P(\varphi)$,
 and
  the variational principle in Theorem 4.4.(1)
  therefore
 simplifies to the usual variational principle, namely,
 $$
  P(\varphi)
=
\sup
\Sb
\mu\in\Cal P_{S}(\Sigma_{\smallsmallG}^{\Bbb N})\\
\endSb
\Bigg(
h(\mu)+\int\varphi\,d\mu
\Bigg)\,.
$$

\bigskip

\proclaim{Corollary 4.5.
The shrinking target
multifractal Bowen equation}
Let $X$ be a metric space and let $U:\Cal P(\Sigma_{\smallG}^{\Bbb N})\to X$ be 
continuous with respect to the weak topology.
Let $C\subseteq X$ be a  subset of $X$.
Fix a continuous function  $\Phi:\Sigma_{\smallG}^{\Bbb N}\to\Bbb R$ 
with
$\Phi<0$
and let
$\,\,\scri f(C)$
be the unique real number
such that
 $$
 \align
 \lim_{r\searrow 0}\overline P_{B(C,r)}^{U}(\,\,\,\scri f(C)\,\Phi)
 &=
 0\,;
 \endalign
 $$ 
alternatively, 
 $\,\,\scri f(C)$
is the unique real number
such that
 $$
 \align
 \lim_{r\searrow 0}
 \sigma_{\radius}
 \big(
 \,
 \zeta_{B(C,r)}^{U}(\,\,\scri f(C)\,\Phi;\cdot)
 \,
 \big)
 &=
 1\,.
 \endalign
 $$ 
Then
 $$
\scri f(C)
=
\sup
\Sb
\mu\in\Cal P_{S}(\Sigma_{\smallsmallG}^{\Bbb N})\\
{}\\
U\mu\in \overline C
\endSb
-
\frac{h(\mu)}{\int\Phi\,d\mu}\,.
$$ 
\endproclaim 
\noindent{\it Proof}\newline
\noindent
It follows from Theorem 4.4 and the definition of
$\,\,\scri f(C)$
that
 $$
 \align
\sup
\Sb
\mu\in\Cal P_{S}(\Sigma_{\smallsmallG}^{\Bbb N})\\
{}\\
U\mu\in \overline C
\endSb
\Bigg(
h(\mu)+\,\,\scri f(C)\int\Phi\,d\mu
\Bigg)
&=
 \lim_{r\searrow 0}\overline P_{B(C,r)}^{U}(\,\,\,\scri f(C)\,\Phi)
 =
0\,.
\tag4.14
\endalign
$$ 
The desired formula for $\,\,\scri f(C)$
follows easily from (4.14).
\hfill$\square$

\bigskip

\newpage

{\bf Statement of  main results in the fixed target setting.}
Of course, if the set $C$ is 
\lq\lq too small",
then
it follows from the discussion 
following
the statement of Proposition 4.3
that we, in general, cannot expect any
meaningful results 
in the fixed target setting.
However, if the set $C$ satisfies a 
non-degeneracy 
condition guaranteeing that
it is not
\lq\lq too small"
(namely condition (4.15) below),
then 
meaningful results can be obtained
in the fixed target setting.
This is the contents of
Theorem 4.6 and Corollary 4.7 below.
Indeed,
Theorem 4.6 and Corollary 4.7 provide
variational principles for the 
multifractal pressure
and for the solution
$\,\,\scri F\,(C)$ to the multifractal Bowen equation (4.6)
 in the
fixed target
setting.


\medskip

\proclaim{Theorem 4.6.
The fixed target variational principle
for the multifractal pressure}
Let $X$ be a normed vector space.
Let $\Gamma:\Cal P(\Sigma_{\smallG}^{\Bbb N})\to X$
be continuous and affine
and let
$\Delta:\Cal P(\Sigma_{\smallG}^{\Bbb N})\to \Bbb R$
be continuous and affine
with
$\Delta(\mu)\not=0$
for all $\mu\in\Cal P(\Sigma_{\smallG}^{\Bbb N})$.
Define 
$U:\Cal P(\Sigma_{\smallG}^{\Bbb N})\to X$
by
$U=\frac{\Gamma}{\Delta}$.
Let $C$ be a closed and convex subset of $X$ and assume that
 $$
 \overset{\,\circ}\to{C}
 \cap
 \,
 U\big(\,\Cal P_{S}(\Sigma_{\smallG}^{\Bbb N})\,\big)
 \not=
 \varnothing\,.
 \tag4.15
 $$
\roster
\item"(1)"
We have
$$
\align
{}\qquad\qquad
  \quad\,\,\,\,
  \,\,\,\,\,\,\,
P_{C}^{U}(\varphi)
&=
\sup
\Sb
\mu\in \Cal P_{S}(\Sigma_{\smallsmallG}^{\Bbb N})\\
{}\\
U\mu\in C
\endSb
\Bigg(
h(\mu)+\int\varphi\,d\mu
\Bigg)
=
\sup
\Sb
\mu\in \Cal P_{S}(\Sigma_{\smallsmallG}^{\Bbb N})\\
U\mu\in \overset{\,\circ}\to{C}
\endSb
\Bigg(
h(\mu)+\int\varphi\,d\mu
\Bigg)\,.
\endalign
$$

 \item"(2)"
We have
$$
\align
-\log
\sigma_{\radius}\big(\,\zeta_{C}^{\dyn,U}(\varphi;\cdot)\,\big)
&=
\sup
\Sb
\mu\in \Cal P_{S}(\Sigma_{\smallsmallG}^{\Bbb N})\\
{}\\
U\mu\in C
\endSb
\Bigg(
h(\mu)+\int\varphi\,d\mu
\Bigg)
=
\sup
\Sb
\mu\in \Cal P_{S}(\Sigma_{\smallsmallG}^{\Bbb N})\\
U\mu\in \overset{\,\circ}\to{C}
\endSb
\Bigg(
h(\mu)+\int\varphi\,d\mu
\Bigg)\,.
\endalign
$$
\endroster
\endproclaim


\medskip

\noindent
Theorem 4.6 is proved in Section 10
using techniques from 
large deviation theory
developed
 in Sections 6--8.
 Again, we
  observe that if we let
 $C=X$ in Theorem 4.6, then
 the multifractal pressure 
 equals the
  usual pressure, i\.e\.
 $P_{C}^{U}(\varphi)
 =
 P(\varphi)$,
 and
  the variational principle in Theorem 4.6.(1)
  therefore
 simplifies to the usual variational principle, namely,
 $$
  P(\varphi)
=
\sup
\Sb
\mu\in\Cal P_{S}(\Sigma_{\smallsmallG}^{\Bbb N})\\
\endSb
\Bigg(
h(\mu)+\int\varphi\,d\mu
\Bigg)\,.
$$


\medskip

\proclaim{Corollary 4.7.
The fixed target
multifractal Bowen equation}
Let $X$ be a normed vector space.
Let $\Gamma:\Cal P(\Sigma_{\smallG}^{\Bbb N})\to X$
be continuous and affine
and let
$\Delta:\Cal P(\Sigma_{\smallG}^{\Bbb N})\to \Bbb R$
be continuous and affine
with
$\Delta(\mu)\not=0$
for all $\mu\in\Cal P(\Sigma_{\smallG}^{\Bbb N})$.
Define 
$U:\Cal P(\Sigma_{\smallG}^{\Bbb N})\to X$
by
$U=\frac{\Gamma}{\Delta}$.
Let $C$ be a closed and convex subset of $X$ and assume that
 $$
 \overset{\,\circ}\to{C}
 \cap
 \,
 U\big(\,\Cal P_{S}(\Sigma_{\smallG}^{\Bbb N})\,\big)
 \not=
 \varnothing\,.
 $$
Let $\Phi:\Sigma_{\smallG}^{\Bbb N}\to\Bbb R$ be continuous
with
$\Phi<0$.
Let 
$\scri F\,(C)$
be the unique real number
such that
 $$
 \align
 P_{C}^{U}(\,\scri F\,(C)\,\Phi)
 &=
 0\,;
 \endalign
 $$ 
 alternatively, 
 $\,\,\scri f(C)$
is the unique real number
such that
 $$
 \align
 \sigma_{\radius}
 \big(
 \,
 \zeta_{C}^{U}(\,\,\scri F\,(C)\,\Phi;\cdot)
 \,
 \big)
 &=
 1\,.
 \endalign
 $$ 
Then
 $$
\scri F\,(C)
=
\sup
\Sb
\mu\in\Cal P_{S}(\Sigma_{\smallsmallG}^{\Bbb N})\\
{}\\
U\mu\in C
\endSb
-
\frac{h(\mu)}{\int\Phi\,d\mu}\,.
$$ 
\endproclaim 
\noindent{\it Proof}\newline
\noindent
The proof is similar to the proof of Corollary 4.5 using
Theorem 4.6 and the definition of
$\scri F\,(C)$.
\hfill$\square$


\medskip

In the next section we will apply
Theorem 4.4, Corollary 4.5, Theorem 4.6 and Corollary 4.7
to
show that in many cases,
 the solutions 
$\,\,\scri f(C)$ 
and
$\,\,\scri F\,(C)$ to the multifractal Bowen equations 
(4.5)
and (4.6) coincide with 
various well-known multifractal spectra.

\bigskip

\newpage

{\bf Outline of the proofs of Theorem 4.4 and Theorem 4.6.}
Below we give an
outline illustrating the key ideas
in the proofs of the two main results, namely, 
Theorem 4.4 and Theorem 4.6.
We first note 
that Theorem 4.6 follows from
 Theorem 4.4
by a 
\lq\lq continuity"
argument; see Section 10 for the details of this argument.
We will now give an outline of the proof of Theorem 4.4.
Recall the statement of Theorem 4.4.
Namely, let
$X$ be a metric space and let $U:\Cal P(\Sigma_{\smallG}^{\Bbb N})\to X$ be 
continuous with respect to the weak topology.
Also, let $C\subseteq X$ be a  subset of $X$.
Theorem 4.4 now says that
if $\varphi:\Sigma_{\smallG}^{\Bbb N}\to\Bbb R$ is a continuous function,
then
 $$
 \align
 \lim_{r\searrow0}
\underline P_{B(C,r)}^{U}(\varphi)
&=
 \lim_{r\searrow0}
\overline P_{B(C,r)}^{U}(\varphi)
=
\sup
\Sb
\mu\in\Cal P_{S}(\Sigma_{\smallsmallG}^{\Bbb N})\\
{}\\
U\mu\in \overline C
\endSb
\Bigg(
h(\mu)+\int\varphi\,d\mu
\Bigg)
\tag4.16
\endalign
$$ 
and
 $$
 \align
 \lim_{r\searrow0}
-\log
\sigma_{\radius}\big(\,\zeta_{B(C,r)}^{\dyn,U}(\varphi;\cdot)\,\big)
&=
\sup
\Sb
\mu\in\Cal P_{S}(\Sigma_{\smallsmallG}^{\Bbb N})\\
{}\\
U\mu\in \overline C
\endSb
\Bigg(
h(\mu)+\int\varphi\,d\mu
\Bigg)\,.
\tag4.17
\endalign
$$ 
It is clear that (4.17) follows from (4.16)
(see Proposition 4.2),
and it therefore suffices to prove (4.16).
Writing
 $$
 D(C)
 =
 \sup
\Sb
\mu\in\Cal P_{S}(\Sigma_{\smallsmallG}^{\Bbb N})\\
{}\\
U\mu\in \overline C
\endSb
\Bigg(
h(\mu)+\int\varphi\,d\mu
\Bigg)\,,
$$
then (4.16)
can be written as
 $$
 \aligned
\lim_{r\searrow 0}
 \,
  \liminf_{n}
 \,
  \,\,
  \frac{1}{n}
 \,\,
 \log
   \sum
  \Sb
  \bold i\in\Sigma_{\smallsmallG}^{n}\\
  {}\\
  UL_{n}[\bold i]\subseteq B(C,r)
  \endSb
   \,\,
  \sup_{\bold u\in[\bold i]}
  \,\,
 \exp
 \,\,
 \sum_{k=0}^{n-1}\varphi S^{k}\bold u
 &=
 D(C)\,,\\
\lim_{r\searrow 0}
\,
  \limsup_{n}
  \,\,
  \frac{1}{n}
 \,\,
 \log
   \sum
  \Sb
 \bold i\in\Sigma_{\smallsmallG}^{n}\\
  {}\\
  UL_{n}[\bold i]\subseteq B(C,r)
  \endSb
 \,\,
  \sup_{\bold u\in[\bold i]}
  \,\,
 \exp
 \,\,
 \sum_{k=0}^{n-1}\varphi S^{k}\bold u
&=
 D(C)\,.
 \endaligned
 \tag4.18
 $$

We will use
techniques from  large deviation
theory in order to analyse the asymptotic behaviour of the sums on the left hand
side
of (4.18).
In particular, we will use one of the most celebrated results from
large deviation
theory, namely, 
Varadhan's integral lemma.
This result says that
 if $X$
is a complete separable metric space
and 
$(P_{n})_{n}$ is a sequence of probability measures on
$X$ satisfying the large deviation property
with
rate constants $a_{n}\in(0,\infty)$ for $n\in\Bbb N$
and rate function $I:\Bbb R\to[-\infty,\infty]$
(this terminology will be explained in Section 7),
then
any
bounded continuous function $F:X\to \Bbb R$
satisfies the following:
$$
\lim_{n}
\frac{1}{a_{n}}
\log
\int
\exp(a_{n}F)\,dP_{n}
=
-
\inf_{x\in X}(\,I(x)-F(x))
\tag4.19
$$
and 
if  we define the probability measure $Q_{n}$ on $X$ by
 $$
 Q_{n}(E)
 =
 \,\,
 \frac
 {\int_{E}\exp(a_{n}F)\,dP_{n}}
 {\int\exp(a_{n}F)\,dP_{n}}
 \quad
 \text{for Borel subsets $E$ of $X$,}
 \tag4.20
 $$
then the sequence $(Q_{n})_{n}$
has the large deviation property with constants $(a_{n})_{n}$
and rate function
$(I-F)-\inf_{x\in X, I(x)<\infty}(I(x)-F(x))$.

In order to use the above results to
analyse the asymptotic
behaviour of the sums
on the left hand side of (4.18) we
must construct 
a continuous function
$\widetilde F_{\varphi}:\Cal P(\Sigma_{\smallG}^{\Bbb N})\to\Bbb R$ 
and
a sequence of measures
$\widetilde\Pi_{n}$
having the large deviation property
with constants $(n)_{n}$
such that the sum
 $$
  \sum
  \Sb
 \bold i\in\Sigma_{\smallsmallG}^{n}\\
  {}\\
  UL_{n}[\bold i]\subseteq B(C,r)
  \endSb
 \,\,
  \sup_{\bold u\in[\bold i]}
  \,\,
 \exp
 \,\,
 \sum_{k=0}^{n-1}\varphi S^{k}\bold u
 \tag4.21
$$
can be expressed 
in terms of the integral
 $$
 \int\exp(n \widetilde F_{\varphi})\,d\widetilde \Pi_{n}
 $$
and the measure $\widetilde Q_{\varphi,n}$ defined by
$$
 \widetilde Q_{\varphi,n}(E)
 =
 \,\,
 \frac
 {\int_{E}\exp(n\widetilde F_{\varphi})\,d\widetilde\Pi_{n}}
 {\int\exp(n\widetilde F_{\varphi})\,d\widetilde\Pi_{n}}
 \quad
 \text{for Borel subsets $E$ of $\Cal P(\Sigma_{\smallG}^{\Bbb N})$.}
 $$
In fact, in stead of working with
the sum in (4.21), we
will (for technical reasons that are explaining in Section 6) work with a slightly modified sum,
namely the sum where the function $L_{n}$
has been replaced
by a slightly different function
$M_{n}$.
We will thus
construct
a continuous function
$F_{\varphi}:\Cal P(\Sigma_{\smallG}^{\Bbb N})\to\Bbb R$ 
and
a sequence of measures
$\Pi_{n}$
having the large deviation property
with constants $(n)_{n}$
such that the sum
 $$
  \sum
  \Sb
 \bold i\in\Sigma_{\smallsmallG}^{n}\\
  {}\\
  UM_{n}[\bold i]\subseteq B(C,r)
  \endSb
 \,\,
  \sup_{\bold u\in[\bold i]}
  \,\,
 \exp
 \,\,
 \sum_{k=0}^{n-1}\varphi S^{k}\bold u
 \tag4.22
$$
can be expressed 
in terms of the integral
 $$
 \int\exp(n F_{\varphi})\,d\Pi_{n}
 $$
and the measure $Q_{\varphi,n}$ defined by
$$
Q_{\varphi,n}(E)
 =
 \,\,
 \frac
 {\int_{E}\exp(n F_{\varphi})\,d\Pi_{n}}
 {\int\exp(n F_{\varphi})\,d\Pi_{n}}
 \quad
 \text{for Borel subsets $E$ of $\Cal P(\Sigma_{\smallG}^{\Bbb N})$.}
 $$
 Finally, we will 
 obtain the asymptotic behaviour
 of the original sum in (4.21) from the 
  the asymptotic behaviour
 of the modified sum in (4.22).
The overall implementation of this strategy 
is divided into  4 steps described  below.

\smallskip

{\bf Step 1 (Section 6):
The map $M_{n}$.}

We define the map
$M_{n}:\Sigma_{\smallG}^{\Bbb N}\to \Cal P(\Sigma_{\smallG}^{\Bbb N})$
and prove a number of 
\lq\lq continuity"
results showing that the maps $M_{n}$ and $L_{n}$
are \lq\lq close" together.

\smallskip

{\bf Step 2 (Section 7):
The measures $\Pi_{n}$.}

Using the map $M_{n}$,
we define the measures $\Pi_{n}$
and prove (using a result from Orey \& Pelikan [OrPe2])
that
the sequence $(\Pi_{n})_{n}$
has the large deviation property
with constants $(n)_{n}$ and rate function $I$ given by
  $$
 I(\mu)
 =
 \cases
 \log\lambda-h(\mu)
\quad
 &\text{for $\mu\in \Cal P_{S}(\Sigma_{\smallG}^{\Bbb N})$;}\\
 \infty
\quad
 &\text{for 
 $\mu\in \Cal P(\Sigma_{\smallG}^{\Bbb N})\setminus \Cal P_{S}(\Sigma_{\smallG}^{\Bbb N})$.
 }
 \endcases
$$
where $\lambda$ is a constant (in fact, $\lambda$
is the entropy of the Parry measure
on $\Sigma_{\smallG}^{\Bbb N}$).

\smallskip

\newpage

{\bf Step 3 (Section 8):
The asymptotic behaviour of the modified sum (4.22).}

Define the function
 $F_{\varphi}:\Cal P_{S}(\Sigma_{\smallG}^{\Bbb N})\to\Bbb R$  by
$F_{\varphi}(\mu)=\int\varphi\,d\mu$.
Let the measure $\Pi_{n}$ be as in Step  2 and define
 the
probability measure
$Q_{\varphi,n}$ on $\Cal P(\Sigma_{\smallG}^{\Bbb N})$ by
 $$
 \align
 Q_{\varphi,n}(E)
&=
 \frac
 {\int_{E}\exp(nF_{\varphi})\,d\Pi_{n}}
 {\int\exp(nF_{\varphi})\,d\Pi_{n}}
 \quad
 \text{for Borel subsets $E$ of $\Cal P(\Sigma_{\smallG}^{\Bbb N})$.}
 \endalign
 $$
We first show that the modified sum (4.22) can be expressed in terms of the 
integral
$\int\exp(nF_{\varphi})\,d\Pi_{n}$
and the 
measure $Q_{\varphi,n}$;
more precisely, we show
that there is a positive 
constant $c$ 
such that
 $$
 \align
 \sum
   \Sb
   \bold k\in\Sigma_{\smallsmallG}^{n}\\
   {}\\
   UM_{n}[\bold k]\subseteq C
	 \endSb
 \sup_{\bold u\in[\bold k]}
 \,\,
 \exp
 \sum_{i=0}^{n-1}\varphi S^{i}\bold u
&\le
c
\,\,
\lambda^{n}
\,\,
Q_{\varphi,n}\Big(\{U\in C\}\Big)
\,\,
\int\exp(nF_{\varphi})\,d\Pi_{n}\,,\\
{}\\
 \sum
   \Sb
    \bold i\in\Sigma_{\smallsmallG}^{n}\\
   {}\\
   UM_{n}[\bold k]\subseteq C
	 \endSb
 \sup_{\bold u\in[\bold k]}
 \,\,	 
 \exp
 \sum_{i=0}^{n-1}\varphi S^{i}\bold u
&\ge
\frac{1}{c}
\,\,
\lambda^{n}
\,\,
Q_{\varphi,n}\Big(\{U\in C\}\Big)
\,\,
\int\exp(nF_{\varphi})\,d\Pi_{n}\,,
\endalign
 $$
 for all $n$.
Using Varadhan's integral lemma, this is easily seen to imply that
if $G$ is an open subset of $X$
with
$U^{-1}G\cap\Cal P_{S}(\Sigma_{\smallG}^{\Bbb N})\not=\varnothing$, then
 $$
 \liminf_{n}
 \,
 \frac{1}{n}
 \,
 \log
 \sum
   \Sb
   \bold k\in\Sigma_{\smallsmallG}^{n}\\
   {}\\
   UM_{n}[\bold k]\subseteq G
	 \endSb
 \sup_{\bold u\in[\bold k]}
 \,\,
 \exp
 \sum_{i=0}^{n-1}\varphi S^{i}\bold u
\ge
D(G)\,,
\tag4.23
$$
and
if $K$ is a closed subset of $X$
with
$U^{-1}K\cap\Cal P_{S}(\Sigma_{\smallG}^{\Bbb N})\not=\varnothing$, then
$$
 \limsup_{n}
 \,
 \frac{1}{n}
 \,
 \log
 \sum
   \Sb
   \bold k\in\Sigma_{\smallsmallG}^{n}\\
   {}\\
   UM_{n}[\bold k]\subseteq K
	 \endSb
 \sup_{\bold u\in[\bold k]}
 \,\,
 \exp
 \sum_{i=0}^{n-1}\varphi S^{i}\bold u
\le
D(K)\,.
\tag4.24
$$

\smallskip

{\bf Step 4 (Section 9):
The asymptotic behaviour of the  (non-modified) sum (4.21).}

Finally, 
using the 
\lq\lq continuity" results from Step 1 
(showing that 
the maps $M_{n}$ and $L_{n}$
are \lq\lq close" together),
the desired
asymptotic behaviour,
i\.e\.  (4.18),
 of the 
(non-modified) sum (4.21)
can be
obtained
from 
the asymptotic behaviour,
i\.e\.  (4.23) and (4.24),
 of the 
modified sum (4.22)
established in Step 3.

  
\newpage


\centerline{\smc 5. Applications: }

\centerline{\smc multifractal spectra of measures}

\centerline{\smc  and}

\centerline{\smc  
multifractal spectra of ergodic Birkhoff averages}

 \medskip

We will now consider several of applications of 
Theorem 4.4, Corollary 4.5, Theorem 4.6 and Corollary 4.7 
to multifractal spectra of measures and ergodic averages.
In particular, we consider the following examples:

\medskip

\roster

\item"$\bullet$"
 Section 5.1: Multifractal spectra of graph-directed self-conformal measures.
 
\medskip

\item"$\bullet$"
Section 5.2: Multifractal spectra of ergodic Birkhoff averages
of continuous functions on 
graph-directed self-conformal sets.

\endroster

\medskip

{\bf 5.1.
Multifractal 
spectra of 
graph directed 
self-conformal measures.}
Let
$(\,
 \V,\,\allowmathbreak
 \E,\,\allowmathbreak
 (V_{\smallvertexi})_{\smallvertexi\in\smallV},\,\allowmathbreak
 (X_{\smallvertexi})_{\smallvertexi\in\smallV},\,\allowmathbreak
 (S_{\smalledge})_{\smalledge\in\smallE},\,\allowmathbreak
 (p_{\smalledge})_{\smalledge\in\smallE}\,
 )$
 be a
 graph-directed conformal iterated function system with probabilities
 (see Section 2.2)
 and let
   $(K_{\smallvertexi})_{\smallvertexi\in\smallV}$
   and 
   $(\mu_{\smallvertexi})_{\smallvertexi\in\smallV}$
   be the list of graph-directed self-conformal 
   sets 
   and the list of graph-direted self-conformal
   measures
   associated with the list 
 $(\,
 \V,\,\allowmathbreak
 \E,\,\allowmathbreak
 (V_{\smallvertexi})_{\smallvertexi\in\smallV},\,\allowmathbreak
 (X_{\smallvertexi})_{\smallvertexi\in\smallV},\,\allowmathbreak
 (S_{\smalledge})_{\smalledge\in\smallE},\,\allowmathbreak
 (p_{\smalledge})_{\smalledge\in\smallE}\,
 )$, respectively,
   i\.e\.
   the sets in the list
  $(K_{\smallvertexi})_{\smallvertexi\in\smallV}$
  are the unique non-empty
  compact sets
  satisfying (2.6)
   and 
   the measures in the list
   $(\mu_{\smallvertexi})_{\smallvertexi\in\smallV}$
   are the unique probability measures
   satisfying (2.7).
Recall that
the Hausdorff multifractal spectrum
$f_{\mu_{\smallsmallvertexi}}$ of $\mu_{\smallvertexi}$ 
is defined
by
 $$
 \align
 f_{\mu_{\smallsmallvertexi}}(\alpha)
&=
   \,\dim_{\Haus}
   \left\{x\in K_{\smallvertexi}
    \,\left|\,
     \lim_{r\searrow0}
     \frac
     {\log\mu_{\smallvertexi} (B(x,r))}{\log r}
     =
     \alpha
      \right.
   \right\}\,,
  \endalign
 $$
for $\alpha\in\Bbb R$, see Section 1.
If the OSC is satisfied,
then
the 
multifractal spectrum
$f_{\mu_{\smallsmallvertexi}}(\alpha)$
can be computed as follows.
Define
$\Phi:\Sigma_{\smallG}^{\Bbb N}\to\Bbb R$ by
 $$
 \align
 \Phi(\bold i)
 &=
 \log p_{\ini(\bold i)}
 \endalign
 $$
for 
$\bold i=\edge_{1}\edge_{2}\ldots
\in\Sigma_{\smallG}^{\Bbb N}$
and recall that the map 
$\Lambda:\Sigma_{\smallG}^{\Bbb N}\to\Bbb R$
is defined by
 $$
 \align
 \qquad\quad\,\,
 \Lambda(\bold i)
&=
 \log
 |DS_{\smalledge_{1}}(\pi S\bold i)|
 \endalign
 $$
for 
$\bold i=\edge_{1}\edge_{2}\ldots
\in\Sigma_{\smallG}^{\Bbb N}$.
Next,
we define the function $\beta:\Bbb R\to\Bbb R$
by
 $$
 P\big(
 \,
 \beta(q)\Lambda
 +
 q\Phi
 \,
 \big)
 =
 0\,;
 \tag5.1
 $$
 alternatively,
 the function
 $\beta:\Bbb R\to\Bbb R$
is defined by
 $$
 \sigma_{\radius}
 \big(
 \,
 \zeta^{\dyn}(q\Phi+\beta(q)\Lambda;\cdot)
 \,
 \big)
 =
 1\,.
 \tag5.2
 $$
 We note that
if the all maps 
$S_{\smalledge}$ are similarities, i\.e\.
if for 
each 
$\edge\in \E$
there is a number $r_{\smalledge}\in(0,1)$ such that
 $$
 |S_{\smalledge}(x)-S_{\smalledge}(y)|
 =
 r_{\smalledge}|x-y|
 $$
 for all
 $x,y\in X_{\termi(\smalledge)}$,
 then
 there is 
an alternative characterisation  
 of the function $\beta$.
 Namely, in this case the function $\beta$
 is given by the following.
For $q,t\in\Bbb R$, define
the matrix 
$A(q,t)=\big(\,a_{\smallvertexi,\smallvertexj}(q,t)\,\big)_{\smallvertexi,\smallvertexj\in\smallV}$ by
 $$
 a_{\smallvertexi,\smallvertexj}(q,t)
 =
 \sum_{\smalledge\in\smallE_{\smallsmallvertexi,\smallsmallvertexj}}
 p_{\smalledge}^{q}r_{\smalledge}^{t}\,.
 $$
 For $q\in\Bbb R$, the number $\beta(q)$ is now the unique
 real number such that
  $$
 \rho_{\specradsmall}\,A(q,\beta(q))
  =
  1\,;
  $$
 here and below we use the following notation, namely,
 if $M$ is a square matrix, then we will write
 $\rho_{\specradsmall} \,M$ for the spectral radius of $M$.
  If the OSC is satisfied, then
the multifractal spectrum $f_{\mu_{\smallsmallvertexi}}$
of $\mu_{\smallvertexi}$
can be computer as follows, see Theorem D  below.
This result
was first established by Edgar \& Mauldin [EdMa] in 1992 
assuming that the maps
$S_{\smalledge}$ were similarities
and was subsequently extended to the conformal case by Cole
[Col1,Col2]
building on earlier results due to 
Arbeiter \& Patzschke [ArPa] and Patzschke [Pa].
Below we use the following notation,
namely,  
 if $\varphi:\Bbb R\to\Bbb R$ is a function,
then
 $\varphi^{*}:\Bbb R\to[\-\infty,\infty]$
 denotes
the Legendre transform 
 of $\varphi$ defined
 by
$$ 
\varphi^{*}(x)
  =
  \inf_{y}
  (xy+\varphi(y))\,.
  $$
We can now state Theorem D.

\bigskip

\proclaim{Theorem D [Col1,Pa]}
Let $(\mu_{\smallvertexi})_{\smallvertexi\in\smallV}$
be the list of graph-directed self-conformal measures
defined by (2.7).
Let $\alpha\in\Bbb R$.
If the OSC is satisfied, then we have
 $$
 f_{\mu_{\smallsmallvertexi}}(\alpha)
 =
 \beta^{*}(\alpha)
 $$
for all $\vertexi\in\V$. 
\endproclaim

\bigskip

As a first application of 
Theorem 4.4, Corollary 4.5, Theorem 4.6 and 
Corollary 4.7 
we obtain a 
dynamical multifractal zeta-function
with an associated
 Bowen equation whose solution
equals the 
multifractal spectrum $f_{\mu_{\smallsmallvertexi}}(\alpha)$
of the graph-directed 
self-conformal measure $\mu_{\smallvertexi}$.
This is the content of the next theorem.

\bigskip

\proclaim{Theorem 5.1. 
Dynamical multifractal zeta-functinons for 
multifractal spectra of graph-directed
self-conformal measures}
Let $(\mu_{\smallvertexi})_{\smallvertexi\in\smallV}$
be the list of graph-directed self-conformal measures
associated with the list
$\big(\,
 \V,\,
 \E,\,
 (V_{\smallvertexi})_{\smallvertexi\in\smallV},\,
 (X_{\smallvertexi})_{\smallvertexi\in\smallV},\,
 (S_{\smalledge})_{\smalledge\in\smallE},\,
 (p_{\smalledge})_{\smalledge\in\smallE}\,
 \big)$, i\.e\.
$\mu_{\smallvertexi}$ is the unique probability measure such that
$ \mu_{\smallvertexi}
 =
 \sum_{\smalledge\in E_{\smallsmallvertexi}}\,
 p_{\smalledge}\,\mu_{\termi(\smalledge)}\circ S_{\smalledge}^{-1}$.

For $C\subseteq\Bbb R$ and 
a continuous function 
$\varphi:\Sigma_{\smallG}^{\Bbb N}\to\Bbb R$,   we
define the dynamical
graph-directed
self-conformal multifractal zeta-function by
  $$
  \zeta_{C}^{\sdyncon}(\varphi;z)
  =
  \sum_{n}
  \frac{z^{n}}{n}
  \left(
  \sum
  \Sb
  \bold i\in\Sigma_{\smallsmallG}^{n}\\
  {}\\
  \frac{\log p_{\bold i}}{\log\diam K_{\bold i}}
 \in  C
  \endSb
  \sup_{\bold u\in[\bold i]}
  \exp
  \sum_{k=0}^{n-1}
  \varphi S^{k}\bold u
  \right)
  $$
 for those complex numbers $z$ 
 for which the series converges.
Let $\Lambda$ be defined by (2.10) and 
let $\beta$ be defined by (5.1)
 (or, alternatively, by (5.2)).
Let $\alpha\in\Bbb R$.

 \roster

  \item"(1)"
 There is a unique real number 
 $\,\,\scr f\,\,(\alpha)$ such that
  $$
  \lim_{r\searrow0}
  \sigma_{\radius}
  \big(
  \,
  \zeta_{B(\alpha,r)}^{\sdyncon}(\,\,\scr f\,\,(\alpha)\,\Lambda;\cdot)
  \,
  \big)
  =
  1\,.
  $$

\item"(2)"
We have
 $$
 \scr f\,\,(\alpha)
 =
 \beta^{*}(\alpha)\,.
 $$

   \item"(3)"
 If the OSC is satisfied,
 then we have
 $$
\align
\quad
  \scr f\,\,(\alpha)
 &=
 f_{\mu_{\smallsmallvertexi}}(\alpha)
 =
 \dim_{\Haus}
 \Bigg\{
 x\in K_{\smallvertexi}
 \,\Bigg|\,
 \lim_{r\searrow 0}
  \frac{\log\mu_{\smallvertexi}(B(x,r))}{\log r}
  =
 \alpha
 \Bigg\}
 \endalign
 $$
for all $\vertexi\in\V$.

 \endroster

\endproclaim

\bigskip

We  will now prove Theorem 5.1.
Recall  that the function $\Lambda:\Sigma_{\smallG}^{\Bbb N}\to\Bbb R$ is defined by
$\Lambda(\bold i)=\log|DS_{\smalledge_{1}}(\pi S\bold i)|$ for 
$\bold i=\edge_{1}\edge_{2}\ldots\in\Sigma_{\smallG}^{\Bbb N}$.
Also,
recall that
$\Phi:\Sigma_{\smallG}^{\Bbb N}\to\Bbb R$ is defined by
$\Phi(\bold i)=\log p_{\ini(\bold i)}$ 
for
$\bold i\in\Sigma_{\smallG}^{\Bbb N}$.
We now introduce the following definition.
Define 
$U:\Cal P(\Sigma_{\smallG}^{\Bbb N})\to\Bbb R$
by 
 $$
 U\mu
  =
 \frac{\int\Phi\,d\mu}{\int\Lambda\,d\mu}\,,
 \tag5.3
 $$
and note that if $\bold i\in\Sigma_{\smallG}^{n}$, then
 $$
 UL_{n}[\bold i]
 =
 \Bigg\{
 \frac{\log p_{\bold i}}{\log |DS_{\bold i}(\pi\bold u)|}
 \,\Bigg|\,
 \bold u\in\Sigma_{\smallG}^{\Bbb N}
 \,\,
 \text{with}
 \,\,
 \termi(\bold i)=\ini(\bold u)
 \Bigg\}\,.
 $$
 It therefore follows that
  $$
  \zeta_{C}^{\dyn,U}(\varphi;z)
  =
  \sum_{n}
  \,\,
  \frac{z^{n}}{n}
  \,\,
  \left(
  \sum
  \Sb
   \bold i\in\Sigma_{\smallsmallG}^{n}\\
   {}\\
   \forall
 \bold u\in\Sigma_{\smallsmallG}^{\Bbb N}
 \,\,
 \text{with}
 \,\,
 \termi(\bold i)=\ini(\bold u)
 \,\,:\,\,
 \frac{\log p_{\bold i}}{\log |DS_{\bold i}(\pi\bold u)|}
\in
C
  \endSb
  \sup_{\bold u\in[\bold i]}
  \,\,
  \exp
  \sum_{k=0}^{n-1}
  \varphi S^{k}\bold u
  \right)
  \,.
  \tag5.4
  $$ 
In order to prove Theorem 5.1  we first prove
three  small auxiliary
results, namely, Proposition 5.2, Proposition 5.3 and Proposition 5.4.

\bigskip

\proclaim{Proposition 5.2}
Let $U$
be defined by (5.3).
Let $\Lambda$ be defined by (2.10)
and let
 $\beta$
 be defined by (5.1)
 (or, alternatively, by (5.2)).
Let $\alpha\in \Bbb R$.
Then
 $$
 \sup
 \Sb
 \mu\in\Cal P_{S}(\Sigma_{\smallsmallG}^{\Bbb N})\\
 {}\\
 U\mu=\alpha
 \endSb
 -
 \frac{h(\mu)}{\int\Lambda\,d\mu}
 =
 \beta^{*}(\alpha)\,.
$$
\endproclaim
\noindent{\it Proof}\newline

\noindent
This result in folk-lore. However, for the sake of completeness we have
 decided to include the brief proof.
We must the following two inequalities, namely
 $$
 \align
 \beta^{*}(\alpha)
&\le
 \sup
 \Sb
 \mu\in\Cal P_{S}(\Sigma_{\smallsmallG}^{\Bbb N})\\
 {}\\
 U\mu=\alpha
 \endSb
 -
 \frac{h(\mu)}{\int\Lambda\,d\mu}\,,
 \tag5.5\\
 \sup
 \Sb
 \mu\in\Cal P_{S}(\Sigma_{\smallsmallG}^{\Bbb N})\\
 {}\\
 U\mu=\alpha
 \endSb
 -
 \frac{h(\mu)}{\int\Lambda\,d\mu}
&\le
 \beta^{*}(\alpha)\,.
 \tag5.6
 \endalign
 $$

{\it Proof of (5.5).}
For $s\in\Bbb R$ and $q\in\Bbb R$,
let $\mu_{s,q}$ denote the
Gibbs state of $s\Lambda+q\Phi$.
We now prove the following two claims.

\medskip

{\it Claim 1. For all $q$, we have
$\frac
 {\int\Phi\,d\mu_{\beta(q),q}}
 {\int\Lambda\,d\mu_{\beta(q),q}} 
 =
 -\beta'(q)$.}

\noindent
{\it Proof of Claim 1.}
Define $F:\Bbb R\times\Bbb R\to\Bbb R$ by
$F(s,q)
 =
 P
 (
 \,
 s\Lambda+q\Phi
 \,
 )$
for $s,q\in\Bbb R$.
It follows from [Rue1]  that
$F$ is real analytic with
$\frac{\partial}{\partial s}F(s,q)
=
 \int\Lambda\,d\mu_{s,q}$
 and
$\frac{\partial}{\partial q}F(s,q)
= 
\int\Phi\,d\mu_{s,q} $.
Next, 
since
$0
 =
 F(\beta(q),q )$
for all $q$, 
it therefore follows from an application of the chain rule that
$$
0
=
\int\Lambda\,d\mu_{\beta(q),q}
\,\,
\beta'(q)
\,\,
+
\,\,
\int\Phi\,d\mu_{\beta(q),q} $$
for all $ q$,
whence
 $-\beta'(q)
 =
 \frac
 {\int\Phi\,d\mu_{\beta(q),q}}
 {\int\Lambda\,d\mu_{\beta(q),q}} $
for all $q$.
This completes the proof of Claim 1.

\medskip

{\it Claim 2. For all $q$, we have
$ -
 \frac{h(\mu_{\beta(q),q})}{\int\Lambda\,d\mu_{\beta(q),q}}
 \ge
 \beta^{*}(-\beta'(q))$.}

\noindent
{\it Proof of Claim 2.}
Since $\mu_{\beta(q),q}$ is a Gibbs state
of 
$\beta(q)\Lambda+q\Phi$
and
$P(\beta(q)\Lambda+q\Phi)=0$,
we deduce that
$$
\align
0
&=
P(\beta(q)\Lambda+q\Phi)\\
&=
h(\mu_{\beta(q),q})
+
\int(\beta(q)\Lambda+q\Phi)\,d\mu_{\beta(q),q}\\
&=
h(\mu_{\beta(q),q})
+
\beta(q)\int\Lambda\,d\mu_{\beta(q),q}
+
q\int\Phi\,d\mu_{\beta(q),q}\,,
\endalign
$$ 
whence
$-
 \frac{h(\mu_{\beta(q),q})}{\int\Lambda\,d\mu_{\beta(q),q}}
 =
 \beta(q)
 +
 q
 \frac{\int\Phi\,d\mu_{\beta(q),q}}{\int\Lambda\,d\mu_{\beta(q),q}}$.
It follows from this and Claim 1 that
$$-
 \frac{h(\mu_{\beta(q),q})}{\int\Lambda\,d\mu_{\beta(q),q}}
=
 \beta(q)
 +
 q
 \frac{\int\Phi\,d\mu_{\beta(q),q}}{\int\Lambda\,d\mu_{\beta(q),q}}
 =
 \beta(q)
 -
 q
 \beta'(q)
 \ge
 \inf_{p} 
 (
  \beta(p)
-
p
\beta'(q)
 )
 =
 \beta^{*}(-\beta'(q)) 
 $$
 for all $q$.
This completes the proof of Claim 2.
 
 \medskip

 We can now prove (5.5).
 If $ \beta^{*}\alpha)=-\infty$, then 
 inequality (5.5)
is clear. 
Hence, we may assume that
$ \beta^{*}(\alpha)>-\infty$.
In this case it follows from the convexity of $\beta$
that there is a point $q_{\alpha}\in\Bbb R$
such that
$\alpha=-\beta'(q_{\alpha})$, see [Ro].
We therefore conclude from Claim 1 that the measure 
$\mu_{\beta(q_{\alpha}),q_{\alpha}}$
satisfies
$U\mu_{\beta(q_{\alpha}),q_{\alpha}}
=
\frac
 {\int\Phi\,d\mu_{\beta(\bold q_{\alpha}),q_{\alpha}}}
 {\int\Lambda\,d\mu_{\beta(\bold q_{\alpha}),\bold q_{\alpha}}} 
 =
 -\beta'(q_{\alpha})
=
\alpha$, whence, using Claim 2,
$\sup_{
 \mu\in\Cal P_{S}(\Sigma_{\smallsmallG}^{\Bbb N})\,,\,
 U\mu=\alpha
 }
 -
 \frac{h(\mu)}{\int\Lambda\,d\mu}
 \ge
 -
 \frac{h(\mu_{\beta(q_{\alpha}),q_{\alpha}})}
 {\int\Lambda\,d\mu_{\beta(q_{\alpha}),q_{\alpha}}}
 \ge
 \beta^{*}(-\beta'(q_{\alpha}))
 =
 \beta^{*}(\alpha)\,.$
This completes the proof of (5.5).

{\it Proof of (5.6).}
Fix $q\in\Bbb R$ and $\mu\in\Cal P_{S}(\Sigma_{\smallG}^{\Bbb N})$ with
$U\mu=\alpha$.
Using the variational principle (see [Wa]) we  conclude  that
$$
\align
P(\beta(q)\Lambda+q\Phi)
&=
\sup_{\nu\in\Cal P(\Sigma_{\smallsmallG}^{\Bbb N})}
(\,h(\nu)+\int(\beta(q)\Lambda+q\Phi)\,d\nu\,)\\
&\ge
h(\mu)+\int(\beta(q)\Lambda+q\Phi)\,d\mu\\
&=
h(\mu)+\beta(q)\int\Lambda\,d\mu+q\int\Phi\,d\mu\,.
\endalign
$$
Since $\Lambda<0$, this implies that
$\beta(q)
\ge
-
 \frac{h(\mu)}{\int\Lambda\,d\mu}
 -
 q
  \frac{\int\Phi\,d\mu}{\int\Lambda\,d\mu}$.
Next, 
as
$\frac{\int\Phi\,d\mu}{\int\Lambda\,d\mu}=U\mu=\alpha$, we therefore conclude  that
 $$
\align
q\alpha+\beta(q)
&\ge
 q\alpha 
-
 \frac{h(\mu)}{\int\Lambda\,d\mu}
 -
 q
  \frac{\int\Phi\,d\mu}{\int\Lambda\,d\mu}\\
 &\ge
 q\alpha 
 +
 \frac{h(\mu)}{\int\Lambda\,d\mu}
 -
 q\alpha\\
&=
-\frac{h(\mu)}{\int\Lambda\,d\mu}\,.
\endalign
$$ 
Finally,
taking infimum over all $q\in\Bbb R$ 
and 
taking
supremum over all
$\mu\in\Cal P_{S}(\Sigma_{\smallG}^{\Bbb N})$
with
$U\mu=\alpha$ 
gives
inequality (5.6).
\hfill$\square$

 \bigskip

\proclaim{Proposition 5.3}
Let $U$
be defined by (5.3).
Let $\Lambda$ be defined by (2.10)
and let
 $\beta$
 be defined by (5.1)
 (or, alternatively, by (5.2)).
Let $\alpha\in\Bbb R$.
Let $t$ be the unique real number such that
 $$
  \lim_{r\searrow0}
 \sigma_{\radius}
 \big(
 \,
 \zeta_{B(\alpha,r)}^{\dyn,U}(t\Lambda;\cdot)
 \,
 \big)
 =
 1\,.
 $$
Then 
 $$
 \sup
 \Sb
 \mu\in\Cal P_{S}(\Sigma_{\smallsmallG}^{\Bbb N})\\
 {}\\
 U\mu=\alpha
 \endSb
 \Bigg(
 h(\mu)+t\int\Lambda\,d\mu
 \Bigg)
=
0\,,
$$
and we have
 $$
 t
 =
 \sup
 \Sb
 \mu\in\Cal P_{S}(\Sigma_{\smallsmallG}^{\Bbb N})\\
 {}\\
 U\mu=\alpha
 \endSb
 -
 \frac{h(\mu)}{\int\Lambda\,d\mu}
 =
 \beta^{*}(\alpha)\,.
$$
\endproclaim
\noindent{\it Proof}\newline

\noindent
We first note that
it follows immediately from
Theorem 4.4 that
 $$
 \sup
 \Sb
 \mu\in\Cal P_{S}(\Sigma_{\smallsmallG}^{\Bbb N})\\
 {}\\
 U\mu=\alpha
 \endSb
 \Bigg(
 h(\mu)+t\int\Lambda\,d\mu
 \Bigg)
 =
   \lim_{r\searrow0}
 \overline P_{B(\alpha,r)}^{U}(t\Lambda) 
 =
   \lim_{r\searrow0}
 -\log
  \sigma_{\radius}
 \big(
 \,
 \zeta_{B(\alpha,r)}^{\dyn,U}(t\Lambda;\cdot)
 \,
 \big)
 =
 0\,.
 \tag5.7
 $$
Equality (5.7) is easily seen to imply that
 $$
 t
 =
 \sup
 \Sb
 \mu\in\Cal P_{S}(\Sigma_{\smallsmallG}^{\Bbb N})\\
 {}\\
 U\mu=\alpha
 \endSb
 -
 \frac{h(\mu)}{\int\Lambda\,d\mu}\,.
 $$
 Finally,
it follows from the above equality
and Proposition 5.2 that
$t=\beta^{*}(\alpha)$.
This completes the proof.
\hfill$\square$

\bigskip

\proclaim{Proposition 5.4}
Let $U$
be defined by (5.3).
Fix a continuous function
$\varphi:\Sigma_{\smallG}^{\Bbb N}\to\Bbb R$.
\roster
\item"(1)"
There is a sequence $(\Delta_{n})_{n}$ with $\Delta_{n}>0$ 
and 
$\Delta_{n}\to 0$ such that
for all
closed subsets $W$ of $\Bbb R$
and
for
all
$n\in\Bbb N$, $\bold i\in\Sigma_{\smallG}^{n}$ and 
$\bold u\in\Sigma_{\smallG}^{\Bbb N}$
with
 $\termi(\bold i)=\ini(\bold u)$, we have
 $$
 \align
 \dist
 \Bigg(
 \,
\frac{\log p_{\bold i}}{\log |DS_{\bold i}(\pi\bold u)|}
 \,,\,
 W
 \,
 \Bigg)
&\le
\dist
\Bigg(
 \,
 \frac{\log p_{\bold i}}{\log \diam K_{\bold i}}
  \,,\,
 W
 \,
 \Bigg)
\,\,
\,\,\,\,\,\,
+
\,\,
\Delta_{n}\,,
\tag5.8\\
&{}\\
\dist
\Bigg(
 \,
\frac{\log p_{\bold i}}{\log \diam K_{\bold i}}
  \,,\,
 W
 \,
 \Bigg)
&\le
\dist
 \Bigg(
  \,
  \frac{\log p_{\bold i}}{\log |DS_{\bold i}(\pi\bold u)|}
 \,,\,
 W
 \,
 \Bigg)
 \,\,
+
\,\,
\Delta_{n}\,.
\tag5.9
\endalign
$$

\item"(2)"
Let $W$ be a closed subset 
of $\Bbb R$.
For all $r>0$, we have
 $$
 \align
 \sigma_{\radius}
 \big(
 \,
 \zeta_{B(W,r)}^{\dyn,U}(\varphi;\cdot)
 \,
 \big)
&\le
 \sigma_{\radius}
 \big(
 \,
 \zeta_{W}^{\dyn\text{-}\scon}(\varphi;\cdot)
 \,
 \big)\,,
 \tag5.10\\
 \sigma_{\radius}
 \big(
 \,
 \zeta_{B(W,r)}^{\dyn\text{-}\scon}(\varphi;\cdot)
 \,
 \big)
&\le
\sigma_{\radius}
 \big(
 \,
 \zeta_{W}^{\dyn,U}(\varphi;\cdot)
 \,
 \big)\,.
 \tag5.11
 \endalign
 $$

\item"(3)"
Let $C$ be a closed subset of $\Bbb R$.
Then we have
 $$
 \lim_{r\searrow0}
 \sigma_{\radius}
 \big(
 \,
 \zeta_{B(C,r)}^{\dyn\text{-}\scon}(\varphi;\cdot) 
 \,
 \big)
 =
 \lim_{r\searrow0}
 \sigma_{\radius}
 \big(
 \,
 \zeta_{B(C,r)}^{\dyn,U}(\varphi;\cdot)
 \,
 \big)\,.
$$

\endroster
\endproclaim
\noindent{\it Proof}\newline
\noindent
(1)
It is well-known and follows from the 
Principle of Bounded Distortion 
(see, for example, [Bar,Fa])
that there is a constant $c>0$ such that
for all integers $n$
and all $\bold i\in\Sigma_{\smallG}^{n}$
and all $\bold u,\bold v\in[\bold i]$,
 we have
$\frac{1}{c}
\le
\frac{|DS_{\bold i}(\pi S^{n}\bold u)|}{\diam K_{\bold i}}
\le
c$
and
$\frac{1}{c}
\le
\frac{|DS_{\bold i}(\pi S^{n}\bold u)|}{|DS_{\bold i}(\pi S^{n}\bold v)|}
\le
c$.
It is not difficult to see that the desired result follows from this.

\noindent
(2)
Fix $r>0$.
Let $(\Delta_{n})_{n}$ be the sequence from Part (1).
Since $\Delta_{n}\to 0$, we can find
a positive integer
$N_{r}$
such that if $n\ge N_{r}$, then $\Delta_{n}< r$.
Consequently, using (5.9) in Part (1),
for all $n\ge N_{r}$, we have
$$
\align
\sum
 \Sb
\bold i\in\Sigma_{\smallsmallG}^{n}\\
 {}\\
UL_{n}[\bold i]\subseteq W
\endSb
  \sup_{\bold u\in[\bold i]}
  \exp
  \sum_{k=0}^{n-1}
  \varphi S^{k}\bold u
 & 
 =
 \qquad
 \sum
 \Sb
 \bold i\in\Sigma_{\smallsmallG}^{n}\\
 {}\\
 \forall
 \bold u\in\Sigma^{\Bbb N}
 \,\,
 \text{with}
 \,\,
 \termi(\bold i)
 =
 \ini(\bold u)
 \,\,:\,\,
 \frac{\log p_{\bold i}}{\log |DS_{\bold i}(\pi\bold u)|}
\in
W
 \endSb
 \qquad
  \sup_{\bold u\in[\bold i]}
  \exp
  \sum_{k=0}^{n-1}
  \varphi S^{k}\bold u\\
 & 
   =
   \sum
 \Sb
\bold i\in\Sigma_{\smallsmallG}^{n}\\
 {}\\
 \forall
 \bold u\in\Sigma^{\Bbb N}
 \,\,
 \text{with}
 \,\,
 \termi(\bold i)
 =
 \ini(\bold u)
 \,\,:\,\,
 \dist
 \big(
 \,
 \frac{\log p_{\bold i}}{\log |DS_{\bold i}(\pi\bold u)|}
  \,,\,
 W
 \,
 \big)
 \,
 =
 0
 \endSb
 \,
  \sup_{\bold u\in[\bold i]}
  \exp
  \sum_{k=0}^{n-1}
  \varphi S^{k}\bold u\\
   &
   \le
   \qquad\quad\,\,\,\,
   \sum
  \Sb
\bold i\in\Sigma_{\smallsmallG}^{n}\\
 {}\\
 \dist
 \big(
  \,
 \frac{\log p_{\bold i}}{\log \diam K_{\bold i}}
 \,,\,
 W
 \,
 \big)
 \,
 \le
 \,
 0
 +
 \Delta_{r}
 \endSb
 \qquad\qquad
  \sup_{\bold u\in[\bold i]}
  \exp
  \sum_{k=0}^{n-1}
  \varphi S^{k}\bold u\\
   &
      \le
   \qquad\qquad\,\,
   \sum
  \Sb
\bold i\in\Sigma_{\smallsmallG}^{n}\\
 {}\\
 \dist
 \big(
 \,
 \frac{\log p_{\bold i}}{\log \diam K_{\bold i}}
 \,,\,
 W
 \,
 \big)
 \,
 <
 \,
 r 
 \endSb
 \qquad\qquad
 \,\,\,\,\,\,
  \sup_{\bold u\in[\bold i]}
  \exp
  \sum_{k=0}^{n-1}
  \varphi S^{k}\bold u\\
   & 
   =
   \qquad\qquad
   \quad\,\,
   \sum
  \Sb
\bold i\in\Sigma_{\smallsmallG}^{n}\\
 {}\\
 \frac{\log p_{\bold i}}{\log \diam K_{\bold i}}
\in
B(W,r)
 \endSb
 \qquad\qquad
 \qquad
  \sup_{\bold u\in[\bold i]}
  \exp
  \sum_{k=0}^{n-1}
  \varphi S^{k}\bold u\,.  \\
  \tag5.12
  \endalign
  $$
A similar argument using (5.8) in Part (1) shows that
 $$
 \sum
  \Sb
\bold i\in\Sigma_{\smallsmallG}^{n}\\
 {}\\
  \frac{\log p_{\bold i}}{\log \diam K_{i}}
\in
W
 \endSb
  \sup_{\bold u\in[\bold i]}
  \exp
  \sum_{k=0}^{n-1}
  \varphi S^{k}\bold u
  \le
  \sum
 \Sb
\bold i\in\Sigma_{\smallsmallG}^{n}\\
 {}\\
UL_{n}[\bold i]\subseteq B(W,r)
\endSb
  \sup_{\bold u\in[\bold i]}
  \exp
  \sum_{k=0}^{n-1}
  \varphi S^{k}\bold u\,.
  \tag5.13
   $$
The desired results follow immediately from inequalities
(5.12) and (5.13).

\noindent
(3)
This result follows easily from Part (2).
\hfill$\square$

\bigskip

We can now prove Theorem 5.1.

\bigskip

\noindent{\it Proof of Theorem 5.1}\newline
Let $U$ be defined by (5.3).

\noindent
(1)
It follows from Proposition 5.4 that
 $$
  \lim_{r\searrow0}
 \sigma_{\radius}
 \big(
 \,
 \zeta_{B(C,r)}^{\dyn\text{-}\scon}(t\Lambda;\cdot) 
 \,
 \big)
 =
 \lim_{r\searrow0}
 \sigma_{\radius}
 \big(
 \,
 \zeta_{B(C,r)}^{\dyn,U}(t\Lambda;\cdot)
 \,
 \big)
 \tag5.14
$$
for all $t$.
Also, it follows from Proposition 5.2 that
there is a unique number
 $\,\,\scr f\,\,(\alpha)$ such that
  $$
  \lim_{r\searrow0}
  \sigma_{\radius}
  \big(
  \,
  \zeta_{B(\alpha,r)}^{\dyn,U}(\,\,\scr f\,\,(\alpha)\,\Lambda;\cdot)
  \,
  \big)
  =
  1\,.
  \tag5.15
  $$
Combining (5.14) and (5.15) shows that
 $\,\,\scr f\,\,(\alpha)$ is the unique real number such we have
 the following
$\lim_{r\searrow0}
  \sigma_{\radius}
  \big(
  \,
  \zeta_{B(\alpha,r)}^{\dyn\text{-}\scon}(\,\,\scr f\,\,(\alpha)\,\Lambda;\cdot)
  \,
  \big)
  =
  \lim_{r\searrow0}
  \sigma_{\radius}
  \big(
  \,
  \zeta_{B(\alpha,r)}^{\dyn,U}(\,\,\scr f\,\,(\alpha)\,\Lambda;\cdot)
  \,
  \big)
  =
  1$.

  \noindent
  (2)
 It follows from (1) and Proposition 5.4 that
 $$
 \align
 1
&=
  \lim_{r\searrow0}
 \sigma_{\radius}
 \big(
 \,
 \zeta_{B(C,r)}^{\dyn\text{-}\scon}(\,\,\scr f\,\,(\alpha)\,\Lambda;\cdot) 
 \,
 \big)
 =
 \lim_{r\searrow0}
 \sigma_{\radius}
 \big(
 \,
 \zeta_{B(C,r)}^{\dyn,U}(\,\,\scr f\,\,(\alpha)\,\Lambda;\cdot)
 \,
 \big)\,,
\endalign
$$
 and Corollary 4.5 therefore implies that
 the number  $\,\,\scr f\,\,(\alpha)$ 
 is given by
 $$
  \scr f\,\,(\alpha)
  =
 \sup
 \Sb
 \mu\in\Cal P_{S}(\Sigma_{\smallsmallG}^{\Bbb N})\\
 {}\\
 U\mu=\alpha
 \endSb
 -
 \frac{h(\mu)}{\int\Lambda\,d\mu}\,.
 \tag5.16
 $$ 
Finally, combining Proposition 5.2 
and (5.16) shows that
$\,\,\scr f\,\,(\alpha)=\beta^{*}(\alpha)$.

\noindent
(3)
This follows immediately from (2) and Theorem D.
\hfill$\square$

\bigskip

{\bf 5.2. Multifractal spectra of 
ergodic Birkhoff averages.}
We first fix $\gamma\in(0,1)$ and define  the metric
 $\distance_{\gamma}$ on $\Sigma_{\smallG}^{\Bbb N}$
 as follows.
For $\bold i,\bold j\in\Sigma_{\smallG}^{\Bbb N}$ 
with $\bold i\not=\bold j$,
we will write
$\bold i\wedge\bold j$ for the longest common 
prefix of 
$\bold i$ and $\bold j$
(i\.e\.
$
\bold i\wedge\bold j
=
\bold u
$
where $\bold u$ is the unique element in $\Sigma_{\smallG}^{*}$
for which there 
are $\bold k',\bold k''\in\Sigma_{\smallG}^{\Bbb N}$
with
$\bold k'=\edge'_{1}\edge'_{2}\ldots$
and
$\bold k''=\edge''_{1}\edge''_{2}\ldots$
such that
$\edge'_{1}
 \not= \edge''_{1}$,
 $\bold i
 =
 \bold u\bold k'$
and 
$\bold j
 =
 \bold u\bold k''$). 
The metric
$\distance_{\gamma}$ is now defined by
 $$
 \distance_{\gamma}(\bold i,\bold j)
 =
 \cases
 0
&\quad 
 \text{if $\bold i=\bold j$;}\\
 \gamma^{|\bold i\wedge\bold j|}
&\quad 
 \text{if $\bold i\not=\bold j$,}
 \endcases
 $$
for $\bold i,\bold j\in\Sigma_{\smallG}^{\Bbb N}$;
throughout this section, we equip
 $\Sigma_{\smallG}^{\Bbb N}$
with the metric  $\distance_{\gamma}$,
and continuity and Lipschitz properties of functions $f:\Sigma_{\smallG}^{\Bbb N}\to\Bbb R$
from $\Sigma_{\smallG}^{\Bbb N}$ to $\Bbb R$
will  always
refer to
the metric  $\distance_{\gamma}$.
Multifractal analysis of Birkhoff
averages has received significant interest
during
the past 10 years, see, for example,
[BaMe,FaFe,FaFeWu,FeLaWu,Oli,Ol1,OlWi].
Fix a positive integer $M$.
The 
multifractal
spectrum
$F_{\bold f}^{\erg}$
of ergodic Birkhoff averages of a vector valued 
continuous function
$\bold f:\Sigma_{\smallG}^{\Bbb N}\to\Bbb R^{M}$ is defined by
 $$
F_{\bold f}^{\erg}(\pmb\alpha)
=
 \dim_{\Haus}
 \pi
 \Bigg\{
 \bold i\in\Sigma_{\smallG}^{\Bbb N}
 \,\Bigg|\,
 \lim_{n}
 \frac{1}{n}\sum_{k=0}^{n-1}\bold f(S^{k}\bold i)
 =
 \pmb\alpha
 \Bigg\}
$$
for $\pmb\alpha\in\Bbb R^{M}$;
recall, that the map $\pi$ is defined in Section 2.
One of the main problems
in
multifractal analysis of Birkhoff
averages
is the detailed study of the multifractal
spectrum
$F_{\bold f}^{\erg}$.
For example,
Theorem E below is
proved
in different settings and at various levels of generality
in [FaFe,FaFeWu,FeLaWu,Oli,Ol1,OlWi].
Before we can state
this result
we introduce the following 
notation.
If $(x_{n})_{n}$ is a sequence of points in a metric space $X$, then we write
$\acc_{n} x_{n}$ for the set of accumulation points 
of the sequence  $(x_{n})_{n}$, i\.e\.
 $$
 \underset{n}\to{\acc} \,x_{n}
=
\Big\{
x\in X
\,\Big|\,
\text{
$x$ is an accumulation point of $(x_{n})_{n}$
}
\Big\}\,.
$$
We will also use the following notation. 
Namely, if 
$\bold f:\Sigma_{\smallG}^{\Bbb N}\to\Bbb R^{M}$
is a continuous function with
$\bold f=(f_{1},\ldots, f_{M})$, then we will write
 $$
 \int\bold f\,d\mu
 =
 \Bigg(
 \int f_{1}\,d\mu
 \,,\,
 \ldots
 \,,\,
 \int f_{M}\,d\mu
 \Bigg)
 $$
for $\mu\in\Cal P(\Sigma_{\smallG}^{\Bbb N})$.
We can now state Theorem E.

\bigskip

\proclaim{Theorem E [FaFe,FaFeWu,FeLaWu,Oli,Ol1,OlWi]}
Fix $\gamma\in(0,1)$
and
let $\bold f:\Sigma_{\smallG}^{\Bbb N}\to\Bbb R^{M}$ be a Lipschitz function
with respect to the metric $\distance_{\gamma}$.
Let $\Lambda:\Sigma_{\smallG}^{\Bbb N}\to\Bbb R$ be defined by (2.10).
Let $C$ be a closed subset of $\Bbb R^{M}$.
If the OSC is satisfied, then
 $$
 \dim_{\Haus}
 \pi
 \Bigg\{
 \bold i\in\Sigma_{\smallG}^{\Bbb N}
 \,\Bigg|\,
 \,\underset{n}\to\acc
 \frac{1}{n}\sum_{k=0}^{n-1}\bold f(S^{k}\bold i)
\subseteq
C
 \Bigg\}
 =
 \sup
 \Sb
 \mu\in\Cal P_{S}(\Sigma_{\smallsmallG}^{\Bbb N})\\
 \int \bold f\,d\mu\in C
 \endSb
 -
 \frac{h(\mu)}{\int\Lambda\,d\mu}\,. 
 $$
In particular, if the OSC is satisfied and $\pmb\alpha\in\Bbb R^{M}$, then we have 
 $$
  \dim_{\Haus}
 \pi
 \Bigg\{
 \bold i\in\Sigma_{\smallG}^{\Bbb N} 
 \,\Bigg|\,
 \lim_{n}
 \frac{1}{n}\sum_{k=0}^{n-1}\bold f(S^{k}\bold i)
 =
 \pmb\alpha
 \Bigg\}
 =
 \sup
 \Sb
 \mu\in\Cal P_{S}(\Sigma_{\smallsmallG}^{\Bbb N})\\
 \int \bold f\,d\mu
 =
 \pmb\alpha
 \endSb
 -
 \frac{h(\mu)}{\int\Lambda\,d\mu}\,.
 $$ 
\endproclaim

\bigskip

As a second
 application of 
Theorem 4.4, Corollary 4.5, Theorem 4.6 and 
Corollary 4.7 
we will now obtain
zeta-functions
whose
radii of convergence equal
different types
of 
multifractal
spectra
of ergodic Birkhoff averages.
We first 
state and prove a 
rather
general 
result, namely Theorem 5.5 below,
from which 
analogous results
for
a number
of
different
multifractal 
spectra
of ergodic Birkhoff averages, including
$F_{\bold f}^{\erg}(\pmb\alpha)$,
can be deduced.
Indeed, 
immediately
after the statement and proof of Theorem 5.5,
we will apply
 Theorem 5.5
 to 
 prove the following results, namely:
Theorem 5.6
on the 
multifractal spectra of
ergodic averages of continuous vector valued functions,
Theorem 5.7
on the
multifractal spectra of
relative
ergodic averages of continuous  functions,
and finally 
Theorem 5.8
on the
multifractal spectra of
a more general type of
relative
ergodic averages of continuous  functions.

\bigskip

\proclaim{Theorem 5.5.
Multifractal zeta-functions for
abstract 
multifractal spectra of ergodic Birkhoff averages}
Fix $\gamma\in(0,1)$
and
$W\subseteq \Bbb R^{I}$
and
let
$\pmb\Phi:\Sigma_{\smallG}^{\Bbb N}\to\Bbb R^{I}$ be a Lipschitz function
with respect to the metric $\distance_{\gamma}$
such that
 $\{
 \int
 \pmb\Phi\,d\mu
 \,|\,
 \mu
 \in\Cal P(\Sigma_{\smallG}^{\Bbb N})
 \}
 \subseteq
 W$;
recall, that if $\pmb\Phi=(\Phi_{1},\ldots,\Phi_{I})$, then we write
$\int\pmb\Phi\,d\mu
=
(\int\Phi_{1}\,d\mu,\ldots,\int\Phi_{I}\,d\mu)$.
Let
 $Q:W\to \Bbb R^{M}$
be a continuous function.

For $C\subseteq \Bbb R^{M}$ and 
a continuous function 
$\varphi:\Sigma_{\smallG}^{\Bbb N}\to\Bbb R$,   we
define the 
abstract dynamical
ergodic multifractal zeta-function associated with $Q$
by
  $$
  \zeta_{C}^{\dyn\text{-}\erg}(\varphi;z)
  =
  \sum_{n}
  \frac{z^{n}}{n}
  \left(
  \sum
  \Sb
 \bold i\in\Sigma_{\smallsmallG}^{n}\\
  {}\\
  \forall\bold u\in[\bold i]
  \,\,:\,\,
  Q\big(
 \frac{1}{n}\sum_{k=0}^{n-1}\pmb\Phi(S^{k}\bold u)
 \big)
 \in  C
  \endSb
  \sup_{\bold u\in[\bold i]}
  \exp
  \sum_{k=0}^{n-1}
  \varphi S^{k}\bold u
  \right)\,.
  $$
Let 
 $\Lambda:\Sigma_{\smallG}^{\Bbb N}\to\Bbb R$ be defined by (2.10).

\noindent
 {\rm (1)} Assume 
 that $C\subseteq \Bbb R^{M}$ is closed.

 \roster
  \item"(1.1)"
 There is a unique real number 
 $\,\,\scr f\,\,(C)$ such that
  $$
  \lim_{r\searrow0}
  \sigma_{\radius}
  \big(
  \,
  \zeta_{B(C,r)}^{\dyn\text{-}\erg}(\,\,\scr f\,\,(C)\,\Lambda;\cdot)
  \,
  \big)
  =
  1\,.
  $$

\item"(1.2)"
We have
$$
\,\,\scr f\,\,(C)
 =
 \sup_{\pmb\alpha\in C}
 \,
 \sup
 \Sb
 \mu\in\Cal P_{S}(\Sigma_{\smallsmallG}^{\Bbb N})\\
Q\big( \int\pmb\Phi\,d\mu\big)=\pmb\alpha
 \endSb
 -
 \frac{h(\mu)}{\int\Lambda\,d\mu}\,.
 $$

\item"(1.3)"
If the OSC is satisfied, then we have
$$
\align
\quad
\,\,\scr f\,\,(C)
 &=
 \dim_{\Haus}
 \pi
 \Bigg\{
 \bold i\in\Sigma_{\smallG}^{\Bbb N}
 \,\Bigg|\,
 \,\underset{n}\to\acc
 \,\,
 Q
 \Bigg(
 \frac{1}{n}\sum_{k=0}^{n-1}\pmb\Phi(S^{k}\bold i)
 \Bigg)
\subseteq
C
 \Bigg\}\,.\\
 \endalign
 $$
\endroster

 \bigskip
 
 \noindent
 {\rm (2)} 
 Assume that there 
 are
 continuous and affine functions
 $\Gamma:W\to\Bbb R^{M}$
 and 
 $\Delta:W\to\Bbb R$
 with
 $\Delta(\bold x)\not=0$
 for all $\bold x$
 such that
 $Q=\frac{\Gamma}{\Delta}$,
 and
 assume that $C\subseteq\Bbb R^{M}$ is  closed 
and
convex 
 with
 $\overset{\,\circ}\to{C}
  \cap
  \{
  \frac{\Gamma(\int\pmb\Phi\,d\mu)}{\Delta(\int\pmb\Phi\,d\mu)}
  \,|\,
  \mu\in\Cal P_{S}(\Sigma_{\smallG}^{\Bbb N})
  \}
    \not=
  \varnothing$.

\roster

 \item"(2.1)"
 There is a unique real number 
 $\scr F\,\,(C)$ such that
  $$
  \sigma_{\radius}
  \big(
  \,
  \zeta_{C}^{\dyn\text{-}\erg}(\scr F\,\,(C)\,\Lambda;\cdot)
  \,
  \big)
  =
  1\,.
  $$

\item"(2.2)"
We have
 $$
 \scr F\,\,(C)
 =
 \sup_{\pmb\alpha\in C}
 \,
 \sup
 \Sb
 \mu\in\Cal P_{S}(\Sigma_{\smallsmallG}^{\Bbb N})\\
Q\big(\int \pmb\Phi\,d\mu\big)=\pmb\alpha
 \endSb
 -
 \frac{h(\mu)}{\int\Lambda\,d\mu}\,.
 $$

\item"(2.3)"
If the OSC is satisfied, then we have
$$
\align
\quad
 \scr F\,\,(C)
 &=
 \dim_{\Haus}
 \pi
 \Bigg\{
 \bold i\in\Sigma_{\smallG}^{\Bbb N}
 \,\Bigg|\,
 \,\underset{n}\to\acc
 \,\,
 Q
 \Bigg(
 \frac{1}{n}\sum_{k=0}^{n-1}\pmb\Phi(S^{k}\bold i)
 \Bigg)
\subseteq
C
 \Bigg\}\,.\\
 \endalign
 $$

\endroster

\endproclaim
\noindent{\it Proof}\newline
We first note that it follows from [Ol1] 
that if $C$ is  
a
closed  subset of $\Bbb R^{M}$, then
 $$
 \dim_{\Haus}
 \pi
 \Bigg\{
 \bold i\in\Sigma_{\smallG}^{\Bbb N}
 \,\Bigg|\,
 \,\underset{n}\to\acc
 \,\,
 Q
 \Bigg(
 \frac{1}{n}\sum_{k=0}^{n-1}\pmb\Phi(S^{k}\bold i)
 \Bigg)
\subseteq
C
 \Bigg\}
 =
  \sup_{\pmb\alpha\in C}
 \,
 \sup
 \Sb
 \mu\in\Cal P_{S}(\Sigma_{\smallsmallG}^{\Bbb N})\\
Q\big(\int \pmb\Phi\,d\mu\big)=\pmb\alpha
 \endSb
 -
 \frac{h(\mu)}{\int\Lambda\,d\mu}\,.
 \tag5.17
 $$
Next,
we
define the function $U:\Cal P(\Sigma_{\smallG}^{\Bbb N})\to\Bbb R^{M}$ 
to be the composite
of the following 2 maps, namely,
$$
 \Cal P(\Sigma_{\smallG}^{\Bbb N})
 @>\,\,\mu\to\int\pmb\Phi\,d\mu\,\,>>
 W
\,,\,\,\,\,
 W
 @>\,\,Q\,\,>>
 \Bbb R^{M}\,,
$$
i\.e\.
  $$
  U\mu
  =
  Q\big(\,{\ssize\int}\pmb\Phi\,d\mu\,\big)
  $$
for $\mu\in\Cal P(\Sigma_{\smallG}^{\Bbb N})$.  
Since clearly
$ \zeta_{C}^{\dyn\text{-}\erg}(\varphi;\cdot)
=
 \zeta_{C}^{\dyn,U}(\varphi;\cdot)$,
 the results now follow from Corollary 4.5,
 Corollary 4.7
 and (5.17).
\hfill$\square$

\bigskip

\newpage

\noindent
Next, we present three 
corollaries of Theorem 5.5.

\bigskip

\proclaim{Theorem 5.6.
Multifractal zeta-functinons for
multifractal spectra of vector valued ergodic Birkhoff averages}
Fix $\gamma\in(0,1)$
and
let $\bold f:\Sigma_{\smallG}^{\Bbb N}\to\Bbb R^{M}$ be a Lipschitz function
with respect to the metric $\distance_{\gamma}$.
For $C\subseteq\Bbb R^{M}$ and 
a continuous function 
$\varphi:\Sigma_{\smallG}^{\Bbb N}\to\Bbb R$,   we
define the dynamical
ergodic multifractal zeta-function by
  $$
  \zeta_{C}^{\dyn\text{-}\vector\text{-}\erg}(\varphi;z)
  =
  \sum_{n}
  \frac{z^{n}}{n}
  \left(
  \sum
  \Sb
 \bold i\in\Sigma_{\smallsmallG}^{n}\\
  {}\\
  \forall\bold u\in[\bold i]
  \,\,:\,\,
 \frac{1}{n}\sum_{k=0}^{n-1}\bold f(S^{k}\bold u)
 \in  C
  \endSb
  \sup_{\bold u\in[\bold i]}
  \exp
  \sum_{k=0}^{n-1}
  \varphi S^{k}\bold u
  \right)\,.
  $$
Let 
 $\Lambda:\Sigma_{\smallG}^{\Bbb N}\to\Bbb R$ be defined by (2.10).

\noindent
 {\rm (1)} Assume 
 that $C\subseteq\Bbb R^{M}$ is closed.

 \roster
  \item"(1.1)"
 There is a unique real number 
 $\,\,\scr f\,\,(C)$ such that
  $$
  \lim_{r\searrow0}
  \sigma_{\radius}
  \big(
  \,
  \zeta_{B(C,r)}^{\dyn\text{-}\vector\text{-}\erg}(\,\,\scr f\,\,(C)\,\Lambda;\cdot)
  \,
  \big)
  =
  1\,.
  $$
If $\pmb\alpha\in\Bbb R^{M}$ and 
$C=\{\pmb\alpha\}$, then we will write
  $\,\,\scr f\,\,(\pmb\alpha)
  =
 \,\,\scr f\,\,(C)$.

\item"(1.2)"
We have
$$
\,\,\scr f\,\,(C)
 =
 \sup_{\pmb\alpha\in C}
 \,
 \sup
 \Sb
 \mu\in\Cal P_{S}(\Sigma_{\smallsmallG}^{\Bbb N})\\
 \int\bold f\,d\mu=\pmb\alpha
 \endSb
 -
 \frac{h(\mu)}{\int\Lambda\,d\mu}\,.
 $$

\item"(1.3)"
If the OSC is satisfied, then we have
$$
\align
\quad
\,\,\scr f\,\,(C)
 &=
 \dim_{\Haus}
 \pi
 \Bigg\{
 \bold i\in\Sigma_{\smallG}^{\Bbb N}
 \,\Bigg|\,
 \,\underset{n}\to\acc
 \frac{1}{n}\sum_{k=0}^{n-1}\bold f(S^{k}\bold i)
\subseteq
C
 \Bigg\}\,.\\
 \endalign
 $$
In particular, if the OSC is satisfied and $\pmb\alpha\in\Bbb R^{M}$, then we have 
$$
\align
\quad
\,\,\scr f\,\,(\pmb\alpha)
 &=
 \dim_{\Haus}
 \pi
 \Bigg\{
 \bold i\in\Sigma_{\smallG}^{\Bbb N}
 \,\Bigg|\,
 \lim_{n}
 \frac{1}{n}\sum_{k=0}^{n-1}\bold f(S^{k}\bold i)
 =
 \pmb\alpha
 \Bigg\}\,.
 \endalign
 $$
\endroster

 \bigskip
 
 \noindent
 {\rm (2)} Assume that $C\subseteq\Bbb R^{M}$ is  closed 
and
convex
 with
 $\overset{\,\circ}\to{C}
  \cap
  \{
 \int\bold f\,d\mu
  \,|\,
  \mu\in\Cal P_{S}(\Sigma_{\smallG}^{\Bbb N})
  \}
    \not=
  \varnothing$.

\roster

 \item"(2.1)"
 There is a unique real number 
 $\scr F\,\,(C)$ such that
  $$
  \sigma_{\radius}
  \big(
  \,
  \zeta_{C}^{\dyn\text{-}\vector\text{-}\erg}(\scr F\,\,(C)\,\Lambda;\cdot)
  \,
  \big)
  =
  1\,.
  $$

\item"(2.2)"
We have
 $$
 \scr F\,\,(C)
 =
 \sup_{\pmb\alpha\in C}
 \,
 \sup
 \Sb
 \mu\in\Cal P_{S}(\Sigma_{\smallsmallG}^{\Bbb N})\\
 \int \bold f\,d\mu=\pmb\alpha
 \endSb
 -
 \frac{h(\mu)}{\int\Lambda\,d\mu}\,.
 $$

\item"(2.3)"
If the OSC is satisfied, then we have
$$
\align
\quad
 \scr F\,\,(C)
 &=
 \dim_{\Haus}
 \pi
 \Bigg\{
 \bold i\in\Sigma_{\smallG}^{\Bbb N}
 \,\Bigg|\,
 \,\underset{n}\to\acc
 \frac{1}{n}\sum_{k=0}^{n-1}\bold f(S^{k}\bold i)
\subseteq
C
 \Bigg\}\,.\\
 \endalign
 $$

\endroster

\endproclaim

\noindent{\it Proof}\newline
\noindent
This follows immediately
by applying
Theorem 5.5
to $W=\Bbb R^{M}$
and the maps
$\pmb\Phi:\Sigma_{\smallG}^{\Bbb N}\to\Bbb R^{M}$
and
 $Q:W\to\Bbb R^{M}$ 
defined by
$\pmb\Phi=\bold f$
and
$Q(\bold x)=\bold x$.
\hfill$\square$

 \bigskip

%

%

%

%

%

\proclaim{Theorem 5.7.
Multifractal zeta-functinons for
multifractal spectra of relative ergodic Birkhoff averages}
Fix $\gamma\in(0,1)$
and let $f,g:\Sigma_{\smallG}^{\Bbb N}\to\Bbb R$
be Lipschitz functions
with respect to the metric
$\distance_{\gamma}$
and
assume that $g(\bold i)\not=0$
for all $\bold i$.
For $C\subseteq\Bbb R$ and 
a continuous function 
$\varphi:\Sigma_{\smallG}^{\Bbb N}\to\Bbb R$,   we
define the dynamical
relative ergodic multifractal zeta-function by
  $$
  \zeta_{C}^{\dyn\text{-}\rel\text{-}\erg}(\varphi;z)
  =
  \sum_{n}
  \frac{z^{n}}{n}
  \left(
  \sum
  \Sb
 \bold i\in\Sigma_{\smallsmallG}^{n}\\
  {}\\
  \forall\bold u\in[\bold i]
  \,\,:\,\,
 \frac
 {\sum_{k=0}^{n-1}f(S^{k}\bold u)}
 {\sum_{k=0}^{n-1}g(S^{k}\bold u)}
 \in  C
  \endSb
  \sup_{\bold u\in[\bold i]}
  \exp
  \sum_{k=0}^{n-1}
  \varphi S^{k}\bold u
  \right)\,.
  $$
Let 
 $\Lambda:\Sigma_{\smallG}^{\Bbb N}\to\Bbb R$ be defined by (2.10).

\noindent
 {\rm (1)} Assume 
 that $C\subseteq\Bbb R^{M}$ is closed.

 \roster
  \item"(1.1)"
 There is a unique real number 
 $\,\,\scr f\,\,(C)$ such that
  $$
  \lim_{r\searrow0}
  \sigma_{\radius}
  \big(
  \,
  \zeta_{B(C,r)}^{\dyn\text{-}\rel\text{-}\erg}(\,\,\scr f\,\,(C)\,\Lambda;\cdot)
  \,
  \big)
  =
  1\,.
  $$
If $\alpha\in\Bbb R$ and 
$C=\{\alpha\}$, then we will write
  $\,\,\scr f\,\,(\alpha)
  =
 \,\,\scr f\,\,(C)$.

\item"(1.2)"
We have
$$
\,\,\scr f\,\,(C)
 =
 \sup_{\alpha\in C}
 \,
 \sup
  \Sb
  \mu\in\Cal P_{S}(\Sigma_{\smallsmallG}^{\Bbb N})\\
  {}\\
  \frac{\int f\,d\mu}{\int g\,d\mu}=\alpha
  \endSb
 -
 \frac{h(\mu)}{\int\Lambda\,d\mu}\,.
 $$

 \item"(1.3)"
If the OSC is satisfied, then we have
$$
\align
\quad
\,\,\scr f\,\,(C)
 &=
 \dim_{\Haus}
 \pi
 \Bigg\{
 \bold i\in\Sigma_{\smallG}^{\Bbb N}
 \,\Bigg|\,
 \,\underset{n}\to\acc
 \,\,
 \frac
 {\sum_{k=0}^{n-1}f(S^{k}\bold i)}
 {\sum_{k=0}^{n-1}g(S^{k}\bold i)}
\subseteq
C
 \Bigg\}\,.\\
 \endalign
 $$
In particular, if the OSC is satisfied and $\alpha\in\Bbb R$, then we have 
$$
\align
\quad
\,\,\scr f\,\,(\alpha)
 &=
 \dim_{\Haus}
 \pi
 \Bigg\{
 \bold i\in\Sigma_{\smallG}^{\Bbb N}
 \,\Bigg|\,
 \lim_{n}
 \frac
 {\sum_{k=0}^{n-1}f(S^{k}\bold i)}
 {\sum_{k=0}^{n-1}g(S^{k}\bold i)} =
 \alpha
 \Bigg\}\,.
 \endalign
 $$
\endroster

 \bigskip
 
 \noindent
 {\rm (2)} Assume that $C\subseteq\Bbb R^{M}$ is  closed 
and
convex
 with
 $\overset{\,\circ}\to{C}
  \cap
   \{
  \frac{\int f\,d\mu}{\int g\,d\mu}
  \,|\,
  \mu\in\Cal P_{S}(\Sigma_{\smallG}^{\Bbb N})
  \}
    \not=
  \varnothing$.

\roster

 \item"(2.1)"
 There is a unique real number 
 $\scr F\,\,(C)$ such that
  $$
  \sigma_{\radius}
  \big(
  \,
  \zeta_{C}^{\dyn\text{-}\rel\text{-}\erg}(\scr F\,\,(C)\,\Lambda;\cdot)
  \,
  \big)
  =
  1\,.
  $$

\item"(2.2)"
We have
 $$
 \scr F\,\,(C)
 =
 \sup_{\alpha\in C}
 \,
 \sup
 \Sb
 \mu\in\Cal P_{S}(\Sigma_{\smallsmallG}^{\Bbb N})\\
   \frac{\int f\,d\mu}{\int g\,d\mu}=\alpha
 \endSb
 -
 \frac{h(\mu)}{\int\Lambda\,d\mu}\,.
 $$

\item"(2.3)"
If the OSC is satisfied, then we have
$$
\align
\quad
 \scr F\,\,(C)
 &=
 \dim_{\Haus}
 \pi
 \Bigg\{
 \bold i\in\Sigma_{\smallG}^{\Bbb N}
 \,\Bigg|\,
 \,\underset{n}\to\acc
 \,\,
 \frac
 {\sum_{k=0}^{n-1}f(S^{k}\bold i)}
 {\sum_{k=0}^{n-1}g(S^{k}\bold i)}
\subseteq
C
 \Bigg\}\,.\\
 \endalign
 $$

\endroster

\endproclaim

\noindent{\it Proof}\newline
\noindent
This follows immediately
by applying
Theorem 5.5
to $W=\Bbb R\times (\Bbb R\setminus\{0\})$
and the maps
$\pmb\Phi:\Sigma_{\smallG}^{\Bbb N}\to\Bbb R^{2}$
and
 $Q:W\to\Bbb R$ 
defined by
$\pmb\Phi=(f,g)$
and
$Q(x,y)=\frac{x}{y}$.
\hfill$\square$

\bigskip

\newpage

As a final application of Theorem 5.5 we
consider
a type of relative ergodic 
multifractal spectra involving
quantities 
similar to those 
appearing in H\"older's inequality;
for this reason we have decided
to refer to these multifractal spectra as
\lq\lq H\"older-like relative ergodic Birkhoff averages".

\bigskip

\proclaim{Theorem 5.8.
Multifractal zeta-functinons for
multifractal spectra of 
H\"older-like relative ergodic Birkhoff averages}
Fix $\gamma\in(0,1)$
and let $f_{1},\ldots,f_{M},g_{1},\ldots,g_{M}:\Sigma_{\smallG}^{\Bbb N}\to\Bbb R$
be Lipschitz functions
with respect to the metric
$\distance_{\gamma}$
and
assume that 
$f_{l}(\bold i)>0$
for all $l$ and all $\bold i$,
and that
$g_{l}(\bold i)>0$
for all $l$ and all $\bold i$.
Fix
 $s_{1},\ldots,s_{M},t_{1},\ldots,t_{M}>0$.
For $C\subseteq\Bbb R$ and 
a continuous function 
$\varphi:\Sigma_{\smallG}^{\Bbb N}\to\Bbb R$,   we
define the dynamical
H\"older-like
relative ergodic multifractal zeta-function by
  $$
  \zeta_{C}^{\dyn\text{-}\Hol\text{-}\erg}(\varphi;z)
  =
  \sum_{n}
  \frac{z^{n}}{n}
  \left(
  \sum
  \Sb
 \bold i\in\Sigma_{\smallsmallG}^{n}\\
  {}\\
  \forall\bold u\in[\bold i]
  \,\,:\,\,
 \frac
 {\prod_{l=1}^{M}\big(\frac{1}{n}\sum_{k=0}^{n-1}f_{l}(S^{k}\bold u)\big)^{s_{l}}}
 {\prod_{l=1}^{M}\big(\frac{1}{n}\sum_{k=0}^{n-1}g_{l}(S^{k}\bold u)\big)^{t_{l}}}
 \in  C
  \endSb
  \sup_{\bold u\in[\bold i]}
  \exp
  \sum_{k=0}^{n-1}
  \varphi S^{k}\bold u
  \right)\,.
  $$
Let 
 $\Lambda:\Sigma_{\smallG}^{\Bbb N}\to\Bbb R$ be defined by (2.10).
Assume 
 that $C\subseteq\Bbb R^{M}$ is closed.

 \roster
  \item"(1)"
 There is a unique real number 
 $\,\,\scr f\,\,(C)$ such that
  $$
  \lim_{r\searrow0}
  \sigma_{\radius}
  \big(
  \,
  \zeta_{B(C,r)}^{\dyn\text{-}\Hol\text{-}\erg}(\,\,\scr f\,\,(C)\,\Lambda;\cdot)
  \,
  \big)
  =
  1\,.
  $$
If $\alpha\in\Bbb R$ and 
$C=\{\alpha\}$, then we will write
  $\,\,\scr f\,\,(\alpha)
  =
 \,\,\scr f\,\,(C)$.

\item"(2)"
We have
$$
\,\,\scr f\,\,(C)
 =
 \sup_{\alpha\in C}
 \,
 \sup
  \Sb
  \mu\in\Cal P_{S}(\Sigma_{\smallsmallG}^{\Bbb N})\\
  {}\\
 \frac
 {\prod_{l=1}^{M}\big(\int f_{l}\,d\mu\big)^{s_{l}}}
 {\prod_{l=1}^{M}\big(\int g_{l}\,d\mu\big)^{t_{l}}}
 =\alpha
  \endSb
 -
 \frac{h(\mu)}{\int\Lambda\,d\mu}\,.
 $$

 \item"(3)"
If the OSC is satisfied, then we have
$$
\align
\quad
\,\,\scr f\,\,(C)
 &=
 \dim_{\Haus}
 \pi
 \Bigg\{
 \bold i\in\Sigma_{\smallG}^{\Bbb N}
 \,\Bigg|\,
 \,\underset{n}\to\acc
 \,\,
 \frac
 {\prod_{l=1}^{M}\big(\frac{1}{n}\sum_{k=0}^{n-1}f_{l}(S^{k}\bold u)\big)^{s_{l}}}
 {\prod_{l=1}^{M}\big(\frac{1}{n}\sum_{k=0}^{n-1}g_{l}(S^{k}\bold u)\big)^{t_{l}}}
\subseteq
C
 \Bigg\}\,.\\
 \endalign
 $$
In particular, if the OSC is satisfied and $\alpha\in\Bbb R$, then we have 
$$
\align
\quad
\,\,\scr f\,\,(\alpha)
 &=
 \dim_{\Haus}
 \pi
 \Bigg\{
 \bold i\in\Sigma_{\smallG}^{\Bbb N}
 \,\Bigg|\,
 \lim_{n}
 \frac
 {\prod_{l=1}^{M}\big(\frac{1}{n}\sum_{k=0}^{n-1}f_{l}(S^{k}\bold u)\big)^{s_{l}}}
 {\prod_{l=1}^{M}\big(\frac{1}{n}\sum_{k=0}^{n-1}g_{l}(S^{k}\bold u)\big)^{t_{l}}} =
 \alpha
 \Bigg\}\,.
 \endalign
 $$
\endroster

\endproclaim

\noindent{\it Proof}\newline
\noindent
This follows immediately
by applying
Theorem 5.5
to $W=\Bbb R^{M}\times (\Bbb R\setminus\{0\})^{M}$
and the maps
$\pmb\Phi:\Sigma_{\smallG}^{\Bbb N}\to\Bbb R^{2M}$
and
 $Q:W\to\Bbb R$ 
defined by
$\pmb\Phi=(f_{1},\ldots,f_{M},g_{1},\ldots,g_{M})$
and
$Q(x_{1},\ldots,x_{M},y_{1},\ldots,y_{M})
=
\frac
{\prod_{l=1}^{M}x_{l}^{s_{l}}}
{\prod_{l=1}^{M}y_{l}^{t_{l}}}$.
\hfill$\square$

  \bigskip


\newpage

\heading
{
6. Proofs. The 
map $M_{n}$}
\endheading

The purpose of this section 
to introduce the 
auxiliary map $M_{n}$
and to prove various 
continuity results
regarding this map.
The two main results are 
Lemma 6.2 and Lemma 6.3.
Loosely speaking these 
lemmas say
that
the maps $L_{n}$ and $M_{n}$ 
behave
asymptotically
in 
the same
way.
We also state and prove a simple but useful lemma
about upper semi-continuous maps, see Lemma 6.4.
All the key lemmas is this section 
(i\.e\. Lemma 6.2, Lemma 6.3 and Lemma 6.4)
play
important roles in the following sections.

\bigskip

{\bf The map $M_{n}$.}
Since the graph
$\G=(\V,\E)$ is strongly connected,
it follows that
for
each $\bold i\in\Sigma_{\smallG}^{*}$, 
 we 
can choose
$\widehat{\bold i}\in\Sigma_{\smallG}^{*}$ 
with
$|\,\widehat{\bold i}\,|\le|\V|$
such that
$\termi(\bold i)=\ini(\,\widehat{\bold i}\,)$
and
$\termi(\,\widehat{\bold i}\,)=\ini(\bold i)$.
Next,
for
$\bold i\in\Sigma_{\smallG}^{*}$,
define
$\overline{\bold i}\in\Sigma_{\smallG}^{\Bbb N}$ by
 $$
 \overline{\bold i}
 =
 \bold i\,\,\widehat{\bold i}\,\,
 \bold i\,\,\widehat{\bold i}\,\,
 \bold i\,\,\widehat{\bold i}\ldots
 $$
 Finally,
 for  a positive integer $n$,
 we define
$M_{n}:\Sigma_{\smallG}^{\Bbb N}\to\Cal P_{S}(\Sigma_{\smallG}^{\Bbb N})$ by
 $$
 \align
 M_{n}\bold i
&=
 L_{n+|\,\widehat{\,\bold i|n\,}\,|}\left(\,\overline{\bold i|n}\,\right)\\
&= 
 \frac{1}{n+|\,\widehat{\,\bold i|n\,}\,|}
 \sum_{k=0}^{n+|\,\widehat{\,\bold i|n\,}\,|-1}
 \delta_{S^{k}(\,\overline{\bold i|n}\,)}
 \tag6.1
 \endalign
 $$
for $\bold i\in\Sigma_{\smallG}^{\Bbb N}$;
recall, that
the map
$L_{n}:\Sigma_{\smallG}^{\Bbb N}\to\Cal P(\Sigma_{\smallG}^{\Bbb N})$ is defined 
by
by
 $$
 L_{n}\bold i
 =
 \frac{1}{n}\sum_{k=0}^{n-1}\delta_{S^{k}\bold i}
 $$
for $\bold i\in\Sigma_{\smallG}^{\Bbb N}$;
see
(4.1).

\bigskip

{\bf Why the map $M_{n}$?}
We will now 
explain the main technical reason for introducing  the map $M_{n}$.
In order to prove Theorem 4.4 we will
use results from large deviation theory.
In particular, we
will use
Varadhan's
integral lemma
(i\.e\. Theorem 8.1)
which says that if $X$
is a complete separable metric space
and $(P_{n})_{n}$ is a sequence of probability measures on
$X$ satisfying the large deviation property
with
rate constants $a_{n}\in(0,\infty)$ for $n\in\Bbb N$
and rate function $I:\Bbb R\to[-\infty,\infty]$
(this terminology will be explained in Section 7),
then
$$
\lim_{n}
\frac{1}{a_{n}}
\log
\int
\exp(a_{n}F)\,dP_{n}
=
-
\inf_{x\in X}(\,I(x)-F(x))
$$
for any bounded continuous function $F:X\to \Bbb R$
(see Section 8 for more a detailed  discussion and statement of this result).

More precisely, in Section 8
we  use 
Varadhan's
integral lemma
to analyse the asymptotic
behaviour of the integral
$$
\frac{1}{n}
\log
\int
\exp(nF_{\varphi}(L_{n}(\overline{\bold i|n})))\,d\Pi(\bold i)
\tag6.2
$$
as $n\to\infty$
where $\Pi$ is the Parry measure on $\Sigma_{\smallG}^{\Bbb N}$
(the Parry measure will be defined in Section 7)
 and the function
 $F_{\varphi}:\Cal P\big(\Sigma_{\smallG}^{\Bbb N}\big)\to\Bbb R$
is given by
 $
 F_{\varphi}(\mu)
 =
\int\varphi\,d\mu$.
Defining
$\Lambda_{n}: \Cal P\big(\Sigma_{\smallG}^{\Bbb N}\big)\to\Bbb R$
by
$\Lambda_{n}(\bold i)
 =
 L_{n}(\,\overline{\bold i|n}\,)$,
 then (6.2) can be written as
$$
\frac{1}{n}
\log
\int
\exp(nF_{\varphi})\,d(\Pi\circ\Lambda_{n}^{-1})\,.
\tag6.3
$$
Consequently,
if the sequence $(\Pi\circ\Lambda_{n}^{-1})_{n}$
satisfied the large deviation property
with rate constants $a_{n}=n$,
then
Varadhan's
integral lemma
could be applied to analyse the asymptotic behaviour of the sequence of integrals in 
(6.3). 
 However, 
 it follows from
 results by Orey \& Pelikan [OrPe1,OrPe2]
 that the sequence
 $(\Pi\circ L_{n})_{n}$
satisfies the large deviation property
with rate constants $a_{n}=n$
and
Varadhan's
integral lemma can therefore be applied to 
provide information about the asymptotic behaviour of the sequence of integral
defined by
$$
\frac{1}{n}
\log
\int
\exp(nF_{\varphi})\,d(\Pi\circ L_{n}^{-1})\,.
\tag6.4
$$
In order to utilise the knowledge of the
asymptotic behaviour of (6.4)
for analysing the asymptotic behaviour of (6.3),
we must therefore show that the measures
 $$
 \Pi\circ\Lambda_{n}^{-1}
 $$
 and
  $$
  \Pi\circ L_{n}^{-1}
  $$
  are
  \lq\lq close".
However, 
for technical reasons we will prove and use a similar 
\lq\lq closeness"
statement involving the measures
  $$
  \Pi\circ M_{n}^{-1}
  $$
  and
  $$
  \Pi\circ L_{n}^{-1}\,.
  $$

   Indeed, below we prove a number
   of
   results
showing that the maps
$M_{n}$
  and
$L_{n}$
(and therefore also the
measures
$\Pi\circ M_{n}^{-1}$
  and
$\Pi\circ L_{n}^{-1}$)
are 
\lq\lq close".
These results
play an important role in
Section 7.
In particular, they
allow us to:
(1)
use 
Orey \& Pelikan's
result from [OrPe1,OrPe2] saying that
the sequence
$(\Pi\circ L_{n}^{-1})_{n}$
satisfies the large deviation property
to
prove
that
the sequence
$(\Pi\circ M_{n}^{-1})_{n}$
also satisfies the large deviation property (see Theorem 7.2),
and (2)
replace all occurencies
of 
$L_{n}(\overline{\bold i|n})$
in the formulas in Section 8 by
$M_{n}\bold i$
allowing us to use
the large deviation property of the sequence
$(\Pi\circ M_{n}^{-1})_{n}$.
 This
 explains the mean reason for introducing the map $M_{n}$
 and the associated measure $\Pi\circ M_{n}^{-1}$.

\bigskip

{\bf Comparing $M_{n}$ and $L_{n}$.}
 We now prove various 
continuity statements saying that the maps
$M_{n}$
  and
$L_{n}$
are 
\lq\lq close".
These statements play an important role in Section 7
where we
apply
Varadhan's
integral lemma 
to prove Theorem 7.2.
We first introduce the
metric $\LDistance$
on 
$\Cal P(\Sigma_{\smallG}^{\Bbb N})$.
Fix $\gamma\in(0,1)$
and let $\distance_{\gamma}$
denote the metric on $\Sigma_{\smallG}^{\Bbb N}$
introduced in Section 5.2.
For a function $f:\Sigma_{\smallG}^{\Bbb N}\to\Bbb R$, we let
$\Lip_{\gamma}(f)$ denote the Lipschitz constant of $f$
with respect to the metric
$\distance_{\gamma}$, 
i\.e\.
$\Lip_{\gamma}(f)
=
\sup_{
\bold i,\bold j\in\Sigma_{\smallsmallG}^{\Bbb N},
\bold i\not=\bold j}
\frac
{|f(\bold i)-f(\bold j)|}
{\distance_{\gamma}(\bold i,\bold j)}$
and
we
define the metric $\LDistance$ in $\Cal P(\Sigma_{\smallG}^{\Bbb N})$ by
$$
\LDistance(\mu,\nu)
=
\sup
\Sb
f:\Sigma_{\smallsmallG}^{\Bbb N}\to\Bbb R\\
\Lip_{\gamma}(f)\le 1
\endSb
\Bigg|
\int f\,d\mu-\int f\,d\nu
\Bigg|;
\tag6.5
$$
we note that it is well-known that $\LDistance$ is a metric
and  that $\LDistance$ induces the weak topology.
Below we will always equip 
the space
$\Cal P(\Sigma_{\smallG}^{\Bbb N})$ 
with the metric $\LDistance$.
Before stating and proving the key results in this section,
 we start by proving a small 
technical auxiliary result.

\bigskip

\proclaim{Lemma 6.1}
Let $(X,\distance)$ be a metric space
and let $U:\Cal P(\Sigma_{\smallG}^{\Bbb N})\to X$ be continuous with respect to the weak topology.
Fix $r>0$.
There is a positive integer $N_{r}$
such that
if 
$n\ge N_{r}$,
$\bold u\in\Sigma_{\smallG}^{n}$ and $\bold k,\bold l\in\Sigma_{\smallG}^{\Bbb N}$
with
$\termi(\bold u)=\ini(\bold k)$
and
$\termi(\bold u)=\ini(\bold l)$, then
we have
 $$
 \distance
 \big(
 \,
 UL_{n}(\bold u\bold k)
 \,,\,
 UM_{n}(\bold u\bold l)
 \,
 \big)
 <
 r\,.
 $$

\endproclaim
\noindent{\it Proof}\newline
\noindent
Fix $\gamma\in(0,1)$ and let $\LDistance$ be the metric in 
$\Cal P(\Sigma_{\smallG}^{\Bbb N})$
defined in (6.5).
Since 
$U:\Cal P(\Sigma_{\smallG}^{\Bbb N})\to X$ is continuous and 
$\Cal P(\Sigma_{\smallG}^{\Bbb N})$ is compact, we
conclude that
$U:\Cal P(\Sigma_{\smallG}^{\Bbb N})\to X$
is uniformly continuous.
This implies that we can choose $\delta>0$
such that
all measures
$\mu,\nu\in\Cal P(\Sigma_{\smallG}^{\Bbb N})$
satisfy the following implication:
 $$
 \LDistance(\mu,\nu)<\delta
 \,\,\,\,
 \Rightarrow
 \,\,\,\,
 \distance(U\mu,U\nu)<\frac{r}{2}\,.
 \tag6.6
 $$

Next,
choose a positive integer $N_{r}$ such that
$2\frac{|\smallV|}{N_{r}}
 <
\delta$
 and
$\frac{1}{N_{r}(1-\gamma)}
 <
 \delta$.

Now, fix
$n\ge N_{r}$,
$\bold u\in\Sigma_{\smallG}^{n}$ and $\bold k,\bold l\in\Sigma_{\smallG}^{\Bbb N}$
with
$\termi(\bold u)=\ini(\bold k)$
and
$\termi(\bold u)=\ini(\bold l)$.
It follows that
$$
 \align
\LDistance
 \big(
 \,
 L_{n}(\bold u\bold k)
 \,,\,
 M_{n}(\bold u\bold l)
 \,
 \big)
&\le
 \LDistance
 \big(
 \,
 L_{n}(\bold u\bold k)
 \,,\,
 L_{n}(\bold u\bold l)
 \,
 \big)
+
\LDistance
 \big(
 \,
 L_{n}(\bold u\bold l)
 \,,\,
 M_{n}(\bold u\bold l)
 \,
 \big)\\
&\le
 \LDistance
 \big(
 \,
 L_{n}(\bold u\bold k)
 \,,\,
 L_{n}(\bold u\bold l)
 \,
 \big)
+
\LDistance
 \big(
 \,
 L_{n}(\bold u\bold l)
 \,,\,
 L_{n+|\,\widehat{\bold u}\,|}\overline{\bold u}
 \,
 \big) 
\endalign
$$

 We first estimate the distance
 $\LDistance
 \big(
 \,
 L_{n}(\bold u\bold k)
 \,,\,
 L_{n}(\bold u\bold l)
 \,
 \big)$.
 Indeed,
 since $\frac{1}{N_{r}(1-\gamma)}<\delta$,
 it follows that
 $$
 \align
\LDistance
 \big(
 \,
 L_{n}(\bold u\bold k)
 \,,\,
 L_{n}(\bold u\bold l)
 \,
 \big)
&=
\sup
\Sb
f:\Sigma_{\smallsmallG}^{\Bbb N}\to\Bbb R\\
\Lip_{\gamma}(f)\le 1
\endSb
\Bigg|
\int f\,d(L_{n}(\bold u\bold k))-\int f\,d(L_{n}(\bold u\bold l))
\Bigg|\\
&=
\sup
\Sb
f:\Sigma_{\smallsmallG}^{\Bbb N}\to\Bbb R\\
\Lip_{\gamma}(f)\le 1
\endSb
\Bigg|
\frac{1}{n}\sum_{i=0}^{n-1} f(S^{i}(\bold u\bold k))
-
\frac{1}{n}\sum_{i=0}^{n-1} f(S^{i}(\bold u\bold l))
\Bigg|\\
&\le
\sup
\Sb
f:\Sigma_{\smallsmallG}^{\Bbb N}\to\Bbb R\\
\Lip_{\gamma}(f)\le 1
\endSb
\frac{1}{n}\sum_{i=0}^{n-1} 
|
f(S^{i}(\bold u\bold k))
-
 f(S^{i}(\bold u\bold l))
|\\
&\le
\frac{1}{n}\sum_{i=0}^{n-1} 
\distance_{\gamma}
\big(
\,
S^{i}(\bold u\bold k)
\,,\,
S^{i}(\bold u\bold l)
\,
\big)\\
&=
\frac{1}{n}\sum_{i=0}^{n-1} 
\gamma^{
|S^{i}(\bold u\bold k)
\wedge
S^{i}(\bold u\bold l)|
}\\
&\le
\frac{1}{N_{r}}\sum_{i=0}^{n-1} 
\gamma^{n-i}\\
&\le
 \frac{1}{N_{r}(1-\gamma)}\\
&<
\delta\,,
\endalign
$$
and we therefore conclude from (6.6) that
 $$
 \distance
 \big(
 \,
 UL_{n}(\bold u\bold k)
 \,,\,
 UL_{n}(\bold u\bold l)
 \,
 \big)
 <
\frac{r}{2}\,.
\tag6.7
$$

 Next, we estimate the distance
 $\LDistance
 \big(
 \,
 L_{n}(\bold u\bold l)
 \,,\,
 L_{n+|\,\widehat{\bold u}\,|}\overline{\bold u}
 \,
 \big)$.
We start by observing that
if we fix $\bold i_{0}\in\Sigma_{\smallG}^{\Bbb N}$, then
$$
\align
\LDistance(\mu,\nu)
&=
\sup
\Sb
f:\Sigma_{\smallsmallG}^{\Bbb N}\to\Bbb R\\
\Lip_{\gamma}(f)\le 1
\endSb
\Bigg|
\int f\,d\mu-\int f\,d\nu
\Bigg|\\
&=
\sup
\Sb
f:\Sigma_{\smallsmallG}^{\Bbb N}\to\Bbb R\\
\Lip_{\gamma}(f)\le 1
\endSb
\Bigg|
\int (f-f(\bold i_{0}))\,d\mu-\int (f-f(\bold i_{0}))\,d\nu
\Bigg|\\
&=
\sup
\Sb
g:\Sigma_{\smallsmallG}^{\Bbb N}\to\Bbb R\\
\Lip_{\gamma}(g)\le 1\\
g(\bold i_{0})=0
\endSb
\Bigg|
\int g\,d\mu-\int g\,d\nu
\Bigg|
\tag6.8
\endalign
$$
for all $\mu,\nu\in\Cal P(\Sigma_{\smallG}^{\Bbb N})$.
It  follows from 
(6.8) that
 $$
 \align
\LDistance
 \big(
 \,
 L_{n}(\bold u\bold l)
 \,,\,
 L_{n+|\,\widehat{\bold u}\,|}\overline{\bold u}
 \,
 \big)
&=
\sup
\Sb
g:\Sigma_{\smallsmallG}^{\Bbb N}\to\Bbb R\\
\Lip_{\gamma}(g)\le 1\\
g(\bold i_{0})=0
\endSb
\Bigg|
\int g\,d(L_{n}(\bold u\bold l))-\int g\,d(L_{n+|\,\widehat{\bold u}\,|}\overline{\bold u})
\Bigg|\\
&=
\sup
\Sb
g:\Sigma_{\smallsmallG}^{\Bbb N}\to\Bbb R\\
\Lip_{\gamma}(g)\le 1\\
g(\bold i_{0})=0
\endSb
\Bigg|
\frac{1}{n}\sum_{i=0}^{n-1} g(S^{i}(\bold u\bold l))
-
\frac{1}{n+|\,\widehat{\bold u}\,|}\sum_{i=0}^{n+|\,\widehat{\bold u}\,|-1} 
g(S^{i}\overline{\bold u})
\Bigg|\\
&\le
\sup
\Sb
g:\Sigma_{\smallsmallG}^{\Bbb N}\to\Bbb R\\
\Lip_{\gamma}(g)\le 1\\
g(\bold i_{0})=0
\endSb
\Bigg(
\,\,
\Bigg|
\frac{1}{n}\sum_{i=0}^{n-1} g(S^{i}(\bold u\bold l))
-
\frac{1}{n+|\,\widehat{\bold u}\,|}\sum_{i=0}^{n-1} g(S^{i}\overline{\bold u})
\Bigg|\\
&\qquad\qquad
\qquad\qquad
\qquad\qquad
\qquad\qquad
+
\Bigg|
\frac{1}{n+|\,\widehat{\bold u}\,|}
\sum_{i=n}^{n+|\,\widehat{\bold u}\,|-1} g(S^{i}\overline{\bold u})
\Bigg|
\,
\Bigg)\\
&\le
\sup
\Sb
g:\Sigma_{\smallsmallG}^{\Bbb N}\to\Bbb R\\
\Lip_{\gamma}(g)\le 1\\
g(\bold i_{0})=0
\endSb
\Bigg(
\frac{|\,\widehat{\bold u}\,|}{n(n+|\,\widehat{\bold u}\,|)}
\sum_{i=0}^{n-1} 
\|g\|_{\infty}
+
\frac{1}{n+|\,\widehat{\bold u}\,|}
\sum_{i=n}^{n+|\,\widehat{\bold u}\,|-1} \|g\|_{\infty}
\Bigg)\\
&=
\sup
\Sb
g:\Sigma_{\smallsmallG}^{\Bbb N}\to\Bbb R\\
\Lip_{\gamma}(g)\le 1\\
g(\bold i_{0})=0
\endSb
2\frac{|\,\widehat{\bold u}\,|}{n+|\,\widehat{\bold u}\,|}
\,
\|g\|_{\infty}\\
&\le
\sup
\Sb
g:\Sigma_{\smallsmallG}^{\Bbb N}\to\Bbb R\\
\Lip_{\gamma}(g)\le 1\\
g(\bold i_{0})=0
\endSb
2\frac{|\V| }{n}
\,
\|g\|_{\infty}\,.
\qquad\qquad
\qquad\qquad
\text{[since $|\,\widehat{\bold u}\,|\le|\V|$]}
\tag6.9
\endalign
$$
However, 
if
$g:\Sigma_{\smallG}^{\Bbb N}\to\Bbb R$ satisfies
$\Lip_{\gamma}(g)\le 1$ and
$g(\bold i_{0})=0$,
then
$|g(\bold i)|
=
|g(\bold i)-g(\bold i_{0})|
\le
\distance_{\gamma}(\bold i,\bold i_{0})
\le
1$
for all $\bold i\in\Sigma_{\smallG}^{\Bbb N}$, whence
$\|g\|_{\infty}\le 1$.
It therefore follows from (6.9) that
if 
$n\ge N_{r}$,
$\bold u\in\Sigma_{\smallG}^{n}$ and $\bold l\in\Sigma_{\smallG}^{\Bbb N}$
with
$\termi(\bold u)=\ini(\bold l)$, then
 $\LDistance
 \big(
 \,
 L_{n}(\bold u\bold l)
 \,,\,
 L_{n+|\,\widehat{\bold u}\,|}\overline{\bold u}
 \,
 \big)
 \le
 2\frac{|\smallV|}{n}
 \le
 2\frac{|\smallV|}{N_{r}}
 <
 \delta$,
 and we therefore conclude from (6.6) that
 $$
 \distance
 \big(
 \,
 UL_{n}(\bold u\bold l)
 \,,\,
 L_{n+|\,\widehat{\bold u}\,|}\overline{\bold u}
 \,
 \big)
 <
 \frac{r}{2}\,.
 \tag6.10
 $$

Finally, combining (6.7) and (6.10) shows 
that
$\LDistance
 \big(
 \,
 L_{n}(\bold u\bold k)
 \,,\,
 M_{n}(\bold u\bold l)
 \,
 \big)
\le
 \LDistance
 \big(
 \,
 L_{n}(\bold u\bold k)
 \,,\,
 L_{n}(\bold u\bold l)
 \,
 \big)
+
\LDistance
 \big(
 \,
 L_{n}(\bold u\bold l)
 \,,\,
 L_{n+|\,\widehat{\bold u}\,|}\overline{\bold u}
 \,
 \big)
 <
 \frac{r}{2}
 +
 \frac{r}{2}
 =
 r$.
 \hfill$\square$

\bigskip

\noindent
We can now state and prove the  first two key 
lemmas in this section.

\bigskip

\proclaim{Lemma 6.2}
Let $(X,\distance)$ be a metric space
and let $U:\Cal P(\Sigma_{\smallG}^{\Bbb N})\to X$ be continuous with respect to the weak topology.
Fix
$r>0$.
There is a positive integer $N_{r}$
such that
if $n\ge N_{r}$
and $C$ is a subset of $X$,
then
we have
 $$
 \aligned
 \Big\{
 \bold u\in\Sigma_{\smallG}^{n}
 \,\Big|\,
 UL_{n}[\bold u]\subseteq C
\Big\}
&\subseteq
\Big\{
 \bold u\in\Sigma_{\smallG}^{n}
 \,\Big|\,
 UM_{n}[\bold u]\subseteq B(C,r)
\Big\}\,,\\
 \Big\{
 \bold i\in\Sigma_{\smallG}^{\Bbb N}
 \,\Big|\,
 UL_{n}\bold i\in C
\Big\}
&\subseteq
\Big\{
 \bold i\in\Sigma_{\smallG}^{\Bbb N}
 \,\Big|\,
 UM_{n}\bold i\in B(C,r)
\Big\}\,,
\endaligned
\tag6.11
$$
and
 $$
 \aligned
\Big\{
 \bold u\in\Sigma_{\smallG}^{n}
 \,\Big|\,
 UM_{n}[\bold u]\subseteq C
\Big\}
&\subseteq
 \Big\{
 \bold u\in\Sigma_{\smallG}^{n}
 \,\Big|\,
 \,\,
 UL_{u}[\bold u]\subseteq B(C,r)\Big\}\,,\\
  \Big\{
 \bold i\in\Sigma_{\smallG}^{\Bbb N}
 \,\Big|\,
 UM_{n}\bold i\in C
\Big\}
&\subseteq
\Big\{
 \bold i\in\Sigma_{\smallG}^{\Bbb N}
 \,\Big|\,
 UL_{n}\bold i\in B(C,r)
\Big\}\,.
\endaligned
\tag6.12
$$
\endproclaim
\noindent{\it Proof}\newline
\noindent
It follows from Lemma 6.1 that
there is a positive integer $N_{r}$
such that
if
$n\ge N_{r}$,
$\bold u\in\Sigma^{n}$ and $\bold k,\bold l\in\Sigma_{\smallG}^{\Bbb N}$
with
$\termi(\bold u)=\ini(\bold k)$
and
$\termi(\bold u)=\ini(\bold l)$, then
 $\distance
 (
 \,
 UL_{n}(\bold u\bold k)
 \,,\,
 UM_{n}(\bold u\bold l)
 \,
 )
<
 r$.
 In particular, this implies that
 if $n\ge N_{r}$
 and
 $\bold u\in\Sigma^{n}$ 
 with
 $UL_{n}(\bold u\bold k)\subseteq C$, then
  $$
 \align
 \dist
 \big(
 \,
 UM_{n}(\bold u\bold l)
 \,,\,
 C
 \,
 \big)
&\le
\distance
\big(
\,
UM_{n}(\bold u\bold l)
\,,\,
UL_{n}(\bold u\bold k)
\,
\big)
+
 \dist
 \big(
 \,
UL_{n}(\bold u\bold k)
 \,,\,
 C
 \,
 \big)\\
 &<
 r+0
 =
 r
  \tag6.13
\endalign
$$
for all
$\bold k,\bold l\in\Sigma_{\smallG}^{\Bbb N}$
with
$\termi(\bold u)=\ini(\bold k)$
and
$\termi(\bold u)=\ini(\bold l)$,
and
if $n\ge N_{r}$
 and
 $\bold u\in\Sigma^{n}$ 
 with
 $UM_{n}(\bold u\bold l)\subseteq C$, then
  $$
 \align
 \dist
 \big(
 \,
 UL_{n}(\bold u\bold k)
 \,,\,
 C
 \,
 \big)
&\le
\distance
\big(
\,
UL_{n}(\bold u\bold k)
\,,\,
UM_{n}(\bold u\bold l)
\,
\big)
+
 \dist
 \big(
 \,
UM_{n}(\bold u\bold l)
 \,,\,
 C
 \,
 \big)\\
 &<
 r+0
 =
 r
  \tag6.14
\endalign
$$
for all
$\bold k,\bold l\in\Sigma_{\smallG}^{\Bbb N}$
with
$\termi(\bold u)=\ini(\bold k)$
and
$\termi(\bold u)=\ini(\bold l)$.
Inclusions (6.11) follow immediately from (6.13),
and inclusions (6.12) follow immediately from (6.14).
\hfill$\square$

\bigskip

\proclaim{Lemma 6.3}
Fix $r>0$.
There is a positive integer $N_{r}$
such that
if $n\ge N_{r}$
and
$C$ be a subset of $\Cal P(\Sigma_{\smallG}^{\Bbb N})$,
then
we have
 $$
 \align
 L_{n}^{-1}(C)
 &\subseteq 
 M_{n}^{-1}(B(C,r))\,,\\
 M_{n}^{-1}(C)
 &\subseteq 
 L_{n}^{-1}(B(C,r))\,.
 \endalign
 $$
\endproclaim
\noindent{\it Proof}\newline
\noindent
This follows by applying the previous lemma to
 $X=\Cal P(\Sigma_{\smallG}^{\Bbb N})$
and the map
$U:\Cal P(\Sigma_{\smallG}^{\Bbb N})\to\Cal P(\Sigma_{\smallG}^{\Bbb N})$ 
defined by 
$U\mu=\mu$ for $\mu\in \Cal P(\Sigma_{\smallG}^{\Bbb N})$.
\hfill$\square$

\bigskip

The final result in this section 
is a simple 
continuity result 
about upper semi-continuous 
maps.
This result is undoubtedly
well-known. However, we have been unable to find 
a reference and for this reason we are including the 
short and simple proof.

\bigskip

\proclaim{Lemma 6.4}
Let $X$ be  a metric space and let $F:X\to\Bbb R$ be an upper 
semi-continuous function.
Let $K_{1},K_{2},\ldots\subseteq X$ be non-empty compact subsets of 
$X$ with
$K_{1}\supseteq K_{2}\supseteq\ldots$. Then
 $$
 \inf_{n}\sup_{x\in K_{n}}F(x)
 =
 \sup_{x\in\bigcap_{n}K_{n}}F(x)\,.
 $$
\endproclaim
\noindent{\it Proof}\newline
First note that it is clear that
$\inf_{n}\sup_{x\in K_{n}}F(x)
 \ge
 \sup_{x\in\cap_{n}K_{n}}F(x)$.
 We will now prove the reverse inequality, namely,
$\inf_{n}\sup_{x\in K_{n}}F(x)
 \le
 \sup_{x\in\cap_{n}K_{n}}F(x)$. 
Let $\varepsilon>0$.
For each $n$, we can choose $x_{n}\in K_{n}$ such that
$F(x_{n})\ge \sup_{x\in K_{n}}F(x)-\varepsilon$.
Next, since $K_{n}$ is compact for all $n$ and
$K_{1}\supseteq K_{2}\supseteq\ldots$, we can find a subsequence $(x_{n_{k}})_{k}$ and a 
point $x_{0}\in\cap_{n}K_{n}$ such that $x_{n_{k}}\to x_{0}$.
Also, since
$K_{n_{1}}\supseteq K_{n_{2}}\supseteq\ldots$, we conclude that
$\sup_{x\in K_{n_{1}}}F(x)
\ge
\sup_{x\in K_{n_{2}}}F(x)
\ge
\ldots$, 
whence
$\inf_{k}\sup_{x\in K_{n_{k}}}F(x)
=
\limsup_{k}\sup_{x\in K_{n_{k}}}F(x)$.
This implies that
$\inf_{n}\sup_{x\in K_{n}}F(x)
\le
\inf_{k}\sup_{x\in K_{n_{k}}}F(x)
=
\limsup_{k}\sup_{x\in K_{n_{k}}}F(x)
\le
\limsup_{k}F(x_{n_{k}})+\varepsilon$.
However, since 
$x_{n_{k}}\to x_{0}$, we deduce from the 
upper semi-continuity of the function $F$, that
$\limsup_{k}F(x_{n_{k}})\le F(x_{0})$.
Consequently
$\inf_{n}\sup_{x\in K_{n}}F(x)
\le
\limsup_{k}F(x_{n_{k}})+\varepsilon
\le
F(x_{0})+\varepsilon
\le
 \sup_{x\in\cap_{n}K_{n}}F(x)+\varepsilon$.
 Finally, letting $\varepsilon\searrow 0$ gives the desired result.
\hfill$\square$

  \bigskip


\heading
{
7. Proofs. The measures $\Pi$ and $\Pi_{n}$}
\endheading

In this section we introduce two technical 
auxiliary measures, namely, the
measures $\Pi$ and $\Pi_{n}$;
see definitions (7.2) and (7.3) below.
The main result is this section 
is Theorem 7.2
showing that the sequence $(\Pi_{n})_{n}$
has the large deviation property.
Theorem 7.2 plays
a
major role in the 
in the proof 
of Theorem 8.3 in Section 8.
We start by introducing some notation.
We first introduce the two
main auxiliary measures
$\Pi$ and $\Pi_{n}$.

\bigskip

%
%

\bigskip

{\bf The measure $\Pi$.}
Let
$B=(b_{\,\smallvertexi,\smallvertexj})_{\smallvertexi,\smallvertexj\in\smallV}$
denote the matrix defined by
 $$
 \align
 b_{\,\smallvertexi,\smallvertexj}
&=
 |\E_{\smallvertexi,\smallvertexj}|\,.
 \endalign
  $$
We denote the spectral radius of $B$ by
$\lambda$.
Since $\G=(\V,\E)$ is strongly connected,
we conclude that the matrix $B$ is irreducible,
and it therefore follows from the Perron-Frobenius theorem that there 
is a unique
right eigen-vector 
$\bold u=(u_{\smallvertexi})_{\smallvertexi\in\smallV}$ 
of $B$ with 
eigen-value $\lambda$
and a 
unique
left eigen-vector
$\bold v=(v_{\smallvertexi})_{\smallvertexi\in\smallV}$ 
of $B$ with 
eigen-value $\lambda$, i\.e\.
$$
 \aligned
 \bold u B
&=
 \lambda\bold u\,,\\
 B\bold v
&=\lambda \bold v\,,
 \endaligned
 \tag7.1
 $$
with
$u_{\smallvertexi},v_{\smallvertexi}>0$
for all $\vertexi$
such that
$\sum_{\smallvertexi}u_{\smallvertexi}v_{\smallvertexi}=1$
and
$\sum_{\smallvertexi}u_{\smallvertexi}=1$.
For $\edge\in\V$, write
 $\pi_{\smalledge}
 =
 v_{\ini(\smalledge)}^{-1}
 \,
 v_{\termi(\smalledge)}
 \,
 \lambda^{-1}$.
A simple calculation shows that
$\sum_{
 \smalledge\in\smallE_{\smallsmallvertexi}
 }
 \pi_{\smalledge} =
 1$
 for all $\vertexi$
and
that
$\sum_{\smallvertexi}
\sum_{
 \smalledge\in\smallE_{\smallsmallvertexi,\smallsmallvertexj}
 }
 u_{\smallvertexi}
 \,
 v_{\smallvertexi}
 \,
 \pi_{\smalledge}
 =
 u_{\smallvertexj}
 \,
 v_{\smallvertexj}$
  for all $\vertexj$.
It follows from this
that there is
 a unique Borel probability measure
 $\Pi\in\Cal P(\Sigma_{\smallG}^{\Bbb N})$
 such that
 $$
 \align
 \Pi[\bold i]
 &=
  u_{\ini(\smalledge_{1})}
 \,
 v_{\ini(\smalledge_{1})}
 \,
  \pi_{\smalledge_{1}}
  \,
  \cdots
  \,
 \pi_{\smalledge_{n}} \\
&=
u_{\ini(\smalledge_{1})}
 \,
 v_{\termi(\smalledge_{n})}
 \,
 \lambda^{-n}\\
&=
u_{\ini(\bold i)}
 \,
 v_{\termi(\bold i)}
 \,
 \lambda^{-n}
 \tag7.2
\endalign 
 $$
for
all
$\bold i
=
\edge_{1}\ldots \edge_{n}\in\Sigma_{\smallG}^{*}$.

\bigskip

{\bf The measure $\Pi_{n}$.}
Finally, for a positive integer $n$, we define the probability measures
$\Pi_{n}\in\Cal P(\Cal P(\Sigma_{\smallG}^{\Bbb N}))$ by
 $$
 \Pi_{n}
 =
 \Pi\circ M_{n}^{-1}\,;
 \tag7.3
 $$
recall that the map
$M_{n}:\Sigma_{\smallG}^{\Bbb N}\to\Cal P_{S}(\Sigma_{\smallG}^{\Bbb N})$ 
is defined 
in (6.1).

\bigskip

We now turn towards the proof of the main result in
  this section, namely, Theorem 7.2 
  showing that the sequence $(\Pi_{n})_{n}$ has the large deviation property.
  The proof of Theorem 7.2
  is based
  on a large deviation theorem by Orey \& Pelikan [OrPe1,OrPe2]
  (see also [JiQiQi, Remark 7.2.2]).
  In particular, Orey \& Pelikan [OrPe1,OrPe2]
  prove that
  various sequences
  of Gibbs measures
  satisfy a large deviation principle, see Theorem 7.1 below.
  However, we begin by 
  stating the definition of the large deviation principle.

\bigskip

\proclaim{Definition}
Let $X$
be a complete separable metric space
and
let $(P_{n})_{n}$ be a sequence of probability measures on $X$. 
Let
$(a_{n})_{n}$ be a sequence of positive numbers with
$a_{n}\to\infty$ 
and let
$I:X\to[0,\infty]$ be a lower semi-continuous function
with compact level sets.
The sequence $(P_{n})_{n}$ is said to have the large deviation property
with constants $(a_{n})_{n}$ and rate function $I$ if the following 
two condistions
hold:
\roster
\item"(i)" For each closed subset $K$ of $X$, we have
 $$
 \limsup_{n}\frac{1}{a_{n}}\log P_{n}(K)\le-\inf_{x\in K}I(x)\,;
 $$
\item"(ii)" For each open subset $G$ of $X$, we have
 $$
 \liminf_{n}\frac{1}{a_{n}}\log P_{n}(G)\ge-\inf_{x\in G}I(x)\,.
 $$
\endroster
\endproclaim

 \bigskip

\proclaim{Theorem 7.1. [OrPe1,OrPe2]}
\roster
\item"(1)"
Let $\Phi:\Sigma_{\smallG}^{\Bbb N}\to\Bbb R$
be a H\"older continous function.
Let $\mu_{\Phi}$ be the 
Gibbs state of $\Phi$ and write
$\mu_{\Phi,n}
 =
 \mu_{\Phi}\circ L_{n}^{-1}$.
Define $I_{\Phi}:\Cal P(\Sigma_{\smallG}^{\Bbb N})\to[0,\infty]$ by
 $$
 I_{\Phi}(\mu)
 =
 \cases
 P(\Phi)-\int\Phi\,d\mu-h(\mu)
&\quad
 \text{for $\mu\in \Cal P_{S}(\Sigma_{\smallG}^{\Bbb N})$;}\\
 \infty
&\quad
 \text{for 
 $\mu\in \Cal P(\Sigma_{\smallG}^{\Bbb N})\setminus \Cal P_{S}(\Sigma_{\smallG}^{\Bbb N})$.
 }
 \endcases
$$
Then the 
sequence
$(\mu_{\Phi,n})_{n}$
has the large deviation deviation property with respect to the sequence
$(n)_{n}$ and rate function $I_{\Phi}$.

\item"(2)"
Write
$\Gamma_{n}
 =
 \Pi\circ L_{n}^{-1}$.
Define $I:\Cal P(\Sigma_{\smallG}^{\Bbb N})\to[0,\infty]$ by
 $$
 I(\mu)
 =
 \cases
 \log\lambda-h(\mu)
\quad
 &\text{for $\mu\in \Cal P_{S}(\Sigma_{\smallG}^{\Bbb N})$;}\\
 \infty
\quad
 &\text{for 
 $\mu\in \Cal P(\Sigma_{\smallG}^{\Bbb N})\setminus \Cal P_{S}(\Sigma_{\smallG}^{\Bbb N})$.
 }
 \endcases
$$
Then the 
sequence
$(\Gamma_{n})_{n}$
has the large deviation deviation property with respect to the sequence
$(n)_{n}$ and rate function $I$.
\endroster
\endproclaim
\noindent{\it Proof}\newline
\noindent
(1) This follows from [OrPe1,OrPe2].

\noindent
(2) 
We use the notation from Part (1), namely, 
if $\Phi:\Sigma_{\smallG}^{\Bbb N}\to\Bbb R$
is a H\"older continous function, then we 
let $\mu_{\Phi}$ denote the 
Gibbs state of $\Phi$ and we write
$\mu_{\Phi,n}$
and
$I_{\Phi}$ for the measure and  function defined 
in the statement of part (1) of the theorem.
Let $\Cal O:\Sigma_{\smallG}^{\Bbb N}\to\Bbb R$ denote the zero-function, i\.e\. 
$\Cal O(\bold i)=0$ for all 
$\bold i\in\Sigma_{\smallG}^{\Bbb N}$.
Noticing that
$\mu_{\Cal O}=\Pi$ 
(by [Wa, Theorem 8.10])
and 
$P(\Cal O)
=
\log\lambda$
(indeed,
the variational principle
implies that
$P(\Cal O)
=
\sup_{\mu\in\Sigma_{\smallsmallG}^{\Bbb N}}h(\mu)$,
and we deduce from [Wa, Theorem 8.10]
that
$\sup_{\mu\in\Sigma_{\smallsmallG}^{\Bbb N}}h(\mu)=h(\Pi)$
and from [Wa, Theorem 7.13]
that
$h(\Pi)=\log\lambda$;
combing these results
show that
$P(\Cal O)=\sup_{\mu\in\Sigma_{\smallsmallG}^{\Bbb N}}h(\mu)=\mu(\Pi)=\log\lambda$),
we conclude that
$\mu_{\Cal O,n}
=
\mu_{\Cal O}\circ L_{n}^{-1}
=
\Pi\circ L_{n}^{-1}
=
\Gamma_{n}$
and $I_{\Cal O}=I$, and the result therefore follows immediately from (1).
\hfill$\square$

\bigskip

We can now state and prove the main result in this section.

\bigskip

\proclaim{Theorem 7.2}
Define $I:\Cal P(\Sigma_{\smallG}^{\Bbb N})\to[0,\infty]$ by
 $$
 I(\mu)
 =
 \cases
 \log\lambda-h(\mu)
\quad
 &\text{for $\mu\in \Cal P_{S}(\Sigma_{\smallG}^{\Bbb N})$;}\\
 \infty
\quad
 &\text{for 
 $\mu\in \Cal P(\Sigma_{\smallG}^{\Bbb N})\setminus \Cal P_{S}(\Sigma_{\smallG}^{\Bbb N})$.
 }
 \endcases
$$
Then the 
sequence
$(\Pi_{n})_{n}$
has the large deviation deviation property with respect to the sequence
$(n)_{n}$ and rate function $I$.
\endproclaim
\noindent{\it Proof}\newline
\noindent
We will use the same notation as in Theorem 7.1, i\.e\. for a positive integer $n$, we write
$\Gamma_{n}=\Pi\circ L_{n}^{-1}$.
We now prove the following two claims.

\bigskip

\noindent{\it Claim 1. For each closed subset
$K$ of $ \Cal P(\Sigma_{\smallG}^{\Bbb N})$,
we have
$$
\limsup_{n}\frac{1}{n}\log \Pi_{n}(K)
\le
-
\inf_{\mu\in K}I(\mu)\,.
$$
}

\noindent{\it Proof of Claim 1.}
Let $K$ be a
close
 subset of
$\Cal P(\Sigma_{\smallG}^{\Bbb N}) $.
Now,
observe that 
it follows from Lemma 6.3
that
for each $r>0$, we can find a positive integer
$N_{r}$ such that
for all $n\ge N_{r}$
we have
 $$
 \Big\{
 \bold i\in\Sigma_{\smallG}^{\Bbb N}
 \,\Big|\,
 M_{n}\bold i\in K
\Big\}
\subseteq
\Big\{
 \bold i\in\Sigma_{\smallG}^{\Bbb N}
 \,\Big|\,
 L_{n}\bold i\in B(K,r)
\Big\}\,.
$$
This clearly implies that 
for all $n\ge N_{r}$
we have
 $$
 \align
 \Pi_{n}(K)
 &\le
 \Gamma_{n}\big(\,B(K,r)\,\big)\\
  &\le
 \Gamma_{n}\big(\,\overline{B(K,r)}\,\big)\,.
 \tag7.4
 \endalign
 $$
Next, since $\overline{B(K,r)}$
is a closed set 
and
since it follows from Theorem 7.1
that
the sequence $(\Gamma_{n})_{n}$ has the 
large deviation property with respect to the sequence $(n)_{n}$
and rate function $I$, we conclude
 from (7.4) that
 $$
\align
\limsup_{n}\frac{1}{n}\log \Pi_{n}(K)
&
\le
\limsup_{n}\frac{1}{n}\log \Gamma_{n}(\, \overline{B(K,r)}\,)\\
&\le
-
\inf_{\mu\in  \overline{B(K,r)}}I(\mu)
\endalign
$$ 
for all $r>0$. Hence
 $$
\align
\limsup_{n}\frac{1}{n}\log \Pi_{n}(K)
&\le
\inf_{r>0}
\Bigg(
-
\inf_{\mu\in  \overline{B(K,r)}}I(\mu)\Bigg)\\
&=
\inf_{r>0}
\,
\sup_{\mu\in  \overline{B(K,r)}}-I(\mu)\\
&=
\inf_{n}
\,
\sup_{\mu\in  \overline{B(K,\frac{1}{n})}}-I(\mu)\,.
\tag7.5
\endalign
$$

However, 
it follows from
Theorem 7.1 that $I$ is a rate function and therefore, in particular, 
lower
semi-continuous.
We deduce from this
that the function
$-I$
is upper semi-continuous.
Since the sets 
$K_{n}=\overline{B(K,\frac{1}{n})}$
are compact 
(because they are closed subsets of
$\Cal P(\Sigma_{\smallG}^{\Bbb N})$
and $\Cal P(\Sigma_{\smallG}^{\Bbb N})$ is compact)
with
$K_{1}\supseteq K_{2}\supseteq K_{3}\ldots$,
we therefore conclude from  Lemma 6.4 and the fact that $-I$ is upper
semi-continuous  that
 $$
 \align
 \inf_{n}
\,
\sup_{\mu\in  \overline{B(K,\frac{1}{n})}}-I(\mu)
&=
\inf_{n}
\,
\sup_{\mu\in K_{n}}-I(\mu)\\
&=
\sup_{\mu\in  \bigcap_{n}K_{n}}-I(\mu)
\tag7.6
\endalign
$$
Next, observing that 
$\cap_{n}K_{n}
=
\cap_{n}\overline{B(K,\frac{1}{n})}
=
K$
(because $K$ is closed), we deduce from (7.6) that
 $$
 \align
 \inf_{n}
\,
\sup_{\mu\in  \overline{B(K,\frac{1}{n})}}-I(\mu)
&=
\sup_{\mu\in K}-I(\mu)\\
&=
-\inf_{\mu\in K}I(\mu)\,.
\tag7.7
\endalign
$$

Finally. combining
(7.5) and (7.7) shows that
$\limsup_{n}\frac{1}{n}\log \Pi_{n}(K)
\le
-\inf_{\mu\in K}I(\mu)$.
This completes the proof of Claim 1.

\bigskip

\noindent{\it Claim 2. For each open subset
$G$ of $ \Cal P(\Sigma_{\smallG}^{\Bbb N})$,
we have
$$
\liminf_{n}\frac{1}{n}\log \Pi_{n}(G)
\ge
-
\inf_{\mu\in G}I(\mu)\,.
$$
}

\noindent{\it Proof of Claim 2.}
Let $G$ be an
open subset of
$\Cal P(\Sigma_{\smallG}^{\Bbb N}) $.
For each $r>0$, we will  write
$I(G,r)
=
\{\mu\in\Cal P(\Sigma_{\smallG}^{\Bbb N})\,|\,
\dist(\mu,\Cal P(\Sigma_{\smallG}^{\Bbb N})\setminus G)>r\}$.
Now,
observe that 
it follows from Lemma 6.3
that
for each $r>0$, we can find a positive integer
$N_{r}$ such that
for all $n\ge N_{r}$
we have
 $$
 \Big\{
 \bold i\in\Sigma_{\smallG}^{\Bbb N}
 \,\Big|\,
 L_{n}\bold i\in I(G,r)
\Big\}
\subseteq
\Big\{
 \bold i\in\Sigma_{\smallG}^{\Bbb N}
 \,\Big|\,
 M_{n}\bold i\in B(\,I(G,r)\,,\,\tfrac{r}{2}\,)
\Big\}\,.
$$
This clearly implies that 
for all $n\ge N_{r}$
we have
 $$
 \Gamma_{n}(\,I(G,r)\,)
 \le
 \Pi_{n}\big(\,B(\,I(G,r)\,,\,\tfrac{r}{2}\,)\,\big)
 \tag7.8
 $$
However, since it is easily seen that
$B(\,I(G,r)\,,\,\frac{r}{2})
\subseteq
G$, it follows that
$\Pi_{n}\big(\,B(\,I(G,r)\,,\,\frac{r}{2})\,\big)
\le
\Pi_{n}(G)$, and (7.8) 
therefore shows that
for all $n\ge N_{r}$
we have
 $$
 \Gamma_{n}(\,I(G,r)\,)
 \le
 \Pi_{n}(G)\,.
 \tag7.9
 $$
Next, since $I(G,r)$ is an open set
and
since it follows from Theorem 7.1
that
the sequence $(\Gamma_{n})_{n}$ has the 
large deviation property with respect to the sequence $(n)_{n}$
and rate function $I$, we conclude
from (7.9)
that
 $$
\align
\liminf_{n}\frac{1}{n}\log \Pi_{n}(G)
&
\ge
\liminf_{n}\frac{1}{n}\log \Gamma_{n}(\,I(G,r)\,)\\
&\ge
-
\inf_{\mu\in I(G,r)}I(\mu)
\endalign
$$ 
for all $r>0$. Hence
 $$
\align
\liminf_{n}\frac{1}{n}\log \Pi_{n}(G)
&\ge
\sup_{r>0}
\Bigg(
-
\inf_{\mu\in I(G,r)}I(\mu)
\Bigg)\\
&=
-
\inf_{r>0}
\,
\inf_{\mu\in I(G,r)}I(\mu)\\
&=
-
\inf_{\mu\in \bigcup_{r>0}I(G,r)}I(\mu)\,.
\tag7.10
\endalign
$$ 
Finally, since $G$ is open, it follows easily that
$\cup_{r>0}I(G,r)=G$, and (7.10) therefore implies that
$\liminf_{n}\frac{1}{n}\log \Pi_{n}(G)
\ge
-
\inf_{\mu\in \cup_{r>0}I(G,r)}I(\mu)
=
-
\inf_{\mu\in G}I(\mu)$. 
This completes the proof of Claim 2.

\bigskip

the desired result follows immediately from Claim 1 and Claim 2.
\hfill$\square$

  \bigskip


%

%

\heading
{
8. Proofs.  The modified multifractal pressure}
\endheading

In this section we introduce our main technical tool, namely, the modified 
multifractal 
pressure;
see definition (8.1) below.
The two main results in 
this section
are Theorem 8.3
and Theorem 8.4.
Theorem 8.3
applies Theorem 7.2
(saying
that
 the sequence $(\Pi_{n})_{n}$
 satisfies the large deviation principle)
 to establish
 a 
variational principle for the modified multifractal pressure
 and Theorem 8.4
shows that
the 
 multifractal pressure and the 
 modified
multifractal
pressure are \lq asymptotically"  the same.

We start be defined the 
modified multifractal pressure.
 For $C\subseteq X$,
we define 
 the 
 modified lower and upper 
  mutifractal pressure
of 
 $\varphi$
associated with the space $X$ and the map $U$ and  by
 $$
 \aligned
  \underline Q_{C}^{U}(\varphi)
 &=
 \,
  \liminf_{n}
 \,
  \,\,
  \frac{1}{n}
 \,\,
 \log
   \sum
  \Sb
  \bold i\in\Sigma_{\smallsmallG}^{n}\\
  {}\\
  UM_{n}[\bold i]\subseteq C
  \endSb
  \sup_{\bold u\in[\bold i]}
  \,\,
   \exp
 \,\,
 \sum_{k=0}^{n-1}\varphi S^{k}\bold u\,,\\
   \overline Q_{C}^{U}(\varphi)
 &=
  \limsup_{n}
  \,\,
  \frac{1}{n}
 \,\,
 \log
   \sum
  \Sb
 \bold i\in\Sigma_{\smallsmallG}^{n}\\
  {}\\
  UM_{n}[\bold i]\subseteq C
  \endSb
   \sup_{\bold u\in[\bold i]}
   \,\,
   \exp
 \,\,
 \sum_{k=0}^{n-1}\varphi S^{k}\bold u\,.
 \endaligned
 \tag8.1
 $$

In order to establish 
a variational principle for the modified multifractal pressure, we introduce the following notation.
 For a continuous function $\varphi:\Sigma_{\smallG}^{\Bbb N}\to\Bbb R$,
 we
define $F_{\varphi}:\Cal P_{S}(\Sigma_{\smallG}^{\Bbb N})\to\Bbb R$ by
 $$
 F_{\varphi}(\mu)=\int\varphi\,d\mu\,.
 $$
Observe that since $\varphi$ is bounded, 
i\.e\. $\|\varphi\|_{\infty}<\infty$,
we conclude that
$\|F_{\varphi}\|_{\infty}
\le
\|\varphi\|_{\infty}<\infty$.
In addition, for a positive integer $n$, define
 the
probability measure
$Q_{\varphi,n}\in\Cal P\big(\,\Cal P(\Sigma_{\smallG}^{\Bbb N})\,\big)$ by
 $$
 \align
 Q_{\varphi,n}(E)
&=
 \frac
 {\int_{E}\exp(nF_{\varphi})\,d\Pi_{n}}
 {\int\exp(nF_{\varphi})\,d\Pi_{n}}
 \quad
 \text{for Borel subsets $E$ of $\Cal P(\Sigma_{\smallG}^{\Bbb N})$;}
 \endalign
 $$
recall, that the measure $\Pi_{n}$ is defined in (7.3).

We now turn towards the proof of the first main result in
  this section, namely, Theorem 8.2 providing 
  a variational principle for the modified multifractal pressure.
The proof of Theorem 8.2 is based on
large deviation theory.
In particular, 
we will
use the fact
that
the sequence $(\Pi_{n})_{n}$
 satisfies the large deviation principle
 together
with Varadhan's [Va]
large deviation theorem (Theorem 8.1.(1) below),
and a non-trivial application of this (namely Theorem 8.1.(2) below) 
providing first order asymptotics of certain
\lq\lq Boltzmann distributions".
Recall, that 
the notion of a large deviation principle is defined in Section 7.


%

\bigskip

\proclaim{Theorem 8.1}
Let $X$ be a complete separable metric space
and
let $(P_{n})_{n}$ be a sequence of probability measures on $X$.
Assume that the sequence $(P_{n})_{n}$
has the large deviation property
with constants $(a_{n})_{n}$ and rate function $I$.
Let $F:X\to\Bbb R$ be a continuous function
satisfying the following two conditions:
\roster
\item"(i)"
For all $n$, we have
 $$
 \int\exp(a_{n}F)\,dP_{n}
 <
 \infty\,.
 $$
\item"(ii)"
We have
 $$
 \lim_{M\to\infty}\,\,
 \limsup_{n}\,\,
 \frac{1}{a_{n}}
 \log
 \int_{\{M\le F\}}
 \exp(a_{n}F)\,dP_{n}
 =
 -\infty\,.
 $$
 \endroster
(Observe that the  Conditions (i)--(ii)
 are satisfied if $F$ is bounded.)
Then the following statements hold.
\roster
\item"(1)" We have
 $$
 \lim_{n}\,\,
 \frac{1}{a_{n}}
 \log
 \int
 \exp(a_{n}F)\,dP_{n}
 =
 -\inf_{x\in X}(I(x)-F(x))\,.
 $$
\item"(2)"
For each $n$ define a probability measure $Q_{n}$ on $X$ by
 $$
 Q_{n}(E)
 =
 \,\,
 \frac
 {\int_{E}\exp(a_{n}F)\,dP_{n}}
 {\int\exp(a_{n}F)\,dP_{n}}\,.
 $$
Then the sequence $(Q_{n})_{n}$
has the large deviation property with constants $(a_{n})_{n}$
and rate function
$(I-F)-\inf_{x\in X, I(x)<\infty}(I(x)-F(x))$.
\endroster
\endproclaim 
\noindent{\it Proof}\newline
\noindent 
Statement (1) follows from [El, Theorem II.7.1] or [DeZe, Theorem 4.3.1], and
statement (2) follows from [El, Theorem II.7.2].
\hfill$\square$

 \bigskip

 \noindent
 Before stating and proving Theorem 8.3, we establish the following auxiliary result.

\bigskip

\proclaim{Theorem 8.2}
Let $X$ be a metric space and let $U:\Cal P(\Sigma_{\smallG}^{\Bbb N})\to X$ be 
continuous with respect to the weak topology.
Let $C\subseteq X$ be a  subset of $X$.
Fix a continuous function $\varphi:\Sigma_{\smallG}^{\Bbb N}\to\Bbb R$.
Then there is a constant $c$
such that
for all positive integers $n$, we have
 $$
 \align
 \sum
   \Sb
   \bold k\in\Sigma_{\smallsmallG}^{n}\\
   {}\\
   UM_{n}[\bold k]\subseteq C
	 \endSb
 \sup_{\bold u\in[\bold k]}
 \,\,
 \exp
 \sum_{i=0}^{n-1}\varphi S^{i}\bold u
&\le
c
\,\,
\lambda^{n}
\,\,
Q_{\varphi,n}\Big(\{U\in C\}\Big)
\,\,
\int\exp(nF_{\varphi})\,d\Pi_{n}\,,\\
{}\\
 \sum
   \Sb
    \bold i\in\Sigma_{\smallsmallG}^{n}\\
   {}\\
   UM_{n}[\bold k]\subseteq C
	 \endSb
 \sup_{\bold u\in[\bold k]}
 \,\,	 
 \exp
 \sum_{i=0}^{n-1}\varphi S^{i}\bold u
&\ge
\frac{1}{c}
\,\,
\lambda^{n}
\,\,
Q_{\varphi,n}\Big(\{U\in C\}\Big)
\,\,
\int\exp(nF_{\varphi})\,d\Pi_{n}\,.
\endalign
 $$
\endproclaim
\noindent{\it Proof}
\newline
\noindent 
For each positive integer $n$ and each 
$ \bold i\in\Sigma_{\smallG}^{n}$, 
we write
$s_{\bold i}
=
\sup_{\bold u\in[\bold i]}
\exp
 \sum_{k=0}^{n-1}\varphi S^{k}\bold u$
 for sake of brevity.
Let $C$ be a subset of $X$.
For each positive integer $n$, we clearly have
 $$
 \align
 &\int
 \limits_{
 \big\{
 \bold j\in\Sigma_{\smallsmallG}^{\Bbb N}
 \,\big|\,
 U M_{n}[\bold j|n]\subseteq C
 \big\}
 }
s_{\bold i|n}
 \,d\Pi(\bold i)\\
&\qquad\qquad
=
\qquad\qquad
\,\,\,\,
 \sum_{ \bold k\in\Sigma_{\smallsmallG}^{n}}
 \,\,\,\,
 \int\limits_{
 [\bold k]
 \,\,\cap\,\,
 \big\{
 \bold j\in\Sigma_{\smallsmallG}^{\Bbb N}
 \,\big|\,
 U M_{n}[\bold j|n]\subseteq C
  \big\}
 }
 s_{\bold i|n}
 \,d\Pi(\bold i)\\
\allowdisplaybreak 
&\qquad\qquad
=
\qquad\qquad
\,\,\,\,
 \sum_{ \bold k\in\Sigma_{\smallsmallG}^{n}}
 s_{\bold k}
 \,\,
 \Pi
 \Big(
 \,
  [\bold k]
 \,\,\cap\,\,
 \Big\{
 \bold j\in\Sigma_{\smallG}^{\Bbb N}
 \,\Big|\,
 U M_{n}[\bold j|n]\subseteq C
 \Big\}
 \,
 \Big)\\
\allowdisplaybreak 
&\qquad\qquad
=
 \sum
   \Sb
    \bold k\in\Sigma_{\smallsmallG}^{n}\\
   {}\\
   [\bold k]
 \,\,\cap\,\,
 \big\{
 \bold j\in\Sigma_{\smallsmallG}^{\Bbb N}
 \,\big|\,
 U M_{n}[\bold j|n]\subseteq C
 \big\}
 \not=\varnothing
	 \endSb
 s_{\bold k}
 \,\,
 \Pi
 \Big(
 \,
  [\bold k]
 \,\,\cap\,\,
 \Big\{
 \bold j\in\Sigma_{\smallG}^{\Bbb N}
 \,\Big|\,
 U M_{n}[\bold j|n]\subseteq C
  \Big\}
 \,
 \Big)\,.\\
& \tag8.3
 \endalign
 $$
Now observe that if
$\bold k\in\Sigma_{\smallG}^{n}$ with 
 $[\bold k]
 \cap
 \{
 \bold j\in\Sigma_{\smallG}^{\Bbb N}
 \,|\,
 U M_{n}[\bold j|n]\subseteq C
 \}
 \not=\varnothing$,
 then
 a simple argument shows  that
 $U M_{n}[\bold k]
 \subseteq 
 C$. 
 We conclude from this  that
 $$
 \align
& \sum
   \Sb
   \bold k\in\Sigma_{\smallsmallG}^{n}\\
   {}\\
   [\bold k]
 \,\,\cap\,\,
 \big\{
 \bold j\in\Sigma_{\smallsmallG}^{\Bbb N}
 \,\big|\,
 U M_{n}[\bold j|n]\subseteq C
 \big\}
 \not=\varnothing
	 \endSb
 s_{\bold k}
 \,\,
 \Pi
 \Big(
 \,
  [\bold k]
 \,\,\cap\,\,
 \Big\{
 \bold j\in\Sigma_{\smallG}^{\Bbb N}
 \,\Big|\,
 U M_{n}[\bold j|n]\subseteq C
  \Big\}
 \,
 \Big)\\
 &\qquad\qquad
 =
 \sum
   \Sb
   \bold k\in\Sigma_{\smallsmallG}^{n}\\
   {}\\
   [\bold k]
 \,\,\cap\,\,
 \big\{
 \bold j\in\Sigma_{\smallsmallG}^{\Bbb N}
 \,\big|\,
 U M_{n}[\bold j|n]\subseteq C
 \big\}
 \not=\varnothing\\
 {}\\
 UM_{n}[\bold k]\subseteq C
	 \endSb
 s_{\bold k}
 \,\,
 \Pi
 \Big(
 \,
  [\bold k]
 \,\,\cap\,\,
 \Big\{
 \bold j\in\Sigma_{\smallG}^{\Bbb N}
 \,\Big|\,
 U M_{n}[\bold j|n]\subseteq C
  \Big\}
 \,
 \Big)\,.\\
 \tag8.4
 \endalign
 $$
Combining (8.3) and (8.4)  gives
 $$
 \align
&\int
 \limits_{
 \big\{
 \bold j\in\Sigma_{\smallsmallG}^{\Bbb N}
 \,\big|\,
 U M_{n}[\bold j|n]\subseteq C
 \big\}
 }
s_{\bold i|n}
 \,d\Pi(\bold i)\\
&\qquad\qquad
=
 \sum
   \Sb
  \bold k\in\Sigma_{\smallsmallG}^{n}\\
   {}\\
   [\bold k]
 \,\,\cap\,\,
 \big\{
 \bold j\in\Sigma_{\smallsmallG}^{\Bbb N}
 \,\big|\,
 U M_{n}[\bold j|n]\subseteq C
 \big\}
 \not=\varnothing\\
  {}\\
 UM_{n}[\bold k]\subseteq C
	 \endSb
 s_{\bold k}
 \,\,
 \Pi
 \Big(
 \,
  [\bold k]
 \,\,\cap\,\,
 \Big\{
 \bold j\in\Sigma_{\smallG}^{\Bbb N}
 \,\Big|\,
 U M_{n}[\bold j|n]\subseteq C
  \Big\}
 \,
 \Big)\,.\\
 &\qquad\qquad
=
\qquad\quad\,\,\,\,
 \sum
   \Sb
   \bold k\in\Sigma_{\smallsmallG}^{n}\\
  {}\\
 UM_{n}[\bold k]\subseteq C
	 \endSb
 s_{\bold k}
 \,\,
 \Pi
 \Big(
 \,
  [\bold k]
 \,\,\cap\,\,
 \Big\{
 \bold j\in\Sigma_{\smallG}^{\Bbb N}
 \,\Big|\,
 U M_{n}[\bold j|n]\subseteq C
  \Big\}
 \,
 \Big)\,.\\
 \tag8.5
 \endalign
 $$

However,
if
$\bold k\in\Sigma_{\smallG}^{n}$ with
  $ U M_{n}[\bold k]\subseteq C$,
  then  it is clear that
  $[\bold k]\subseteq
 \{
 \bold j\in\Sigma_{\smallG}^{\Bbb N}
 \,|\,
  U M_{n}[\bold j|n]\subseteq C
 \}$,
 whence
 $[\bold k]\cap
 \{
 \bold j\in\Sigma_{\smallG}^{\Bbb N}
 \,|\,
  U M_{n}[\bold j|n]\subseteq C
 \}
 =
 [\bold k]$.
 This and (8.5) now imply that
 $$
 \align
  \int
 \limits_{
 \big\{
 \bold j\in\Sigma_{\smallsmallG}^{\Bbb N}
 \,\big|\,
 U M_{n}[\bold j|n]\subseteq C
 \big\}
 }
&s_{\bold i|n}
 \,d\Pi(\bold i)\\
&=
 \sum
   \Sb
  \bold k\in\Sigma_{\smallsmallG}^{n}\\
   {}\\
   U M_{n}[\bold k]\subseteq C
	 \endSb
 s_{\bold k}
 \,\,
 \Pi
 \Big(
 \,
  [\bold k]
 \,\,\cap\,\,
 \Big\{
 \bold j\in\Sigma_{\smallG}^{\Bbb N}
 \,\Big|\,
 U M_{n}[\bold j|n]\subseteq C
  \Big\}
 \,               
 \Big)\\ 
&=
 \sum
   \Sb
  \bold k\in\Sigma_{\smallsmallG}^{n}\\
   {}\\
 U M_{n}[\bold k]\subseteq C
	 \endSb
 s_{\bold k}
 \,\,
 \Pi\big(\,
 [\bold k]
 \,\big)\\
\allowdisplaybreak 
&=
 \sum
   \Sb
  \bold k\in\Sigma_{\smallsmallG}^{n}\\
  {}\\
   U M_{n}[\bold k]\subseteq C
	 \endSb
 s_{\bold k}
 \,
 u_{\ini(\bold k)}
 \,
 v_{\termi(\bold k)}
 \,
 \lambda^{-n}\,.
 \tag8.6
 \endalign
 $$

It follows from the 
Principle of
Bounded Distortion 
(see, for example, [Bar,Fa])
that there is a constant $c_{0}>0$
such that
if $n\in\Bbb N$, $\bold i\in\Sigma_{\smallG}^{n}$
and $\bold u,\bold v\in[\bold i]$, then
$\frac{1}{c}
\le
\frac
{
\exp\sum_{k=0}^{n-1}\varphi S^{k}\bold u
}
{
\exp\sum_{k=0}^{n-1}\varphi S^{k}\bold v
}
\le
c$.
In particular, this implies that
for all 
$n\in\Bbb N$ and for all $\bold i\in\Sigma_{\smallG}^{n}$, we have
 $$
 \frac{1}{c_{0}}
\exp
\sum_{k=0}^{n-1}\varphi S^{k}\overline{\bold i}
\le
s_{\bold i}
\le
c_{0}
\exp
\sum_{k=0}^{n-1}\varphi S^{k}\overline{\bold i}\,.
\tag8.7
$$
We can also find a constant $c_{1}>0$
such that
$\frac{1}{c_{1}}
\le
u_{\smallvertexi}
v_{\smallvertexi}
\le
c_{1}$
for all $\vertexi$.
Now put $c=c_{0}c_{1}$.

\medskip

{\it Claim 1. For all positive integers $n$, we have
 $$
 \align
  \sum
   \Sb
\bold k\in\Sigma_{\smallsmallG}^{n}\\
  {}\\
   U M_{n}[\bold k]\subseteq C
	 \endSb
 s_{\bold k}
&\le
\,
 c
 \,
 \lambda^{n}
 \int\limits_{
 \big\{
 \bold j\in\Sigma_{\smallsmallG}^{\Bbb N}
 \,\big|\,
 U M_{n}[\bold j|n]\subseteq C
 \big\}
 }
 \exp\left(n F_{\varphi}(M_{n}\bold i)\right)
 \,d\Pi(\bold i)\,,
 \tag8.8\\
  \sum
   \Sb
\bold k\in\Sigma_{\smallsmallG}^{n}\\
  {}\\
   U M_{n}[\bold k]\subseteq C
	 \endSb
 s_{\bold k}
&\ge
 \frac{1}{c}
 \,
 \lambda^{n}
 \int\limits_{
 \big\{
 \bold j\in\Sigma_{\smallsmallG}^{\Bbb N}
 \,\big|\,
 U M_{n}[\bold j|n]\subseteq C
 \big\}
 }
 \exp\left(n F_{\varphi}(M_{n}\bold i)\right)
 \,d\Pi(\bold i)\,.
 \tag8.9
 \endalign
 $$
}

\noindent
{\it Proof of Claim 1.}
It follows from (8.6) and (8.7) that
if $n$ is a positive integer, then
we have
 $$
 \align
  \sum
   \Sb
\bold k\in\Sigma_{\smallsmallG}^{n}\\
  {}\\
   U M_{n}[\bold k]\subseteq C
	 \endSb
 s_{\bold k}
&\le 
 c
 \,
 \lambda^{n}
 \int\limits_{
 \big\{
 \bold j\in\Sigma_{\smallsmallG}^{\Bbb N}
 \,\big|\,
U M_{n}[\bold j|n]\subseteq C
 \big\}
 }
 \exp
 \Bigg(
 \sum_{k=0}^{n-1}\varphi S^{k}\left(\,\overline{\bold i|n}\,\right)
 \Bigg)
 \,d\Pi(\bold i)\\
\allowdisplaybreak 
&=
 c
 \,
 \lambda^{n}
 \int\limits_{
 \big\{
 \bold j\in\Sigma_{\smallsmallG}^{\Bbb N}
 \,\big|\,
U M_{n}[\bold j|n]\subseteq C
 \big\}
  }
 \exp
 \Bigg(
 n\int\varphi\,d(M_{n}\bold i)
 \Bigg)
 \,d\Pi(\bold i)\\
\allowdisplaybreak 
&=
 c
 \,
 \lambda^{n}
 \int\limits_{
 \big\{
 \bold j\in\Sigma_{\smallsmallG}^{\Bbb N}
 \,\big|\,
 U M_{n}[\bold j|n]\subseteq C
 \big\}
 }
 \exp\left(n F_{\varphi}(M_{n}\bold i)\right)
 \,d\Pi(\bold i)\,.
 \endalign
 $$
This proves inequality (8.8).
Inequality (8.9) is proved similarly.
This completes the proof of Claim 1.

\medskip

{\it Claim 2.
For all positive integers $n$, we have
$\{
 \bold j\in\Sigma_{\smallG}^{\Bbb N}
 \,|\,
 U M_{n}[\bold j|n]\subseteq C
 \}
=
\{
 \bold j\in\Sigma_{\smallG}^{\Bbb N}
 \,|\,
 U M_{n}\bold j\subseteq C
 \}$.
}
 
\noindent
{\it Proof of Claim 2.}
This is easily seen; however, for 
the convenience of the reader we have decided to state the result explicitly.
 This completes the proof of Claim 2.

 \medskip

 For all positive integers $n$,
 we now  deduce from 
 Claim 1 and Claim 2
 that
 $$
 \align
  \sum
   \Sb
\bold k\in\Sigma_{\smallsmallG}^{n}\\
  {}\\
   U M_{n}[\bold k]\subseteq C
	 \endSb
 s_{\bold k}
&\le
 c
 \,
 \lambda^{n}
 \int\limits_{
 \big\{
 \bold j\in\Sigma_{\smallsmallG}^{\Bbb N}
 \,\big|\,
 U M_{n}[\bold j|n]\subseteq C
 \big\}
 }
 \exp\left(n F_{\varphi}(M_{n}\bold i)\right)
 \,d\Pi(\bold i)\\ 
&=
 c
 \,
 \lambda^{n}
 \int\limits_{
 \big\{U M_{n}\in C\big\}
 }
 \exp\left(nF_{\varphi}(M_{n}\bold i)\right)
 \,d\Pi(\bold i)\\
\allowdisplaybreak 
&=
c
 \,
 \lambda^{n}
 \int\limits_{
 \big\{U\in C\big\}
 }
 \exp\left(nF_{\varphi}\right)
 \,d\Pi_{n}\\
\allowdisplaybreak 
&=
 c
 \,
 \lambda^{n}
 \,\,
 Q_{\varphi,n}\Big(\{U\in C\}\Big)
 \,\,
 \int\exp\left(nF_{\varphi}\right)\,d\Pi_{n}\,.
 \endalign
 $$

Similarly, we prove that
for all positive integers $n$, we have
  $$
 \align
  \sum
   \Sb
\bold k\in\Sigma_{\smallsmallG}^{n}\\
  {}\\
   U M_{n}[\bold k]\subseteq C
	 \endSb
 s_{\bold k}
&\ge
\frac{1}{c}
 \,
 \lambda^{n}
 \,\,
 Q_{\varphi,n}\Big(\{U\in C\}\Big)
 \,\,
 \int\exp\left(nF_{\varphi}\right)\,d\Pi_{n}\,.
 \endalign
 $$
This completes the proof of Theorem 8.2.
\hfill$\square$

\bigskip

\proclaim{Theorem 8.3. 
The variational principle for the modified
multifractal pressure}
Let $X$ be a metric space and let $U:\Cal P(\Sigma_{\smallG}^{\Bbb N})\to X$ be 
continuous with respect to the weak topology.
Fix a continuous function $\varphi:\Sigma_{\smallG}^{\Bbb N}\to\Bbb R$.
\roster
\item"(1)"
If $G$ is an open subset of $X$
with
$U^{-1}G\cap\Cal P_{S}(\Sigma_{\smallG}^{\Bbb N})\not=\varnothing$, then
 $$
\underline Q_{G}^{U}(\varphi)
\ge
\sup
\Sb
\mu\in\Cal P_{S}(\Sigma_{\smallsmallG}^{\Bbb N})\\
{}\\
U\mu\in G
\endSb
\Bigg(
h(\mu)+\int\varphi\,d\mu
\Bigg)\,.
\tag8.10
$$
\item"(2)"
If $K$ is a closed subset of $X$
with
$U^{-1}K\cap\Cal P_{S}(\Sigma_{\smallG}^{\Bbb N})\not=\varnothing$, then
$$
\overline Q_{K}^{U}(\varphi)
\le
\sup
\Sb
\mu\in\Cal P_{S}(\Sigma_{\smallsmallG}^{\Bbb N})\\
{}\\
U\mu\in K
\endSb
\Bigg(
h(\mu)+\int\varphi\,d\mu
\Bigg)\,.
\tag8.11
$$ 
\endroster
\endproclaim

\noindent{\it Proof}\newline
\noindent
We introduce the 
simplified notation from the proof of Theorem 8.2, i\.e\.
for each positive integer $n$ and each $\bold i\in\Sigma_{\smallG}^{n}$,
we write
$s_{\bold i}
=
\sup_{\bold u\in[\bold i]}
\exp
 \sum_{k=0}^{n-1}\varphi S^{k}\bold u$.
First
note that it follows immediately from Theorem 8.2 that
 $$
 \aligned
 \liminf_{n}
 \frac{1}{n}
 \log
 \sum
   \Sb
\bold i\in\Sigma_{\smallsmallG}^{n}\\
   {}\\
   U M_{n}[\bold i]\subseteq G
	 \endSb
 s_{\bold i}
&\ge
 \log \lambda
 \,+\,
 \liminf_{n}
 \frac{1}{n}
 \log Q_{\varphi,n}\Big(\{U\in G \}\Big)\\
&\qquad\qquad
  \qquad\qquad
 \,+\,
 \liminf_{n}
 \frac{1}{n}
 \log 
 \int
 \exp
 \left(
 nF_{\varphi}
 \right)
 \,d\Pi_{n}\,,\\
\limsup_{n}
 \frac{1}{n}
 \log
 \sum
   \Sb
\bold i\in\Sigma_{\smallsmallG}^{n}\\
   {}\\
   U M_{n}[\bold i]\subseteq K
	 \endSb
 s_{\bold i}
&\le
 \log \lambda
 \,+\,
 \limsup_{n}
 \frac{1}{n}
 \log Q_{\varphi,n}\Big(\{U\in K \}\Big)\\
&\qquad\qquad
  \qquad\qquad
 \,+\,
 \limsup_{n}
 \frac{1}{n}
 \log 
 \int
 \exp
 \left(
 nF_{\varphi}
 \right)
 \,d\Pi_{n}\,.
 \endaligned
 \tag8.12
 $$

Define $I:\Cal P(\Sigma_{\smallG}^{\Bbb N})\to[0,\infty]$ by
 $$
 I(\mu)
 =
 \cases
 \log\lambda-h(\mu)
\quad
 &\text{for $\mu\in \Cal P_{S}(\Sigma_{\smallG}^{\Bbb N})$;}\\
 \infty
\quad
 &\text{for 
 $\mu\in \Cal P(\Sigma_{\smallG}^{\Bbb N})\setminus \Cal P_{S}(\Sigma_{\smallG}^{\Bbb N})$.
 }
 \endcases
$$

Next, we observe that it follows from Theorem 7.2
that the sequence 
$(\Pi_{n})_{n}
\subseteq
\Cal P\big(\,\Cal P(\Sigma_{\smallG}^{\Bbb N})\,\big)$ 
has the large deviation property with respect to
the sequence
$(n)_{n}$ and rate function $I$.
We therefore 
conclude from Part (1) of Theorem 8.1 that
 $$
 \align
  \lim_{n}
 \frac{1}{n}
 \log 
 \int
 \exp
 \left(
 nF_{\varphi}
 \right)
 \,d\Pi_{n}
&=
 -
 \,
  \inf_{\nu\in\Cal P(\Sigma_{\smallsmallG}^{\Bbb N})}(I(\nu)-F_{\varphi}(\nu))\\
&=
 -
  \inf_{\nu\in\Cal P_{S}(\Sigma_{\smallsmallG}^{\Bbb N})}(I(\nu)-F_{\varphi}(\nu))\,.
  \tag8.13
  \endalign
 $$

 Note that it follows that
 $$
 (I-F_{\varphi})
-
\inf
\Sb
\nu\in\Cal P(\Sigma_{\smallsmallG}^{\Bbb N})\\
{}\\
I(\nu)<\infty
\endSb
(I(\nu)-F_{\varphi}(\nu))
=
(I-F_{\varphi})
-
\inf_{\nu\in\Cal P_{S}(\Sigma_{\smallsmallG}^{\Bbb N})}(I(\nu)-F_{\varphi}(\nu))\,;
\tag8.14
$$
indeed,
if 
$\nu\in\Cal P(\Sigma_{\smallsmallG}^{\Bbb N})$,
then
it is clear
from the definition of $I$
that
$I(\nu)=\infty$
if and only if
$\nu\in \Cal P(\Sigma_{\smallG}^{\Bbb N})\setminus \Cal P_{S}(\Sigma_{\smallG}^{\Bbb N})$,
i\.e\.
$\{\nu\in \Cal P(\Sigma_{\smallG}^{\Bbb N})\,|\,I(\nu)<\infty\}
=
\Cal P_{S}(\Sigma_{\smallG}^{\Bbb N})$, 
whence
$\inf_{\nu\in\Cal P(\Sigma_{\smallsmallG}^{\Bbb N}),I(\nu)<\infty}(I(\nu)-F_{\varphi}(\nu))
=
\inf_{\nu\in\Cal P_{S}(\Sigma_{\smallsmallG}^{\Bbb N})}(I(\nu)-F_{\varphi}(\nu))$,
and equality  (8.14) follows immediately from this.
Also,
since the sequence 
$(\Pi_{n})_{n}
\subseteq
\Cal P\big(\,\Cal P(\Sigma_{\smallG}^{\Bbb N})\,\big)$ 
has the large deviation property with respect to
the sequence
$(n)_{n}$ and rate function 
$I$, 
we conclude from Part (2) of Theorem 8.1
that the sequence $(Q_{\varphi,n})_{n}$ has the large deviation property with 
respect to 
the sequence
$(n)_{n}$ and rate function
$(I-F_{\varphi})
-
\inf_{\nu\in\Cal P(\Sigma_{\smallsmallG}^{\Bbb N}),I(\nu)<\infty}(I(\nu)-F_{\varphi}(\nu))
=
(I-F_{\varphi})
-
\inf_{\nu\in\Cal P_{S}(\Sigma_{\smallsmallG}^{\Bbb N})}(I(\nu)-F_{\varphi}(\nu))$
(where we have used (8.14)).
As the set
$\{U\in G\}
=
U^{-1}G$
is open
and
the set $\{U\in K\}
=
U^{-1}K$
is closed,
 it therefore follows from the large deviation property that
 $$
 \aligned
 \limsup_{n}
 \frac{1}{n}
&\log Q_{\varphi,n}\Big(\{U\in G \}\Big)\\
&\ge
-
 \inf
  \Sb
	\mu\in\Cal P(\Sigma_{\smallsmallG}^{\Bbb N})\\
	{}\\
	U\mu\in G
	\endSb
 \Bigg(
 (I(\mu)-F_{\varphi}(\mu))
 -
 \inf_{\nu\in\Cal P_{S}(\Sigma_{\smallsmallG}^{\Bbb N})}
 (I(\nu)-F_{\varphi}(\nu))
 \Bigg)\,.\\
&=
-
 \inf
  \Sb
	\mu\in\Cal P_{S}(\Sigma_{\smallsmallG}^{\Bbb N})\\
	{}\\
	U\mu\in G
	\endSb
 \Bigg(
 (I(\mu)-F_{\varphi}(\mu))
 -
 \inf_{\nu\in\Cal P_{S}(\Sigma_{\smallsmallG}^{\Bbb N})}
 (I(\nu)-F_{\varphi}(\nu))
 \Bigg)\,.\\
 \limsup_{n}
 \frac{1}{n}
&\log Q_{\varphi,n}\Big(\{U\in K \}\Big)\\
&\le
-
 \inf
  \Sb
	\mu\in\Cal P(\Sigma_{\smallsmallG}^{\Bbb N})\\
	{}\\
	U\mu\in K
	\endSb
 \Bigg(
 (I(\mu)-F_{\varphi}(\mu))
 -
 \inf_{\nu\in\Cal P_{S}(\Sigma_{\smallsmallG}^{\Bbb N})}
 (I(\nu)-F_{\varphi}(\nu))
 \Bigg)\\
&=
-
 \inf
  \Sb
	\mu\in\Cal P_{S}(\Sigma_{\smallsmallG}^{\Bbb N})\\
	{}\\
	U\mu\in K
	\endSb
 \Bigg(
 (I(\mu)-F_{\varphi}(\mu))
 -
 \inf_{\nu\in\Cal P_{S}(\Sigma_{\smallsmallG}^{\Bbb N})}
 (I(\nu)-F_{\varphi}(\nu))
 \Bigg)\,. 
 \endaligned
 \tag8.15
 $$

Combining (8.12). (8.13) and (8.15) now yields
 $$
 \align
 \limsup_{n}
 \frac{1}{n}
 \log
 \sum
   \Sb
   \bold i\in\Sigma_{\smallsmallG}^{n}\\
   {}\\
   U M_{n}[\bold i]\subseteq G
	 \endSb
 s_{\bold i}
&\ge
 \log \lambda
 \,+\,
 \limsup_{n}
 \frac{1}{n}
 \log Q_{\varphi,n}\Big(\{U\in G \}\Big)\\
&\qquad\qquad
  \qquad\qquad
 \,+\,
 \limsup_{n}
 \frac{1}{n}
 \log 
 \int
 \exp
 \left(
 nF_{\varphi}
 \right)
 \,d\Pi_{n}\\
&\ge
 \log \lambda\\
&\qquad
-
 \inf
  \Sb
	\mu\in\Cal P_{S}(\Sigma_{\smallsmallG}^{\Bbb N})\\
	{}\\
	U\mu\in G
	\endSb
 \Bigg(
 (I(\mu)-F_{\varphi}(\mu))
 -
 \inf_{\nu\in\Cal P_{S}(\Sigma_{\smallsmallG}^{\Bbb N})}
 (I(\nu)-F_{\varphi}(\nu))
 \Bigg)\\
&\qquad\qquad
  \qquad\qquad
 \,-\,
 \inf_{\nu\in\Cal P_{S}(\Sigma_{\smallsmallG}^{\Bbb N})}(I(\nu)-F_{\varphi}(\nu))\\
&=
 \log \lambda
 \,+\,
 \sup
  \Sb
	\mu\in\Cal P_{S}(\Sigma_{\smallsmallG}^{\Bbb N})\\
	{}\\
	U\mu\in G
	\endSb
 (F_{\varphi}(\mu)-I(\mu))\\
&=
 \sup
  \Sb
	\mu\in\Cal P_{S}(\Sigma_{\smallsmallG}^{\Bbb N})\\
	{}\\
	U\mu\in G
	\endSb
 \left(h(\mu)+\int\varphi\,d\mu\right)\,. 
 \endalign
 $$
This completes the proof of inequality (8.10). 
Inequality (8.11) is proved similarly. 
\hfill$\square$

\bigskip

We now turn towards the second main result in this section, namely, 
Theorem 8.4 showing that
 the 
 multifractal pressure and the 
 modified
multifractal
pressure are \lq asymptotically"  the same.

\bigskip

\proclaim{Theorem 8.4}
Let $X$ be a metric space and let $U:\Cal P(\Sigma_{\smallG}^{\Bbb N})\to X$ be 
continuous with respect to the weak topology.
Let $C\subseteq X$ be a  subset of $X$ and $r>0$.
Fix a continuous function $\varphi:\Sigma_{\smallG}^{\Bbb N}\to\Bbb R$.
Then we have
 $$
 \gather
 \underline P_{C}^{U}(\varphi)
 \le
 \underline Q_{B(C,r)}^{U}(\varphi)\,,\,\,\,\,
 \underline Q_{C}^{U}(\varphi)
 \le
 \underline P_{B(C,r)}^{U}(\varphi)\,,\\
 \overline P_{C}^{U}(\varphi)
 \le
 \overline Q_{B(C,r)}^{U}(\varphi)\,,\,\,\,\,
 \overline Q_{C}^{U}(\varphi)
 \le
 \overline P_{B(C,r)}^{U}(\varphi)\,.
 \endgather
 $$
\endproclaim
\noindent{\it Proof}\newline
\noindent
This follows immediately from Lemma 6.2.
\hfill$\square$

  \bigskip


\heading
{
9. Proof of Theorem 4.4}
\endheading

The purpose of this section is to prove Theorem 4.4.

\bigskip

\noindent{\it Proof of Theorem 4.4}\newline
\noindent (1)
We must prove the following two inequalities, namely,
 $$
 \gather
  \sup
\Sb
\mu\in\Cal P_{S}(\Sigma_{\smallsmallG}^{\Bbb N})\\
{}\\
U\mu\in\overline  C
\endSb
\Bigg(
h(\mu)+\int\varphi\,d\mu
\Bigg)
\le
\inf_{r>0}
\,\,
\underline P_{B(C,r)}^{U}(\varphi)\,,
\tag9.1\\
 \inf_{r>0}
\,\,
\overline P_{B(C,r)}^{U}(\varphi)
\le
 \sup
\Sb
\mu\in\Cal P_{S}(\Sigma_{\smallsmallG}^{\Bbb N})\\
{}\\
U\mu\in\overline  C
\endSb
\Bigg(
h(\mu)+\int\varphi\,d\mu
\Bigg)\,.
\tag9.2
 \endgather
 $$

{\it Proof of (9.1).}
Since $B(C,r)$ is open with 
$\overline C\subseteq B(C,r)$, 
we conclude from Theorem 8.3  that
 $$
 \align
 \sup
\Sb
\mu\in\Cal P_{S}(\Sigma_{\smallsmallG}^{\Bbb N})\\
{}\\
U\mu\in\overline  C
\endSb
\Bigg(
h(\mu)+\int\varphi\,d\mu
\Bigg)
&\le
\sup
\Sb
\mu\in\Cal P_{S}(\Sigma_{\smallsmallG}^{\Bbb N})\\
{}\\
U\mu\in B(C,r)
\endSb
\Bigg(
h(\mu)+\int\varphi\,d\mu
\Bigg)\\
&{}\\
&\le
\underline Q_{B(C,r)}^{U}(\varphi)\,.
\tag9.3
\endalign
$$
Taking infimum over all $r>0$ in (9.3) gives
  $$
 \align
 \sup
\Sb
\mu\in\Cal P_{S}(\Sigma_{\smallsmallG}^{\Bbb N})\\
{}\\
U\mu\in \overline C
\endSb
\Bigg(
h(\mu)+\int\varphi\,d\mu
\Bigg)
&\le
\inf_{r>0}
\,\,
\underline Q_{B(C,r)}^{U}(\varphi)\,.
\tag9.4
\endalign
$$
Next, we note that it follows from Theorem 8.4
that
$\underline Q_{B(C,r)}^{U}(\varphi)
 \le
 \underline P_{B(\,B(C,r)\,,\,r\,)}^{U}(\varphi)$.
Combining this inequality with and (9.4) and using the fact that 
$B(\,B(C,r)\,,\,r\,)\subseteq B(C,2r)$, we now conclude
from Theorem 8.4  that
 $$
 \align
 \sup
\Sb
\mu\in\Cal P_{S}(\Sigma_{\smallsmallG}^{\Bbb N})\\
{}\\
U\mu\in \overline C
\endSb
\Bigg(
h(\mu)+\int\varphi\,d\mu
\Bigg)
&\le
\inf_{r>0}
\,\,
\underline Q_{B(C,r)}^{U}(\varphi)\\
&\le
\inf_{r>0}
\,\,
\underline P_{B(\,B(C,r)\,,\,r\,)}^{U}(\varphi)
\qquad\qquad
\text{[by Theorem 8.4]}\\
&\le
\inf_{r>0}
\,\,
\underline P_{B(C,2r)}^{U}(\varphi)\\
&\le
\inf_{s>0}
\,\,
\underline P_{B(C,s)}^{U}(\varphi)\,.
\endalign
$$
This completes the proof of inequality (9.1).
 
\bigskip

{\it Proof of (9.2).}
Since
$B(\,B(C,r)\,,\,r\,)
\subseteq 
B(C,2r)
\subseteq
\overline{B(C,2r)}$
and
$\overline{B(C,2r)}$ is closed,
we conclude from Theorem 8.3 and Theorem 8.4
that
  $$
  \align
  \inf_{r>0}
\,\,
\overline P_{B(C,r)}^{U}(\varphi)
&\le
  \inf_{r>0}
\,\,
\overline Q_{B(\,B(C,r)\,,\,r\,)}^{U}(\varphi)
\qquad\qquad
\qquad\qquad
\text{[by Theorem 8.4]}\\
&\le
  \inf_{r>0}
\,\,
\overline Q_{\overline{B(C,2r)}}^{U}(\varphi)\\
&\le
  \inf_{r>0}
\,\,
 \sup
\Sb
\mu\in\Cal P_{S}(\Sigma_{\smallsmallG}^{\Bbb N})\\
{}\\
U\mu\in \overline{B(C,2r)}
\endSb
\Bigg(
h(\mu)+\int\varphi\,d\mu
\Bigg)\,.
\qquad\,
\text{[by Theorem 8.3]}
\endalign
$$
Letting $U_{S}:\Cal P_{S}(\Sigma_{\smallG}^{\Bbb N})\to X$ denote the restriction of $U$ to 
$\Cal P_{S}(\Sigma_{\smallG}^{\Bbb N})$,
the above inequality can be written as
    $$
  \align
  \inf_{r>0}
\,\,
\overline P_{B(C,r)}^{U}(\varphi)
&\le
  \inf_{r>0}
\,\,
 \sup
\Sb
\mu\in U_{S}^{-1}\overline{B(C,2r)}
\endSb
\Bigg(
h(\mu)+\int\varphi\,d\mu
\Bigg)\\
&=
  \inf_{n}
\,\,
 \sup
\Sb
\mu\in U_{S}^{-1}\overline{B(C,\frac{1}{n})}
\endSb
\Bigg(
h(\mu)+\int\varphi\,d\mu
\Bigg)\,.
\tag9.5
\endalign
$$

Next, note that
since
$  \overline{B(C,\frac{1}{n})}$ is closed and $U_{S}$ is continuous,
the set
$U_{S}^{-1}\overline{B(C,\frac{1}{n})}$
is a closed 
  subset of 
 $\Cal P_{S}(\Sigma_{\smallG}^{\Bbb N})$.
 As  $\Cal P_{S}(\Sigma_{\smallG}^{\Bbb N})$ is compact, we therefore
 deduce that
  $U_{S}^{-1}\overline{B(C,\frac{1}{n})}$
  is compact.
 Also, note that it follows from [Wa, Theorem 8.2]
 that the
 entropy function $h:\Cal P_{S}(\Sigma_{\smallG}^{\Bbb N})\to\Bbb R$
 is upper semi-continuous.
 We conclude from this
 that the
 map $F:\Cal P_{S}(\Sigma_{\smallG}^{\Bbb N})\to\Bbb R$
 defined by
 $F(\mu)
 =
  h(\mu)+\int\varphi\,d\mu$ is upper semi-continuous.
  Finally, since the sets 
  $K_{n}=U_{S}^{-1}\overline{B(C,\frac{1}{n})}$
  are compact with 
  $K_{1}\supseteq K_{2}\supseteq K_{3}\supseteq\ldots$
  and
  $F$ is upper semi-continuous, we deduce from Lemma 6.4 that
   $$
  \align
  \inf_{n}
  \,\,
  \sup
\Sb
\mu\in U_{S}^{-1}\overline{B(C,\frac{1}{n})}
\endSb
\Bigg(
h(\mu)+\int\varphi\,d\mu
\Bigg)
&=
\inf_{n}
\sup
\Sb
\mu\in K_{n}
\endSb
F(\mu)\\
&=
\sup
\Sb
\mu\in \bigcap_{n}K_{n}
\endSb
F(\mu)\\
&=
\sup
\Sb
\mu\in \bigcap_{n}U_{S}^{-1}\overline{B(C,\frac{1}{n})}
\endSb
\Bigg(
h(\mu)+\int\varphi\,d\mu
\Bigg)\,.
\tag9.6
\endalign
$$

 Next, observe that
  $\bigcap_{n}U_{S}^{-1}\overline{B(C,\frac{1}{n})}
  \subseteq
  U_{S}^{-1}(\,\bigcap_{n}\overline{B(C,\frac{1}{n})})
  =
  U_{S}^{-1}\overline C$,
  whence
 $$
 \align
\sup
\Sb
\mu\in \bigcap_{n}U_{S}^{-1}\overline{B(C,\frac{1}{n})}
\endSb
\Bigg(
h(\mu)+\int\varphi\,d\mu
\Bigg)
&\le
\sup
\Sb
\mu\in \bigcap_{n}U_{S}^{-1}\overline C
\endSb
\Bigg(
h(\mu)+\int\varphi\,d\mu
\Bigg)\\
&=
\sup
\Sb
\mu\in \Cal P_{S}(\Sigma_{\smallsmallG}^{\Bbb N})\\
{}\\
U\mu\in\overline C
\endSb
\Bigg(
h(\mu)+\int\varphi\,d\mu
\Bigg)\,.
\tag9.7
\endalign
 $$

  Finally, combining (9.5), (9.6) and (9.7) gives inequality (9.2).

  \noindent
  (2)
  This part follows immediately from Part (1) and Proposition 4.2.
  \hfill$\square$

  \bigskip


\heading
{10. Proof of Theorem 4.6}
\endheading

The purpose of this section is to prove Theorem 5.5.
We first prove two small lemmas.

\bigskip

\proclaim{Lemma 10.1}
let
$\Delta:\Cal P(\Sigma_{\smallG}^{\Bbb N})\to \Bbb R$
be continuous 
with
$\Delta(\mu)\not=0$
for all $\mu\in\Cal P(\Sigma_{\smallG}^{\Bbb N})$.
The either  $\Delta<0$ or $\Delta>0$.
\endproclaim
\noindent{\it Proof}\newline
\noindent
This is clear since $\Cal P(\Sigma_{\smallG}^{\Bbb N})$
is convex and therefore, in particular, connected. 
\hfill$\square$

\bigskip

\proclaim{Lemma 10.2}
Let $X$ be a normed vector space.
Let $\Gamma:\Cal P(\Sigma_{\smallG}^{\Bbb N})\to X$
be continuous and affine
and let
$\Delta:\Cal P(\Sigma_{\smallG}^{\Bbb N})\to \Bbb R$
be continuous and affine
with
$\Delta(\mu)\not=0$
for all $\mu\in\Cal P(\Sigma_{\smallG}^{\Bbb N})$.
Define 
$U:\Cal P(\Sigma_{\smallG}^{\Bbb N})\to X$
by
$U=\frac{\Gamma}{\Delta}$.
Let $C$ be a closed and convex subset of $X$ and assume that
 $$
 \overset{\,\circ}\to{C}
 \cap
 \,
 U\big(\,\Cal P_{S}(\Sigma_{\smallG}^{\Bbb N})\,\big)
 \not=
 \varnothing\,.
 $$
Then
$$
\sup
\Sb
\mu\in \Cal P_{S}(\Sigma_{\smallsmallG}^{\Bbb N})\\
{}\\
U\mu\in C
\endSb
\Bigg(
h(\mu)+\int\varphi\,d\mu
\Bigg)
=
\sup
\Sb
\mu\in \Cal P_{S}(\Sigma_{\smallsmallG}^{\Bbb N})\\
U\mu\in \overset{\,\circ}\to{C}
\endSb
\Bigg(
h(\mu)+\int\varphi\,d\mu
\Bigg)\,.
$$
\endproclaim
\noindent{\it Proof}\newline
\noindent
For brevity define 
$F:\Cal P_{S}(\Sigma_{\smallG}^{\Bbb N})\to\Bbb R$ by
$F(\mu)
=
h(\mu)+\int\varphi\,d\mu$.
It clearly suffices to show that
$$
\sup
\Sb
\mu\in \Cal P_{S}(\Sigma_{\smallsmallG}^{\Bbb N})\\
{}\\
U\mu\in C
\endSb
F(\mu)
\le
\sup
\Sb
\mu\in \Cal P_{S}(\Sigma_{\smallsmallG}^{\Bbb N})\\
U\mu\in \overset{\,\circ}\to{C}
\endSb
F(\mu)\,.
\tag10.1
$$

We will now prove inequality (10.1).
Write
$s
 =
\sup_{
\mu\in \Cal P_{S}(\Sigma_{\smallsmallG}^{\Bbb N})\,,\,
U\mu\in C
}
F(\mu)$.
Fix $\varepsilon>0$.
It follows from the definition of $s$ that we can choose 
$\lambda\in\Cal P_{S}(\Sigma_{\smallsmallG}^{\Bbb N})$ with
$U\lambda\in C$
and
$F(\lambda)>s-\varepsilon$.
Also, since
$\overset{\,\circ}\to{C}
 \cap
 \,
 U\big(\,\Cal P_{S}(\Sigma_{\smallsmallG}^{\Bbb N})\,\big)
 \not=
 \varnothing$, we can find
 $\nu\in\Cal P_{S}(\Sigma_{\smallsmallG}^{\Bbb N})$, with
 $U\nu\in \overset{\,\circ}\to{C}$.
For $t\in(0,1)$ we now define $\gamma_{t}\in\Cal P_{S}(\Sigma_{\smallsmallG}^{\Bbb N})$
by
$\gamma_{t}=t\nu+(1-t)\lambda$.
Next, we prove the following two claims.

\smallskip

{\it Claim 1. 
For all $t\in(0,1)$, we have
$U\gamma_{t}\in \overset{\,\circ}\to{C}$.}

\noindent
{\it Proof of Claim 1.}
Fix $t\in(0,1)$.
Write
$a
=
\frac
{t\Delta(\nu)}
{t\Delta(\nu)+(1-t)\Delta(\lambda)}$
and 
$b
=
\frac
{(1-t)\Delta(\lambda)}
{t\Delta(\nu)+(1-t)\Delta(\lambda)}$.
We now make a few observations.
We first observe that
 it follows from Lemma 10.1
 that
either  $\Delta<0$ or $\Delta>0$.
This clearly implies that
 $a,b\in(0,1)$.
Next, we  note that
 $
  U\gamma_{t}
=
 \frac{\Gamma(t\nu+(1-t)\lambda)}{\Delta(t\nu+(1-t)\lambda)}
 =
 \frac
 {t\Gamma(\nu)+(1-t)\Gamma(\lambda)}
 {t\Delta(\nu)+(1-t)\Delta(\lambda)}
  =
  \frac
 {t\Gamma(\nu)}
 {t\Delta(\nu)+(1-t)\Delta(\lambda)} 
 +
  \frac
 {(1-t)\Gamma(\lambda)}
 {t\Delta(\nu)+(1-t)\Delta(\lambda)} 
 =
  \frac
 {t\Delta(\nu)}
 {t\Delta(\nu)+(1-t)\Delta(\lambda)} U\nu
 +
  \frac
 {(1-t)\Delta(\lambda)}
 {t\Delta(\nu)+(1-t)\Delta(\lambda)} U\lambda
 =
 aU\nu+bU\lambda$.
 We can now prove that
 $U\gamma_{t}\in \overset{\,\circ}\to{C}$.
 Indeed,
 since
 $a,b\in(0,1)$ with $a+b=1$
 and
 $U\lambda\in C$
and
$U\nu\in \overset{\,\circ}\to{C}$,
we conclude from [Con, p\. 102, Proposition 1.11]
 that
  $
  U\gamma_{t}
 =
 aU\nu+bU\lambda\in\overset{\,\circ}\to{C}$.
 This completes the proof of Claim 1.

\smallskip

{\it Claim 2. 
There is $\mu_{0}\in \Cal P_{S}(\Sigma_{\smallG}^{\Bbb N})$ with 
$U\mu_{0}\in \overset{\,\circ}\to{C}$
such that
$F(\mu_{0})>s-\varepsilon$.
}

\noindent
{\it Proof of Claim 2.}
Since
the entropy function $h:\Cal P_{S}(\Sigma_{\smallG}^{\Bbb N})$ is affine (see [Wa]),
we conclude that
 $F$ is affine, and so
 $F(\gamma_{t})
 =
 F(t\nu+(1-t)\lambda)
=
tF(\nu)+(1-t)F(\lambda)
\to
F(\lambda)
>
s-\varepsilon$.
This implies that there is $t_{0}\in(0,1)$
with
$F(\gamma_{t_{0}})>s-\varepsilon$.
Now put $\mu_{0}=\gamma_{t_{0}}$.
Then
$F(\mu_{0})=F(\gamma_{t_{0}})>s-\varepsilon$
and Claim 1 implies that
$U\mu_{0}=U\gamma_{t_{0}}\in\overset{\,\circ}\to{C}$.
This completes the proof of Claim 2.

\smallskip

We can now prove inequality (10.1). 
Indeed, it follows from Claim 2 that there is 
$\mu_{0}\in \Cal P_{S}(\Sigma_{\smallG}^{\Bbb N})$ with 
$U\mu_{0}\in \overset{\,\circ}\to{C}$
such that
$F(\mu_{0})>s-\varepsilon$,
whence
$s-\varepsilon
<
F(\mu_{0})
\le
\sup_{
\mu\in \Cal P_{S}(\Sigma_{\smallsmallG}^{\Bbb N})\,,\,U\mu\in \overset{\,\circ}\to{C}
}
F(\mu)$.
Finally, letting $\varepsilon\searrow 0$ gives
$s
\le
\sup_{
\mu\in \Cal P_{S}(\Sigma_{\smallsmallG}^{\Bbb N})
\,,\,
U\mu\in \overset{\,\circ}\to{C}
}
F(\mu)$.
\hfill$\square$

\bigskip

\noindent
We can now prove Theorem 4.6.

\bigskip

\noindent{\it Proof of Theorem 4.6}\newline
\noindent
In view of Lemma 10.2, 
it suffices
to prove the following two inequalities, namely,
 $$
\gather
\sup
\Sb
\mu\in \Cal P_{S}(\Sigma_{\smallsmallG}^{\Bbb N})\\
{}\\
U\mu\in \overset{\,\circ}\to{C}
\endSb
\Bigg(
h(\mu)+\int\varphi\,d\mu
\Bigg)
\le
\underline P_{C}^{U}(\varphi)\,,
\tag10.2\\
{}\\
\overline P_{C}^{U}(\varphi)
\le
\sup
\Sb
\mu\in \Cal P_{S}(\Sigma_{\smallsmallG}^{\Bbb N})\\
{}\\
U\mu\in C
\endSb
\Bigg(
h(\mu)+\int\varphi\,d\mu
\Bigg)\,.
\tag10.3
\endgather
$$

{\it Proof of inequality (10.2).}
For $r>0$, let
$G_{r}
=
\{x\in C\,|\,\dist(x,X\setminus C)>r\}$,
and
note that
$G_{r}$ is open 
with
$B(G_{r},\rho)\subseteq C$
for all $0<\rho<r$.
We therefore conclude from Theorem 8.3
and Theorem 8.4
 that if $0<\rho<r$, then
 $$
 \align
 \underline P_{C}^{U}(\varphi)
&\ge
  \underline P_{ B(G_{r},\rho)}^{U}(\varphi)\\
 &\ge
  \underline Q_{ G_{r}}^{U}(\varphi)
  \qquad\qquad
  \qquad\qquad
  \qquad\,\,\,\,\,
  \text{[by Theorem 8.4]}\\
 &\ge
\sup
\Sb
\mu\in \Cal P_{S}(\Sigma_{\smallsmallG}^{\Bbb N})\\
{}\\
U\mu\in G_{r}
\endSb
\Bigg(
h(\mu)+\int\varphi\,d\mu
\Bigg)\,.
  \qquad
  \text{[by Theorem 8.3]}
  \tag10.4
 \endalign
 $$
Taking supremum over all $r>0$ in (10.4)  yields
 $$
 \align
 \underline P_{C}^{U}(\varphi)
 &\ge
 \sup_{r>0}
 \,\,
\sup
\Sb
\mu\in \Cal P_{S}(\Sigma_{\smallsmallG}^{\Bbb N})\\
{}\\
U\mu\in G_{r}
\endSb
\Bigg(
h(\mu)+\int\varphi\,d\mu
\Bigg)\,.
 \endalign
 $$
Letting $U_{S}:\Cal P_{S}(\Sigma_{\smallG}^{\Bbb N})\to(0,\infty)$ denote the restriction of $U$ to
$\Cal P_{S}(\Sigma_{\smallG}^{\Bbb N})$,
the previous inequality can be written as
 $$
\align
 \underline P_{C}^{U}(\varphi)
 &
 \ge
 \sup_{r>0}
 \,\,
\sup
\Sb
\mu\in U_{S}^{-1} G_{r}
\endSb
\Bigg(
h(\mu)+\int\varphi\,d\mu
\Bigg)
\\
&
=
\sup
\Sb
\mu\in \bigcup_{r>0}U_{S}^{-1} G_{r}
\endSb
\Bigg(
h(\mu)+\int\varphi\,d\mu
\Bigg)\,.
\tag10.5
 \endalign
 $$
However, it is easily seen that
$\bigcup_{r>0} G_{r}
=
\overset{\,\circ}\to{C}$,
whence
$\bigcup_{r>0}U_{S}^{-1} G_{r}
=
U_{S}^{-1}(\,\bigcup_{r>0} G_{r})
=
U_{S}^{-1}\overset{\,\circ}\to{C}$.
We conclude from this inclusion and  inequality (10.5) that
  $$
 \align
 \underline P_{C}^{U}(\varphi)
 &\ge
\sup
\Sb
\mu\in U_{S}^{-1}\overset{\,\circ}\to{C}
\endSb
\Bigg(
h(\mu)+\int\varphi\,d\mu
\Bigg)\\
&=
\sup
\Sb
\mu\in \Cal P_{S}(\Sigma_{\smallsmallG}^{\Bbb N})\\
U\mu\in \overset{\,\circ}\to{C}
\endSb
\Bigg(
h(\mu)+\int\varphi\,d\mu
\Bigg)\,.
 \endalign
 $$
This proves inequality (10.2).

{\it Proof of inequality (10.3).}
Since $C$ is closed we immediately 
 conclude from Theorem 8.3
and Theorem 8.4
 that
 $$
 \align
\overline P_{C}^{U}(\varphi)
&\le
\overline Q_{C}^{U}(\varphi)
  \qquad\qquad
  \qquad\qquad
  \qquad\,\,\,\,\,\,\,\,
  \text{[by Theorem 8.4]}\\
&\le
\sup
\Sb
\mu\in \Cal P_{S}(\Sigma_{\smallsmallG}^{\Bbb N})\\
{}\\
U\mu\in C
\endSb
\Bigg(
h(\mu)+\int\varphi\,d\mu
\Bigg)\,.
  \qquad
  \text{[by Theorem 8.3]}
\endalign
$$
This proves inequality (10.3).
\hfill$\square$



\heading{Acknowledgements.}\endheading

We thank an anonymous
referee for detailed and useful
comments.



\end